\documentclass[11pt]{article}
\usepackage{amsfonts}

\usepackage{graphics}
\usepackage{indentfirst}
\usepackage{cite}
\usepackage{latexsym}
\usepackage{amsmath}
\usepackage{amssymb}
\usepackage[dvips]{epsfig}
\usepackage{amscd}

\let\oldsection\section
\renewcommand\section{\setcounter{equation}{0}\oldsection}

\allowdisplaybreaks
\def\pf{\it{Proof.}\rm\quad}

\def\R{\mathbb{R}}\def\pa{\partial}

\newcommand\divg{{\text{div}}}

\newtheorem{thm}{Theorem}[section]
\newtheorem{pro}{Proposition}[section]
\newtheorem{lem}{Lemma}[section]

\newtheorem{theorem}{Theorem}[section]

\newtheorem{lemma}[theorem]{Lemma}

\renewcommand{\div}{ {\rm div }  }

\newcommand{\na}{\nabla }

\renewcommand{\r}{\mathbb{R}}
\newcommand{\bi}{\bibitem}

\newcommand{\bl}{\begin{lemma}}
\newcommand{\el}{\end{lemma}}
\newcommand{\et}{\end{theorem}}

\newcommand{\te}{\theta}

\newcommand{\la}{\label}

\newcommand{\bn}{\begin{eqnarray}}
\newcommand{\en}{\end{eqnarray}}
\newcommand{\bnn}{\begin{eqnarray*}}
\newcommand{\enn}{\end{eqnarray*}}

\newcommand{\bnnn}{\begin{eqnarray*}}
\newcommand{\ennn}{\end{eqnarray*}}

\newcommand{\ba}{\begin{aligned}}
\newcommand{\ea}{\end{aligned}}
\newcommand{\be}{\begin{equation}}
\newcommand{\ee}{\end{equation}}

\def\norm[#1]#2{\|#2\|_{#1}}

\newcommand{\n}{\rho}
\newcommand{\si}{\sigma}

\def\la{\label}

\def\na{\nabla}

\begin{document}
\title{\bf Global Well-Posedness with Large Oscillations and Vacuum to the Three-Dimensional Equations of Compressible Nematic Liquid Crystal Flows\thanks{This work was partially
supported by
  NNSFC (Grant Nos.  10801111, 10971215 \& 10971171),
  the Fundamental Research Funds for the Central Universities (Grant No. 2010121006),
  and the Natural Science Foundation of Fujian Province of China (Grant No.
  2010J05011).
      }}
\author{{Jing Li}\\
\small Institute of Applied Mathematics, AMSS \& Hua Loo-Keng Key Laboratory of Mathematics, \\
\small Chinese Academy of Sciences,  Beijing 100190, P. R. China\\
\small E-mail: ajingli@gmail.com\\[2mm]
{Zhonghai Xu}\\
\small College of Science, Northeast Dianli University, Jilin
132013, P. R. China\\
\small E-mail: xuzhonghai@163.com\\[2mm]
{Jianwen Zhang}\\
\small School of Mathematical Sciences, Xiamen University, Xiamen
361005, P. R. China\\
\small E-mail: jwzhang@xmu.edu.cn}

\date{}

\maketitle \noindent{\bf Abstract.} This paper is concerned with the
three-dimensional equations of a simplified hydrodynamic flow
modeling the motion of compressible, nematic liquid crystal
materials. The authors  establish the global existence of
classical solution to the Cauchy problem with smooth initial data
which are of small energy but possibly large oscillations with
constant state as far-field condition which could be either vacuum
or non-vacuum. The initial density is allowed to vanish and the
spatial measure of the set of vacuum can be arbitrarily large, in
particular, the initial density can even have compact support. As a
byproduct, the large-time behavior of the solution is also studied.

\vskip 2mm


\section{Introduction}
Nematic liquid crystals are aggregates of molecules which possess
same orientational order and are made of elongated, rod-like
molecules. The continuum theory of liquid crystals was established
by Ericksen \cite{E1962} and Leslie \cite{Le1968} during the period
of 1958 through 1968, see also the book by de Gennes \cite{Ge1974}.
Since then, there have been some important research developments in
liquid crystals from both theoretical and applied aspects. When the
fluid containing nematic liquid crystal materials is at rest, we
have the well-known Ossen-Frank theory for static nematic liquid
crystals, see the pioneering work by Hardt-Lin- Kinderlehrer
\cite{HKL1986} on the analysis of energy minimal configurations of
nematic liquid crystals. In general, the motion of fluid always
takes place. The so-called Ericksen-Leslie system is a macroscopic
continuum description of the time evolution of the materials under
the influence of both the flow velocity field $u$ and the
macroscopic description of the microscopic orientation
configurations $d$ of rod-like liquid crystals.

When the fluid is an incompressible viscous fluid, Lin \cite{Li1989}
first derived a simplified Ericksen- Leslie equations modeling the
liquid crystal flows in 1989. Later, Lin and Liu
\cite{LL1995,LL1998} made some important analytic studies, such as
the existence of weak/strong solutions and the partial regularity of
suitable solutions of the simplified Ericksen-Leslie system, under
the assumption that the liquid crystal director field is of either
varying length by Leslie's terminology or variable degree of
orientation by Ericksen's terminology. Recently, Wang \cite{Wa2011} obtained the global well-posedness of the hydrodynamic flow
of nematic liquid crystals in the entire space under some small conditions. However, when the fluid is
compressible, the Ericksen-Leslie system becomes more complicate and
there are very few analytic works available yet. The readers can
refer to the recent works due to Morro \cite{Mo2009} and
Zakharov-Vakulenko \cite{ZV2009} on the modeling and numerical
studies, respectively.

In this paper, we consider a simplified version of Ericksen-Leslie
equations modeling the hydrodynamic flow of compressible, nematic
liquid crystals in the whole spatial domain $\R^3$:
\begin{eqnarray}
\rho_t+\divg (\rho u)&=&0,\label{1.1}
\\
(\rho u)_t+\divg(\rho u\otimes u)+\nabla P(\rho)&=&\mu\Delta
u+(\mu+\lambda)\nabla\divg u-\nabla d\cdot\Delta d,\label{1.2}\\
d_t+u\cdot\nabla d&=&\Delta d+|\nabla d|^2d,\label{1.3}
\end{eqnarray}
where $\rho : \R^3\times[0,\infty)\to \R^+$ is the density of the
fluid, $u :  \R^3\times[0,\infty)\to\R^3$ is the velocity field of
the fluid, $P :  \R^3\times[0,\infty)\to\R^+$ is the pressure of the
fluid, and $d :  \R^3\times[0,\infty)\to \mathbb{S}^2$ (the unit
sphere in $\R^3$) is the unit-vector field that represents the
macroscopic/continuum molecular orientations. The viscosity
coefficients $\mu,\lambda$ satisfy the physical conditions:
\begin{equation}
\mu>0,\quad \lambda+\frac{2}{3}\mu\geq 0.\label{1.4}
\end{equation}
The pressure $P(\rho)$ is usually determined through the equation of
states. Here, we focus our interest on the case of isentropic flows
and assume that
\begin{equation}
P(\rho)\triangleq A\rho^\gamma\quad{\rm with}\quad
A>0,\gamma>1.\label{1.5}
\end{equation}

Though the system (\ref{1.1})--(\ref{1.3}) is a simplified version,
it still retains most of the interesting mathematical properties
(without destroying the basic nonlinear structure) of the origin
Ericksen-Leslie model for the hydrodynamics of nematic liquid
crystals (cf. \cite{Er1961,Er1987,HK1987,Le1968,Li1989}). From the
viewpoint of partial differential equations,
(\ref{1.1})--(\ref{1.3}) is a coupled hyperbolic-parabolic system
with strong nonlinearities, and thus, its mathematical analysis is
full of challenge, especially in the case that vacuum states are
allowed (i.e. the density $\rho$ may vanish). It is worth pointing
out that when $d \equiv{\it Constant}$, (\ref{1.1})--(\ref{1.3})
reduces to the famous Navier-Stokes equations which describe the
three-dimensional motion of compressible viscous isentropic flows
and have been studied by many authors, see, for example,
\cite{Li1998,HLX2010,Fe2004,MN1980}, and the references cited
therein.

We shall look for the solutions, $(\rho(x,t),u(x,t), d(x,t))$, to
the Cauchy problem of (\ref{1.1})--(\ref{1.5}) with the far-field
behavior
\begin{equation}
(\rho,u,d)(x,t)\to (\tilde\rho,0,{\bf 1} )  \quad\mbox{ as }\;\;
|x|\to\infty,\;\; t>0,\label{1.6}
\end{equation}
and the initial data
\begin{equation}
(\rho,u,d)(x,0)=(\rho_0,u_0,d_0)(x),\quad x\in\R^3,\label{1.7}
\end{equation}
where $\tilde\rho\geq0$ is a given nonnegative constant, `` ${\bf
1}$ " is a given unit vector and $|d_0|\equiv1$.

The local-in-time strong solutions to the initial value or
initial-boundary value problem of (\ref{1.1})--(\ref{1.3}) with
nonnegative initial density were studied by Ding-Lin-Wang-Wen
\cite{DLWW2011} and Huang-Wang \cite{HWW2011} in  1-D and 3-D
spatial domain, respectively. Recently, based on the arguments in
\cite{HLX2010-1,HLX2010-2} for the compressible Navier-Stokes
equations, some interesting results on the blow-up criterion of
strong solutions of (\ref{1.1})--(\ref{1.3}) were obtained (see
\cite{HW2011,HWW2011-1,LL2011}). However, to the author's knowledge,
the global existence of strong/classical solutions with vacuum is
still open even that the initial data are suitably small in some
sense.

  Recently,  Huang-Li-Xin \cite{HLX2010}
established the global existence of classical solution with large
oscillations and vacuum to the compressible isentropic Navier-Stokes
equations.
Motivated by \cite{HLX2010}, in  this paper, we will study the global
well-posedness of strong/classical solutions to the Cauchy problem
(\ref{1.1})--(\ref{1.7}) when the initial data are sufficiently
smooth and are suitably small in some energy-norm.

Before formulating our main results, we first explain the notations
and conventions used throughout this paper. For simplicity, we set
$$
\int fdx=\int_{\R^3} fdx.
$$
For $1\leq r\leq\infty$ and $\beta>0$, the standard homogeneous and
inhomogeneous Sobolev spaces are simply denoted as follows:
$$
\begin{cases}
L^r=L^r(\R^3),\quad D^{k,r}=\left\{u\in L_{\rm loc}^2 \; |\;
\|\nabla ^k
u\|_{L^r}<\infty\right\},\quad \|u\|_{D^{k,r}}=\|\nabla^k u\|_{L^r},\\[2mm]
W^{k,r}=L^r\cap D^{k,r},\quad H^{k}=W^{k,2},\quad D^k=D^{k,2},\quad
D^1=\left\{u\in L^6\; |\; \|\nabla u\|_{L^2}<\infty\right\},\\[2mm]
\dot H^\beta=\left\{f : \R^3\to\R\; |\; \|f\|_{\dot
H^\beta}^2=\int|\xi|^{2\beta}|{\hat f}(\xi)|^2d\xi<\infty\right\},
\end{cases}
$$
where $\hat f$ denotes the Fourier transform of $f$.

For given initial data $(\rho_0,u_0,d_0)$, we define
$$
C_0\triangleq\int\left(G(\rho_0)+\frac{1}{2}\rho_0|u_0|^2+\frac{1}{2}|\nabla
d_0|^2\right)dx,
$$
where $G(\cdot)$ is the potential energy density given by
$$
G(\rho)\triangleq\rho\int_{\tilde\rho}^\rho\frac{P(s)-P(\tilde\rho)}{s^2}ds.
$$
It is clear that
$$
\begin{cases}
G(\rho)=\frac{1}{\gamma-1}P(\rho),\quad &{\rm
if}\quad\tilde\rho=0,\\[2mm]
c_1(\tilde\rho,\bar\rho)(\rho-\tilde\rho)^2\leq G(\rho)\leq
c_2(\tilde\rho,\bar\rho)(\rho-\tilde\rho)^2,\quad &{\rm
if}\quad\tilde\rho>0,\; 0\leq\rho\leq2\bar\rho,
\end{cases}
$$
for positive constants $c_1(\tilde\rho,\bar\rho),
c_2(\tilde\rho,\bar\rho)$ depending on $\tilde\rho$ and $\bar\rho$.

Our main results in this paper now can be stated as follows.
\begin{thm}\label{thm1.1} For given numbers $M_1,M_2>0$ (not
necessarily small), $\bar\rho\geq \tilde\rho+1$, $\beta\in(1/2,1]$,
and $q\in(3,6)$, assume that the initial data $(\rho_0,u_0, d_0)$
satisfy
\begin{equation}
\begin{cases}
G(\rho_0)+\rho_0|u_0|^2+|\nabla d_0|^2\in L^1,\quad
0\leq\rho_0(x)\leq\bar\rho,\\[2mm]
\left(\rho_0-\tilde\rho,P(\rho_0)-P(\tilde\rho)\right)\in H^2\cap
W^{2,q},\\[2mm]
u_0 \in \dot H^\beta\cap D^1\cap D^2,\quad \nabla d_0\in \dot
H^\beta\cap
D^1\cap D^2,\\[2mm]
\|u_0\|_{\dot H^\beta}\leq M_1,\quad \|\nabla d_0\|_{\dot
H^\beta}\leq M_2,\quad |d_0|=1,
\end{cases}\label{1.11}
\end{equation}
and that the compatibility condition
\begin{equation}
-\mu\Delta u_0-(\lambda+\mu)\nabla{\rm div}u_0 +\nabla
P(\rho_0)+\nabla d_0\cdot\Delta d_0=\rho_0^{1/2} g\label{1.12}
\end{equation}
holds for some $g\in L^2$. Then there exists a positive constant
$\varepsilon>0$, depending only on $\mu$, $\lambda$, $A$, $\gamma$,
$\tilde\rho$, $\tilde\rho$, $M_1$ and $M_2$, such that if
\begin{equation}
C_0 \leq \varepsilon,\label{1.13}
\end{equation}
the Cauchy problem (\ref{1.1})--(\ref{1.7}) has a unique global
classical solution $(\rho,u,d)$ satisfying
\begin{equation}
0\leq\rho(x,t)\leq 2\bar\rho\quad{for\;\;all}\quad
x\in\R^3,\;t\geq0,\label{1.14}
\end{equation}
and
\begin{equation}
\left\{
\begin{array}{lll}
\left(\rho-\tilde\rho,P(\rho)-P(\tilde\rho)\right)\in C([0,T];H^2\cap W^{2,q}),\\[2mm]
u\in C([0,T];D^1\cap D^2)\cap L^\infty(\tau,T; D^3\cap D^{3,q}),\\[2mm]
u_t\in L^\infty(\tau,T;D^1\cap D^2)\cap H^1(\tau,T;D^1),\\[2mm]
\nabla d\in C([0,T]; H^2)\cap L^\infty(\tau,T; H^3),\\[2mm]
\nabla d_t\in C([0,T]; L^2)\cap H^1(\tau,T; L^2)
\end{array}\right.\label{1.15}
\end{equation}
for any $0<\tau<T<\infty$. Moreover, the following large-time
behavior:
\begin{equation}\label{1.16} \lim_{t\rightarrow
\infty}\left(\|\rho(\cdot,t)-\tilde\rho\|_{L^p} +\int\rho^{1/2}|u|^4dx+\|\nabla
u(\cdot,t)\|_{L^r }+\|\nabla d(\cdot,t)\|_{W^{1,k}} \right)=0,
\end{equation}
holds for $r\in [2,6)$, $k\in(2,6)$, and
\be\la{a116}
p\in
\begin{cases}
 (\gamma,\infty)\quad &{if}\quad \tilde\rho=0,\\
 (2 ,\infty)\quad &{if} \quad \tilde\rho>0.
\end{cases}
\ee
\end{thm}

\noindent{\it \bf Remark 1.1.} The solution obtained in Theorem
\ref{thm1.1} becomes a classical one away from the initial time.
Moreover, the oscillations of $(\rho,u,\nabla d)$ could be
arbitrarily large and the interior and far field vacuum states are
allowed.

\vskip 2mm

\noindent{\it\bf Remark 1.2.} When $d\equiv {\it Constant}$, the
system (\ref{1.1})--(\ref{1.3}) reduces to the well-known
Navier-Stokes equations for compressible isentropic flows. So,
Theorem \ref{thm1.1} improves the result due to Huang-Li-Xin
\cite{HLX2010}, since the compatibility condition (\ref{1.12})
imposed on the initial data is much weaker than the one in
\cite{HLX2010} (see also \cite{CK2006}). Indeed, to prove the
existence of a classical solution of the compressible Navier-Stokes
equations, the authors in \cite{HLX2010,CK2006} had to assume that
\begin{equation}
-\mu\Delta u_0-(\mu+\lambda)\nabla \divg u_0+\nabla P(\rho_0)=\rho
g\quad{\rm with}\quad g\in D^1,\;\sqrt\rho_0 g\in L^2.\label{1.17}
\end{equation}

\noindent{\it \bf Remark 1.3.} In \cite{FK1964,Ka1984}, Fujita-Kato
and Kato proved that the incompressible Navier-Stokes system is
globally well-posed for small initial data in the homogeneous
Sobolev spaces $\dot H^{1/2}$ or in $L^3$. In our case, since the
initial energy is small, we need the boundedness assumptions on the
$H^\beta$-norm on the initial velocity which is analogous to the one
for the compressible Navier-Stokes equations \cite{HLX2010}. Note
that $\dot H^\beta\hookrightarrow L^{6/(3-2\beta)}$ and
$6/(3-2\beta)>3$ for $\beta>1/2$. Thus, compared with the results in
\cite{FK1964,Ka1984}, the conditions on the initial velocity may be
optimal under the smallness conditions on the initial energy.

\vskip 2mm

\noindent{\it \bf Remark 1.4.} Recently, Wang \cite{Wa2008,Wa2011}
proved that the heat flow of harmonic maps (\ref{1.3}) (with $u=0$)
and the incompressible liquid crystal flow are globally well-posed
provided that $\|\nabla d_0\|_{L^3}$ and $\|u_0\|_{L^3}+\|\nabla
d_0\|_{L^3}$ are sufficiently small, respectively. In view of these
results in \cite{Wa2008,Wa2011}, the conditions imposed on the
initial director field $d_0$ may also be optimal under the smallness
conditions on the initial energy.

We now comment on the analysis of this paper. For large initial data
satisfying (\ref{1.11}) and (\ref{1.12}), one can utilize the
Galerkhin approximation method to construct the local classical
solutions in a similar manner as that in \cite{HW2011} (see,
Proposition \ref{pro5.1} below). So, to extend the classical
solution globally in time, we need some global a priori estimates on
the solutions $(\rho,u,d)$ in suitable higher norms. To do so, we
notice that (\ref{1.1})--(\ref{1.3}) are indeed a coupled system
between the Navier-Stokes equations and the equations for heat flow
of harmonic maps, so that, we shall make use of some ideas developed
in \cite{Ho2002,HLX2010}. However, compared with the compressible
Navier-Stokes equations, some new difficulties arise due to the
additional presence of the liquid crystal director field $d$ and the
weaker compatibility condition (\ref{1.12}) (cf. (\ref{1.17})).
Especially, the super critical nonlinearity $|\nabla d|^2d$ in the
transported heat flow of harmonic map equation (\ref{1.3}) and the
strong coupling nonlinear term $\nabla d\cdot\Delta d$ in the
momentum equations (\ref{1.2}) will cause serious difficulties in
the proofs of the time-independent global energy estimates.

As that in \cite{HLX2010}, it turns out that the key issue in this
paper is to prove both the time-independent upper bound for the
density and the time-dependent higher norm estimates of the
solutions $(\rho,u,d)$. For this purpose, as usual we start with the
basic energy estimate (see Lemma \ref{lem3.1}). To overcome the
difficulties induced by the nonlinearities $|\nabla d|^2d$ and
$\nabla d\cdot\Delta d$, we succeed in deriving an estimate on the
spatial $L^3$-norm of the gradient of $d$ and the estimates on the
second-order spatial derivatives of $d$, which are indeed suitably
controlled by the initial energy $C_0$ and the $\dot H^\beta$-norm
of $\nabla d$ (see Lemmas \ref{lem3.2}, \ref{lem3.3}). Then, basing
on these estimates, we can utilize the techniques in
\cite{Ho2002,hof2002,HLX2010} to obtain an estimate on the spatial
$L^3$-norm of the velocity, and to carry out some careful
initial-layer analysis which is concerned with the elegant estimates
on the gradient of the velocity, the second-order and third-order
derivatives of the director field, and the material derivatives of
the velocity as well (see Lemmas \ref{lem3.4}--\ref{lem3.7}). It is
worth pointing out that similar to that for the compressible
Navier-Stokes equations, the effective viscous flux and the
vorticity (see (\ref{2.3}) for the definition) play a very
mathematically important role in the entire analysis, which are
instrumental in controlling $\|\nabla u\|_{L^p}$ ($2\leq p\leq 6$)
by the $L^2$-norm of the gradient and the material derivative of the
velocity (see Lemma \ref{lem2.2}). With these estimates, we then can
obtain the desired estimates on both $L^1(0,\min\{1,T\};L^\infty)$
and $L^{8/3}(\min\{1,T\},T;L^\infty)$ norms of the effective viscous
flux, and thus, it follows from Zlotnik's inequality (see Lemma
\ref{lem2.3}) that the density admits a uniform-in-time upper bound
which is the cornerstone for the global classical estimates of the
solutions.

The next main step is to estimate the derivatives of the solutions.
Indeed, to achieve these, we first apply the Beale-Kato-Majda type
inequality (see Lemma \ref{lem2.4}) to prove the important estimates
on the gradients of the density and velocity by solving a logarithm
Gronwall inequality in a similar manner as that in
\cite{HLX2010,HLX2010-2}. As a result, one can easily obtain the
$L^2$-estimates for the second-order derivatives of density,
pressure and velocity. However, due to the weaker compatibility
condition in (\ref{1.12}) (cf. (\ref{1.17})), the method used in
\cite{HLX2010,HLX2010-2} cannot be applied any more to obtain
further estimates. Motivated by \cite{lzz,HuLi}, instead of the
$L^2$-method, we succeed in achieving these classical estimates by
proving some desired $L^q$-estimates ($3 <q <6$) on the higher-order
time-spatial derivatives of the density and velocity, basing on some
careful initial-layer analysis (see Lemmas
\ref{lem4.4}--\ref{lem4.6}).

The rest of this paper is organized as follows. In Sect. 2, we first
collect some known inequalities and facts which will be frequently
used later. In Sect. 3, we derive the time-independent (weighted)
energy estimates of the solutions and the key pointwise upper bound
of the density, which will be used to study the large-time behavior.
In Sect. 4, we establish the time-dependent estimates on the
higher-order norms of the solutions, which are needed for the
existence of classical solutions. Finally, the proof of the main
result (i.e. Theorem 1.1) will be done in Sect. 5

\section{Auxiliary lemmas}
In this section, we recall some elementary inequalities and known
results which will be used frequently later. We start with the
well-known Sobolev inequalities (see, for example, \cite{LSU1968}).
\begin{lem}\label{lem2.1} For $2\leq p\leq6$, $1<q<\infty$ and $3<r<\infty$, there
exists a generic constant $C>0$, depending only on $q$ and $r$, such
that for $f\in H^1$ and $g\in L^q\cap D^{1,r}$, we have
\begin{eqnarray}
\|f\|_{L^p}&\leq& C\|f\|_{L^2}^{(6-p)/(2p)}\|\nabla
f\|_{L^2}^{(3p-6)/(2p)},\label{2.1}
\\[2mm]
\|g\|_{L^\infty}&\leq& C\|g\|_{L^q}^{q(r-3)/(3r+q(r-3))}\|\nabla
g\|_{L^r}^{3r/(3r+q(r-3))}.\label{2.2}
\end{eqnarray}
\end{lem}

As that for compressible Navier-Stokes equations (see, for example,
\cite{Li1998,Ho2002,HLX2010}), the connections among the effective
viscous flux, the vorticity and the physical quantities will play an
important role in the entire analysis of the present paper. So, to
be continued,  we set
\begin{equation}
F\triangleq(2\mu+\lambda)\divg u-(P(\rho)-P(\tilde\rho)),\;\;
\omega\triangleq\nabla\times u,\;\; M(d)\triangleq\nabla
d\odot\nabla d-\frac{1}{2}|\nabla d|^2\mathbb{I}_3,\label{2.3}
\end{equation}
where $F$ is the so-called effective viscous flux, $\omega$ is the
vorticity, $\mathbb{I}_3$ is the $3\times3$ unit matrix, and
$$
\nabla d\odot\nabla
d=\sum_{k=1}^3\partial_id^k\partial_jd^k.
$$
Then, it follows from Lemma \ref{lem2.1} and the standard
$L^p$-estimates of elliptic system that
\begin{lem}\label{lem2.2}Let $(\rho,u,d)$ be a smooth solution of
(\ref{1.1})--(\ref{1.7}) on $\R^3\times(0,T]$. Then there exists a
generic constant $C>0$, which may depend on $\mu$ and $\lambda$,
such that for any $p\in[2,6]$,
\begin{eqnarray}
\|\nabla F\|_{L^p}+\|\nabla\omega\|_{L^p}&\leq& C\left(\|\rho\dot
u\|_{L^p}+\||\nabla d||\nabla^2d|\|_{L^p}\right),\label{2.4}\\
\|F\|_{L^6}+\|\omega\|_{L^6}&\leq& C\left(\|\rho\dot
u\|_{L^2}+\||\nabla d||\nabla^2d|\|_{L^2}\right),\label{2.5}
\end{eqnarray}
and
\begin{equation}
\|\nabla u\|_{L^6}\leq C\left(\|\rho\dot
u\|_{L^2}+\|P(\rho)-P(\tilde\rho)\|_{L^6}+\||\nabla d||\nabla^2d|\|_{L^2}\right).\label{2.6}
\end{equation}
\end{lem}
\pf Indeed, due to (\ref{1.2}), one has
$$
\rho\dot u=\nabla F-\mu\nabla\times\omega-\divg(M(d)),
$$
and hence,
\be\la{a26}
\Delta F=\divg(\rho\dot u)+\divg\divg(M(d)),\quad \mu\Delta
\omega=\nabla\times(\rho\dot u+\divg(M(d))),
\ee
where $\dot f\triangleq f_t+u\cdot\nabla f$ denotes the material
derivative. Thus, an application of the standard $L^p$-estimate of
elliptic system leads to (\ref{2.4}), which, together with
(\ref{2.1}), gives (\ref{2.5}). Using (\ref{2.3}), (\ref{2.5}) and
the standard $L^p$-estimate, we obtain (\ref{2.6}).\hfill$\square$

\vskip 2mm

To prove the uniform-in-time upper bound of density, we need the
following Zlotnik inequality, whose proof can be found in
\cite{Zl2000}.
\begin{lem}\label{lem2.3} Assume that the function $y\in W^{1,1}(0,T)$
solves the ODE system:
$$
y'=g(y)+b'(t)\quad{on}\quad[0,T],\quad y(0)=y_0,
$$
where $b\in W^{1,1}(0,T)$ and $g\in C(\R)$. If $g(\infty)=-\infty$
and
\begin{equation}
b(t_2)-b(t_1)\leq N_0+N_1(t_2-t_1)\label{2.7}
\end{equation}
for all $0\leq t_1\leq t_2\leq T$ with some $N_0\geq0$ and
$N_1\geq0$, then one has
\begin{equation}
y(t)\leq \max\{y_0,\xi^*\}+N_0<+\infty\quad {
on}\quad[0,T],\label{2.8}
\end{equation}
where $\xi^*\in\R$ is a constant such that
\begin{equation}
g(\xi)\leq -N_1\quad for\quad \xi\geq\xi^*.\label{2.9}
\end{equation}
\end{lem}

Finally, we recall the following Beale-Kato-Majda type inequality
(cf. \cite{B1,HLX2010-2}), which is an essential tool for
the estimates of the gradients of $(\rho, u)$.
\begin{lem}\label{lem2.4} For $q\in(3,\infty)$, there exists a constant $C(q)>0$ such that for all $\nabla u\in L^2\cap
D^{1,q}$,
\begin{equation}
\|\nabla u\|_{L^\infty}\leq C\left(\|{\rm div}u\|_{L^\infty}+
\|\nabla\times u\|_{L^\infty} \right)\ln\left(e+\|\nabla^2
u\|_{L^q}\right)+C\|\nabla u\|_{L^2} +C.\label{2.10}
\end{equation}
\end{lem}


\section{Time-independent lower-order estimates}
This section is devoted to proving the time-independent (weighted)
energy estimates and the uniform upper bound of density. Assume that
$(\rho,u,d)$ is a smooth solution of (\ref{1.1})--(\ref{1.7}) on
$\R^3\times(0,T)$ with some positive time $T$. To estimate the
solutions, we set
$$ \sigma(t)\triangleq\{1,t\}
$$
and define the following functionals:
\begin{eqnarray}
A_1(T)&\triangleq&\sup_{0\leq t\leq T}\sigma\left(\|\nabla
u\|_{L^2}^2+\|\nabla^2d\|_{L^2}^2\right)\nonumber\\
&&\qquad+\int_0^T\sigma\left(\|\rho^{1/2}\dot u\|_{L^2}^2+\|\nabla
d_t\|_{L^2}^2+\|\nabla^3d\|_{L^2}^2\right)dt,\label{3.1}
\\[2mm]
A_2(T)&\triangleq&\sup_{0\leq t\leq T}\sigma^2\left(\|\rho^{1/2}\dot
u\|_{L^2}^2+\|\nabla
d_t\|_{L^2}^2+\|\nabla^3d\|_{L^2}^2\right)\nonumber\\
&&\qquad+\int_0^T\sigma^2\left(\|\nabla \dot
u\|_{L^2}^2+\|d_{tt}\|_{L^2}^2+\|\nabla^2d_t\|_{L^2}^2\right)dt,\label{3.2}\\[2mm]
A_3(T)&\triangleq& \sup_{0\leq t\leq T}\|\nabla
d\|_{L^3}^3+\int_0^T\int|\nabla d||\nabla^2d|^2dxdt\label{3.3}\\[2mm]
A_4(T)&\triangleq&\sup_{0\leq t\leq T}
\sigma^{(3-2\beta)/4}\left(\|\nabla
u\|_{L^2}^2+\|\nabla^2d\|_{L^2}^2\right)\nonumber\\
&&\qquad +\int_0^T\sigma^{(3-2\beta)/4}\left(\|\rho^{1/2}\dot
u\|_{L^2}^2+\|\nabla^3d\|_{L^2}^2\right)dt\label{3.4}
\end{eqnarray}
and
\begin{equation}
A_5(T)\triangleq\sup_{0\leq t\leq T}\int\rho |u|^3dx\label{3.5},
\end{equation}
where $\beta\in(1/2,1]$ and $\dot v=v_t+u\cdot\nabla v$ is the
material derivative.

We shall prove the following key a priori estimates on the solutions
$(\rho, u,d)$.
\begin{pro}\label{pro3.1} For given numbers $M_1,M_2>0$, $\bar\rho\geq
\tilde\rho+1$ and $\beta\in(1/2,1]$, assume that
\begin{equation}
\begin{cases}
G(\rho_0)+\rho_0|u_0|^2+|\nabla d_0|^2\in L^1,\\[2mm]
0\leq\inf\rho_0 \leq\sup\rho_0\leq\bar\rho,\quad |d_0|=1,\\[2mm]
\|u_0\|_{\dot H^\beta}\leq M_1,\quad \|\nabla d_0\|_{\dot
H^\beta}\leq M_2.
\end{cases}\label{3.6}
\end{equation}
Suppose that $(\rho,u,d)$ is a smooth solution of
(\ref{1.1})--(\ref{1.7}) on $\R^3\times(0,T]$ satisfying
\begin{equation}\label{3.7}
\begin{cases}
0\leq\inf\limits_{(x,t)\in\R^3\times[0,T]}\rho(x,t)\leq\sup\limits_{(x,t)\in\R^3\times[0,T]}\rho(x,t)\leq
2\bar\rho,\\[4mm]
A_1(T)+A_2(T)\leq 2C_0^{1/2}, \;\; A_3(T)\leq 2C_0^{\delta_0},\;\;
A_4(\sigma(T))+A_5(\sigma(T))\leq 2C_0^{\delta_0},
\end{cases}
\end{equation}
where
\begin{equation}
\delta_0\triangleq\frac{2\beta-1}{9\beta}\in(0,1/9]\quad{for}\quad
\beta\in(1/2,1].\label{3.8}
\end{equation}
Then there exists a positive constant $\varepsilon>0$, which depends
only on $\mu,\lambda,A,\gamma,\bar\rho,\tilde\rho,\beta,M_1,$ and
$M_2$, such that
\begin{equation}
\begin{cases}
0\leq\inf\limits_{(x,t)\in\R^3\times[0,T]}\rho(x,t)\leq\sup\limits_{(x,t)\in\R^3\times[0,T]}\rho(x,t)\leq
\frac{7}{4}\bar\rho,\\[4mm]
A_1(T)+A_2(T)\leq C_0^{1/2},\;\; A_3(T)\leq C_0^{\delta_0},\;\;
A_4(\sigma(T))+A_5(\sigma(T))\leq C_0^{\delta_0}, \label{3.9}
\end{cases}
\end{equation}
provided
\begin{equation}
C_0\leq\varepsilon.\label{3.10}
\end{equation}
\end{pro}

The proof of Proposition \ref{pro3.1} is based on a series of lemmas
and is postponed to the end of this section.
 For notational convenience, throughout this section we denote by
$C$ or $C_i$ ($i=1,2,\ldots$) the generic positive constants which
may depend on $\mu,\lambda,A,\gamma,\bar\rho,\tilde\rho,\beta,M_1,$
and $M_2$, but are independent of $T$. We also sometimes write
$C(\alpha)$ to emphasize the dependence on $\alpha$.

We begin with the following estimates.
\begin{lem}\label{lem3.1} Let $(\rho,u,d)$ be a smooth solution of (\ref{1.1})--(\ref{1.7}) on $\R^3\times(0,T]$ satisfying (\ref{3.7}).
Then there exists a constant $\varepsilon_1>0$, depending only on
$\mu,\lambda,A,\gamma,\bar\rho,\tilde\rho,M_1,$ and $M_2$, such that
\begin{eqnarray}
&&\sup_{0\leq t\leq T}\int \left(\rho|u|^2+G(\rho)+|\nabla
d|^2\right)dx\nonumber\\
&&\qquad+\int_0^T\left(\|\nabla
u\|_{L^2}^2+\|d_t\|_{L^2}^2+\|\nabla^2d\|_{L^2}^2\right)dt\leq
CC_0,\label{3.11}
\end{eqnarray}
and moreover,
\begin{equation}
\int_0^T\left(\|u\|_{L^6}^4+\|\nabla
u\|_{L^2}^4+\|\nabla^2d\|_{L^2}^4\right)dt\leq C
C_0^{2\delta_0},\label{3.12}
\end{equation}
provided $C_0\leq \varepsilon_1$.
\end{lem}
\pf Multiplying (\ref{1.1}) and (\ref{1.2}) by $G'(\rho)$ and $u$ in
$L^2$ respectively, integrating by parts and adding them together,
by (\ref{2.1}) we have
\begin{eqnarray} &&\frac{d}{dt}\int\left(G(\rho)+\frac{1}{2}\rho
|u|^2\right)dx+\int\left(\mu |\nabla
u|^2+(\mu+\lambda)(\divg u)^2\right)dx\nonumber\\
&&\quad=-\int (\nabla d\cdot\Delta d)\cdot u dx\leq C\|u\|_{L^6}\|\nabla d\|_{L^3}\|\nabla^2d\|_{L^2}\nonumber\\
&&\quad\leq \frac{\mu}{4}\|\nabla u\|_{L^2}^2+C\|\nabla
d\|_{L^3}^2\|\nabla^2d\|_{L^2}^2.\label{3.13}
\end{eqnarray}
Using the fact that $|d|=1$ and integrating by parts, we infer from
(\ref{1.3}) and (\ref{2.1}) that
\begin{eqnarray*}
&& \frac{d}{dt}\int |\nabla
d|^2dx+\int\left(|d_t|^2+|\nabla^2d|^2\right)dx\nonumber\\
&&\quad=\int|d_t-\Delta d|^2dx=\int\left|u\cdot\nabla
d-|\nabla d|^2d\right|^2dx\nonumber\\
&&\quad\leq C\left(\|u\|_{L^6}^2\|\nabla d\|_{L^3}^2+\|\nabla
d\|_{L^4}^4\right)\nonumber\\
&&\quad\leq C\|\nabla d\|_{L^3}^2\left(\|\nabla
u\|_{L^2}^2+\|\nabla^2 d\|_{L^2}^2\right),
\end{eqnarray*}
which, together with (\ref{3.13}) and (\ref{3.7}), gives
\begin{eqnarray}
&&\frac{d}{dt}\int \left(G(\rho)+\frac{1}{2}\rho
|u|^2+ |\nabla
d|^2\right)dx+\left(\frac{3\mu}{4}\|\nabla u\|_{L^2}^2+\|d_t\|_{L^2}^2+\|\nabla^2d\|_{L^2}^2\right)\nonumber\\
&&\quad\leq C\|\nabla d\|_{L^3}^2\left(\|\nabla
u\|_{L^2}^2+\|\nabla^2 d\|_{L^2}^2\right)\leq
C_1C_0^{2\delta_0/3}\left(\|\nabla u\|_{L^2}^2+\|\nabla^2
d\|_{L^2}^2\right).\label{3.14}
\end{eqnarray}
Thus if $C_0$ is chosen to be such that
$$
C_0\leq \varepsilon_1\triangleq  \min\left\{1,
\left(\frac{\mu}{4C_1}\right)^{3/(2\delta_0)},\left(\frac{1}{2C_1}\right)^{3/(2\delta_0)}\right\},
$$
then integrating (\ref{3.14}) in $t$ over $[0,T]$ immediately leads
to (\ref{3.11}).

To prove (\ref{3.12}), we utilize (\ref{2.1}), (\ref{3.7}) and
(\ref{3.11}) to deduce that
\begin{eqnarray*}
&&\int_0^T\left(\|u\|_{L^6}^4 +\|\nabla u\|_{L^2}^4+\|\nabla^2d\|_{L^2}^4\right)dt\\
&&\quad\leq C\int_0^{\sigma(T)} \left(\|\nabla
u\|_{L^2}^4+\|\nabla^2d\|_{L^2}^4\right) dt+\int_{\sigma(T)}^T\sigma
\left(\|\nabla
u\|_{L^2}^4+\|\nabla^2d\|_{L^2}^4\right) dt\nonumber\\
&&\quad\leq C\sup_{0\leq
t\leq\sigma(T)}\left(\sigma^{(3-2\beta)/4}\left(\|\nabla
u\|_{L^2}^2+\|\nabla^2d\|_{L^2}^2\right)\right)^2\int_0^{\sigma(T)}\sigma^{(2\beta-3)/2}dt\nonumber\\
&&\qquad+C\sup_{\sigma(T)\leq t\leq T}\sigma\left(\|\nabla
u\|_{L^2}^2+\|\nabla^2d\|_{L^2}^2\right)
\int_{\sigma(T)}^T\left(\|\nabla
u\|_{L^2}^2+\|\nabla^2d\|_{L^2}^2\right)dt\nonumber\\
&&\quad\leq CC_0^{2\delta_0}+C C_0^{3/2}\leq C C_0^{2\delta_0},
\end{eqnarray*}
since $\beta\in(1/2,1]$ and $\delta_0\in(0,1/9]$. The proof of Lemma
\ref{lem3.1} is therefore completed.\hfill$\square$

\vskip 2mm

The next lemma is concerned with the estimate of $A_3(T)$.
\begin{lem}\label{lem3.2} Let $(\rho,u,d)$ be a smooth solution of (\ref{1.1})--(\ref{1.7}) on $\R^3\times(0,T]$ satisfying (\ref{3.7}).
Then there exists a constant $\varepsilon_2>0$, depending only on
$\mu,\lambda,A,\gamma,\bar\rho,\tilde\rho,\beta,M_1$ and $M_2$, such
that
\begin{equation}
A_3(T)+\int_0^T \|\nabla d\|_{L^9}^3 dt\leq
C_0^{\delta_0},\label{3.15}
\end{equation}
provided $C_0\leq \varepsilon_2$.
\end{lem}
\pf Operating $\nabla$ to both sides of (\ref{1.3}) gives
\begin{equation}
\nabla d_t-\nabla \Delta d=\nabla (|\nabla d|^2d)-\nabla
(u\cdot\nabla d).\label{3.16}
\end{equation}
Thus, multiplying (\ref{3.16}) by $3|\nabla d|\nabla d$ and
integrating by parts over $\R^3$, we obtain
\begin{eqnarray}
&&\frac{d}{dt}\|\nabla d\|_{L^3}^3+3\int\left(|\nabla d||\nabla^2
d|^2+|\nabla d|\left|\nabla(|\nabla d|)\right|^2\right)dx\nonumber\\
&&\quad\leq C\int\left(|\nabla d|^3|\nabla^2d|+|u||\nabla
d|^2|\nabla^2d|\right)dx\nonumber\\
&&\quad\leq\frac{1}{2}\int|\nabla d||\nabla^2d|^2dx+C\left(\|\nabla
d\|_{L^5}^5+\|\nabla u\|_{L^2}^2\|\nabla
d\|_{L^{9/2}}^3\right).\label{3.17}
\end{eqnarray}

To deal with the right-hand side of (\ref{3.17}), we first use Lemma
\ref{lem2.1} to get that
\begin{equation}
\left\{
\begin{array}{l}
\|\nabla d\|_{L^5}\leq C\|\nabla d\|_{L^3}^{2/5}\|\nabla
d\|_{L^9}^{3/5},\quad \|\nabla d\|_{L^{9/2}}\leq C\|\nabla
d\|_{L^3}^{1/2}\|\nabla d\|_{L^9}^{1/2},\\[3mm]
\|\nabla d\|_{L^9}^3=\||\nabla d|^{3/2}\|_{L^6}^2\leq
C\|\nabla(|\nabla d|^{3/2})\|_{L^2}^2\leq C\||\nabla
d|^{1/2}\nabla^2d\|_{L^2}^2,\label{3.18}
\end{array}\right.
\end{equation}
so that, by virtue of (\ref{3.7}) and Cauchy-Schwarz inequality we
find
\begin{eqnarray}
{\rm R.H.S.\; of\; (3.17)}&\leq& \left(1+C\|\nabla
d\|_{L^3}^2\right)\||\nabla d|^{1/2}\nabla^2d\|_{L^2}^2+C\|\nabla
d\|_{L^3}^3\|\nabla
u\|_{L^2}^4\nonumber\\
&\leq& \left(1+C_1C_0^{2\delta_0/3}\right)\||\nabla
d|^{1/2}\nabla^2d\|_{L^2}^2+CC_0^{\delta_0}\|\nabla
u\|_{L^2}^4.\label{3.19}
\end{eqnarray}
Thus if $C_0$ is chosen to be such that
$$
C_0\leq
\varepsilon_{2,1}\triangleq\min\left\{\varepsilon_1,C_1^{-3/(2\delta_0)}
\right\},
$$
then we obtain after putting (\ref{3.19}) into (\ref{3.17}) and
integrating it over $(0,T)$ that
\begin{eqnarray}
&&\sup_{0\leq t\leq T}\|\nabla d\|_{L^3}^3+\int_0^T\left( \|\nabla
d\|_{L^9}^3+ \||\nabla
d|^{1/2}\nabla^2d\|_{L^2}^2\right)dt\nonumber\\
&&\quad\leq C\|\nabla d_0\|_{L^3}^3+CC_0^{\delta_0}\int_0^T\|\nabla
u\|_{L^2}^4dt\nonumber\\
&&\quad\leq C\|\nabla d_0\|_{L^2}^{3(2\beta-1)/(2\beta)}\|\nabla
d_0\|_{\dot H^\beta}^{3/(2\beta)}+CC_0^{3\delta_0}\nonumber\\
&&\quad\leq CC_0^{3(2\beta-1)/(4\beta)} +CC_0^{3\delta_0}\leq
C_2C_0^{3\delta_0},\label{3.20}
\end{eqnarray}
where we have used (\ref{3.12}), (\ref{3.18}),  (\ref{2.1}), and the
Sobolev embedding inequality $\dot H^\beta \hookrightarrow
L^{6/(3-2\beta)}$.

Therefore, by choosing $C_0$ sufficiently small to be such that
$$
C_0\leq
\varepsilon_{2 }\triangleq\min\left\{\varepsilon_{2,1},C_2^{-2\delta_0}
\right\},
$$
we immediately obtain the desired estimate in (\ref{3.15}) from
(\ref{3.20}).\hfill$\square$

\vskip 2mm

The following initial-layer estimate of $d$ is crucial for the
subsequent analysis.
\begin{lem}\label{lem3.3} Let $(\rho,u,d)$ be a smooth solution of (\ref{1.1})--(\ref{1.7}) on $\R^3\times(0,T]$ satisfying (\ref{3.7}).
Then there exists a positive constant $\varepsilon_3$, depending
only on $\mu,\lambda,A,\gamma,\bar\rho,\tilde\rho,\beta,M_1$ and
$M_2$, such that for any $\te\in[0,1]$,
\begin{equation}
\sup_{0\leq t\leq
T}\left(\si^{1-\te}\|\nabla^2d\|_{L^2}^2\right)+\int_0^{T}\si^{1-\te}\left(\|\nabla
d_t\|_{L^2}^2+\|\nabla^3d\|_{L^2}^2\right)dt\leq C\|\nabla
d_0\|_{\dot H^{\te}}^2,\label{3.21}
\end{equation}
provided $C_0\leq \varepsilon_3$. In particular, if
$C_0\leq\varepsilon_3$, then
\begin{equation}
\sup_{0\leq t\leq
T}\left(\si\|\nabla^2d\|_{L^2}^2\right)+\int_0^{T}\si\left(\|\nabla
d_t\|_{L^2}^2+\|\nabla^3d\|_{L^2}^2\right)dt\leq C\|\nabla
d_0\|_{L^2}^2\leq CC_0.\label{3.22}
\end{equation}
\end{lem}
\pf Let $(\rho,u,d)$ be the smooth solution of
(\ref{1.1})--(\ref{1.7}) on $\R^3\times(0,T]$. Consider the
following Cauchy problem of linear parabolic equations:
\begin{equation}
v_t-\Delta v=-u\cdot\nabla v+d\nabla d:\nabla v,\quad
v(x,0)=v_0(x),\label{3.23}
\end{equation}
where $A:B\triangleq \sum_{i,j=1}^3a_{i,j}b_{i,j}$ for
$A=(a_{i,j})_{3\times3}$ and $B=(b_{i,j})_{3\times3}$.

By (\ref{2.1}) and H\"{o}lder inequality, we easily deduce from
(\ref{3.23}) that
\begin{eqnarray*}
&&\frac{d}{dt}\|\nabla
v\|_{L^2}^2+\left(\|v_t\|_{L^2}^2+\|\nabla^2v\|_{L^2}^2\right)\nonumber\\
&&\quad=\int|v_t-\Delta v|^2dx=\int|u\cdot\nabla v-d\nabla
d:\nabla v|^2dx\nonumber\\
&&\quad \leq C\|\nabla u\|_{L^2}^2\|\nabla v\|_{L^2}\|\nabla
v\|_{L^6}+C\|\nabla d\|_{L^9}^2\|\nabla v\|_{L^2}^{4/3}\|\nabla
v\|_{L^6}^{2/3}\nonumber\\
&&\quad\leq\frac{1}{2}\|\nabla^2v\|_{L^2}^2+C\left(\|\nabla
u\|_{L^2}^4+\|\nabla d\|_{L^9}^3\right)\|\nabla v\|_{L^2}^2,
\end{eqnarray*}
which, together with Gronwall's inequality, (\ref{3.12}) and
(\ref{3.15}), gives
\begin{equation}
\sup_{0\leq t\leq T}\|\nabla
v\|_{L^2}^2+\int_0^T\left(\|v_t\|_{L^2}^2+\|\nabla^2v\|_{L^2}^2\right)dt\leq
C\|\nabla v_0\|_{L^2}^2.\label{3.25}
\end{equation}

Applying $\nabla$ to both sides of (\ref{3.23}) and taking the
$L^2$-inner product, we get
\begin{eqnarray}
&&\frac{d}{dt}\|\nabla^2v\|_{L^2}^2+\left(\|\nabla
v_t\|_{L^2}^2+\|\nabla^3v\|_{L^2}^2\right)\nonumber\\
&&\quad\leq\int\left(|\nabla (d\nabla d:\nabla v)|^2+|\nabla
(u\cdot\nabla v)|^2\right)dx\triangleq I.\label{3.26}
\end{eqnarray}
Using Lemma \ref{lem2.1}, (\ref{3.15}) and Cauchy-Schwarz
inequality, we know that
\begin{eqnarray}
I&\leq& C\int\left[\left(|\nabla d|^2+|u|^2\right)|\nabla^2v|^2+
\left(|\nabla^2 d|^2+|\nabla
d|^4+|\nabla u|^2\right)|\nabla v|^2\right] dx\nonumber\\
&\leq& C \|\nabla d\|_{L^3}^2\|\nabla^2v\|_{L^6}^2+C\|u\|_{L^6}^2\|\nabla^2v\|_{L^2}\|\nabla^2v\|_{L^6}\nonumber\\
&&+C\|\nabla v\|_{L^\infty}^2\left(\|\nabla^2 d\|_{L^2}^2+\|\nabla
d\|_{L^3}^2\|\nabla
d\|_{L^6}^2+\|\nabla u\|_{L^2}^2\right)\nonumber\\
&\leq&\left(\frac{1}{4}+C_2C_0^{2\delta_0/3}\right)\|\nabla^3v\|_{L^2}^2+C\|\nabla^2
v\|_{L^2}^2\left(\|\nabla^2 d\|_{L^2}^4+\|\nabla
u\|_{L^2}^4\right).\label{3.27}
\end{eqnarray}
Thus if $C_0$ is chosen to be such that
$$
C_0\leq
\varepsilon_{3 }\triangleq (4C_2)^{-3/(2\delta_0)} ,
$$
then substitution of (\ref{3.27}) into (\ref{3.26}) results in
\begin{equation}
\frac{d}{dt}\|\nabla^2v\|_{L^2}^2+\left(\|\nabla
v_t\|_{L^2}^2+\|\nabla^3v\|_{L^2}^2\right)\leq C\|\nabla^2
v\|_{L^2}^2\left(\|\nabla^2 d\|_{L^2}^4+\|\nabla
u\|_{L^2}^4\right),\label{3.28}
\end{equation}
which, together with  Gronwall's inequality and (\ref{3.12}),
implies that
\begin{equation}
\sup_{0\leq t\leq T}\|\nabla^2v\|_{L^2}^2+\int_0^T\left(\|\nabla
v_t\|_{L^2}^2+\|\nabla^3v\|_{L^2}^2\right)dt \leq
C\|\nabla^2v_0\|_{L^2}^2.\label{3.29}
\end{equation}

On the other hand, multiplying (\ref{3.28}) by $\si$ and integrating
it over $(0,T)$, we infer from (\ref{3.25}) that if
$C_0\leq\varepsilon_3$, then
\begin{eqnarray}
&&\sup_{0\leq t\leq
T}\left(\si\|\nabla^2v\|_{L^2}^2\right)+\int_0^T\si\left(\|\nabla
v_t\|_{L^2}^2+\|\nabla^3v\|_{L^2}^2\right)dt\nonumber\\
&&\quad \leq C\int_0^T\|\nabla^2v\|_{L^2}^2dt\leq C \|\nabla
v_0\|_{L^2}^2.\label{3.30}
\end{eqnarray}

Since the solution operator $v_{0}\mapsto v(\cdot,t)$ is linear, one
may apply the standard Riesz-Thorin interpolation argument (see
\cite{BL1976}) to (\ref{3.29}) and (\ref{3.30}) to get that  for any
$\te\in[0,1]$,
\begin{equation}
\sup_{0\leq t\leq
T}\left(t^{1-\te}\|\nabla^2v\|_{L^2}^2\right)+\int_0^Tt^{1-\te}\left(\|\nabla
v_t\|_{L^2}^2+\|\nabla^3v\|_{L^2}^2\right)dt \leq C \|\nabla
v_0\|_{\dot H^\te}^2,\label{3.31}
\end{equation}
with a uniform constant $C$ independent of $\te$. Thus, choosing
$v_0=d_0$ so that $v=d$, then (\ref{3.21}) follows from (\ref{3.31})
directly.\hfill$\square$

\vskip 2mm

By Lemmas \ref{lem3.1}--\ref{lem3.3}, we can now derive preliminary
bounds for $A_1(T)$ and $A_2(T)$.
\begin{lem}\label{lem3.4} Let $(\rho,u,d)$ be a smooth solution of (\ref{1.1})--(\ref{1.7}) on $\R^3\times(0,T]$ satisfying (\ref{3.7}).
Then there exists a constant $\varepsilon_4>0$, depending only on
$\mu,\lambda,A,\gamma,\bar\rho,\tilde\rho,\beta,M_1$ and $M_2$, such
that
\begin{equation}
A_1(T)\leq CC_0+C\int_0^T\sigma^2\left(\|\nabla
u\|_{L^4}^4+\|P(\rho)-P(\tilde\rho)\|_{L^4}^4\right)dt,\label{3.32}
\end{equation}
and
\begin{equation}
A_2(T)\leq C C_0^{1/2+2\delta_0}+ C A_1(T)+C
\int_0^T\sigma^2\|\nabla u\|_{L^4}^4dt.\label{3.33}
\end{equation}
provided $C_0\leq\varepsilon_4$.
\end{lem}
\pf Multiplying (\ref{1.2}) by $\sigma u_t$ and integrating by parts
over $\R^3$ give
\begin{eqnarray}
&&\frac{1}{2}\frac{d}{dt}\int\sigma\left(\mu|\nabla
u|^2+(\mu+\lambda)(\divg
u)^2\right)dx+ \sigma\int\rho|\dot u|^2dx \nonumber\\
&&\quad=\frac{d}{dt}\int
\sigma\left(\left(P(\rho)-P(\tilde\rho)\right)\divg u+ M(d):\nabla
u\right)dx\nonumber\\
&&\qquad + \sigma\int\rho u\cdot\nabla u\cdot\dot u dx  -
\sigma\int\left(P(\rho)-P(\tilde\rho)\right)_t\divg u dx -
\sigma\int M(d)_t:\nabla u
dx \nonumber\\
&&\qquad+ \sigma'\int\left(\frac{1}{2}\left(\mu|\nabla
u|^2+(\mu+\lambda)(\divg
u)^2\right)-\left(P(\rho)-P(\tilde\rho)\right)\divg
u- M(d):\nabla u\right) dx \nonumber\\
&&\quad \triangleq \frac{d}{dt}I_0+ \sum_{i=1}^4I_i.\label{3.34}
\end{eqnarray}

We are now in a position of estimating the terms on the right-hand
side of (\ref{3.34}). First, by H\"{o}lder and Cauchy-Schwarz
inequalities we easily see that
\begin{eqnarray}
I_1&\leq&C \sigma\|\rho^{1/2} u\|_{L^4}\|\nabla
u\|_{L^4}\|\rho^{1/2}\dot u\|_{L^2} \nonumber\\
&\leq&\frac{1}{8} \sigma\|\rho^{1/2}\dot u\|_{L^2}^2+C
\sigma^2\|\nabla u\|_{L^4}^4 +C\|\rho^{1/2}u\|_{L^4}^4.\label{3.35}
\end{eqnarray}

To deal with $I_2$, we first deduce from (\ref{1.1}) that
\begin{equation}
\left(P(\rho)-P(\tilde\rho)\right)_t+
u\cdot\nabla\left(P(\rho)-P(\tilde\rho)\right) + \gamma P(\rho)
\divg u=0,\label{3.36}
\end{equation}
which, together with the effective viscous flux $F$ in (\ref{2.3}),
yields that
\begin{eqnarray}
I_2&=&\sigma \int \gamma P(\rho) (\divg u)^2dx+\sigma\int
u\cdot\nabla \left(P(\rho)-P(\tilde\rho)\right) \divg u
dx\nonumber\\
&\leq&C \|\nabla u\|_{L^2}^2+\frac{\sigma}{2\mu+\lambda}\int
u\cdot\nabla
\left(P(\rho)-P(\tilde\rho)\right)\left(F+(P(\rho)-P(\tilde\rho))\right)
dx\nonumber\\
&\leq&C\|\nabla
u\|_{L^2}^2+C\sigma^2\|P(\rho)-P(\tilde\rho)\|_{L^4}^4+C\sigma\int
|P(\rho)-P(\tilde\rho)|\left(|\nabla u||F|+|u||\nabla
F|\right)dx\nonumber\\
&\leq&C\|\nabla
u\|_{L^2}^2+C\sigma^2\|P(\rho)-P(\tilde\rho)\|_{L^4}^4+C\sigma\|P(\rho)-P(\tilde\rho)\|_{L^3}\|\nabla
u\|_{L^2}\|\nabla F\|_{L^2}\nonumber\\
&\leq&C\|\nabla
u\|_{L^2}^2+C\sigma^2\|P(\rho)-P(\tilde\rho)\|_{L^4}^4+C\sigma
\|\nabla u\|_{L^2}\left(\|\rho^{1/2}\dot u\|_{L^2}+\||\nabla
d||\nabla^2d|\|_{L^2}\right)\nonumber\\
&\leq&\frac{1}{8}\sigma\|\rho^{1/2}\dot u\|_{L^2}^2+C\|\nabla
u\|_{L^2}^2+C\sigma\|\nabla^3d\|_{L^2}^2+C\sigma^2\|P(\rho)-P(\tilde\rho)\|_{L^4}^4,\label{3.37}
\end{eqnarray}
where we have also used (\ref{2.1}), (\ref{3.11}), (\ref{3.15}) and
Cauchy-Schwarz inequality.

It follows from (\ref{2.1}), (\ref{3.15}) and Cauchy-Schwarz
inequality that
\begin{eqnarray}
I_3&\leq& C \sigma\|\nabla d\|_{L^3}^{1/2}\|\nabla
d\|_{L^6}^{1/2}\|\nabla
d_t\|_{L^2}\|\nabla u\|_{L^4} \nonumber\\
&\leq& C \|\nabla^2d\|_{L^2}^{2}+C\sigma\|\nabla
d_t\|_{L^2}^2+C\sigma^2\|\nabla u\|_{L^4}^4,\label{3.38}
\end{eqnarray}
and finally, it is easily seen that
\begin{eqnarray}
I_4&\leq& C \sigma'\left(\|P(\rho)-P(\tilde\rho)\|_{L^2}^2+\|\nabla
u\|_{L^2}^2+\|\nabla d\|_{L^4}^4\right)\nonumber\\
&\leq& C \sigma'\left(C_0+\|\nabla
u\|_{L^2}^2+C_0^{2\delta_0/3}\|\nabla^2
d\|_{L^2}^2\right).\label{3.39}
\end{eqnarray}

Due to (\ref{3.11}), (\ref{3.15}) and (\ref{3.22}), we have
\begin{eqnarray}
I_0&\leq&\frac{\mu}{4}\sigma\|\nabla
u\|_{L^2}^2+C\sigma\left(\|P(\rho)-P(\tilde\rho)\|_{L^2}^2+\|\nabla
d\|_{L^3}^2\|\nabla^2
d\|_{L^2}^2\right)\nonumber\\
&\leq&\frac{\mu}{4}\sigma\|\nabla u\|_{L^2}^2+CC_0,\label{3.40}
\end{eqnarray}
so that, putting (\ref{3.35}) and (\ref{3.37})--(\ref{3.40}) into
(\ref{3.34}) and integrating it over $(0,T)$, we obtain
\begin{eqnarray}
&&\sup_{0\leq t\leq T}\left(\sigma\|\nabla
u\|_{L^2}^2\right)+\int_0^T\sigma\|\rho^{1/2}\dot
u\|_{L^2}^2dt\nonumber\\
&&\quad\leq CC_0+C\int_0^T\left(\|\nabla
u\|_{L^2}^2+\|\nabla^2d\|_{L^2}^2+\sigma
\|\nabla^3d\|_{L^2}^2+\sigma\|\nabla
d_t\|_{L^2}^2\right)dt\nonumber\\
&&\qquad+C\int_0^T\sigma^2\left(\|\nabla
u\|_{L^4}^4+\|P(\rho)-P(\tilde\rho)\|_{L^4}^4\right)dt+C\int_0^T\|\rho^{1/2}u\|_{L^4}^4dt\nonumber\\
&&\quad\leq CC_0+C\int_0^T\sigma^2\left(\|\nabla
u\|_{L^4}^4+\|P(\rho)-P(\tilde\rho)\|_{L^4}^4\right)dt+C\int_0^T\|\rho^{1/2}u\|_{L^4}^4dt,\label{3.41}
\end{eqnarray}
provided $C_0\leq\varepsilon_3$. Here we have also used (\ref{3.11})
and (\ref{3.22}).

To estimate the last term on the right-hand side of (\ref{3.41}), we
observe from (\ref{2.1}), (\ref{3.7}) and (\ref{3.11}) that
$$
\begin{cases}
\|\rho^{1/2}u\|_{L^4}^4\leq C  \|\rho^{1/3} u\|_{L^3}^2\|\nabla
u\|_{L^2}^2\leq C \|\nabla u\|_{L^2}^2&\quad{\rm
for}\quad 0\leq t\leq\sigma(T),
\\[2mm]
\|\rho^{1/2}u\|_{L^4}^4\leq C \|\rho^{1/2} u\|_{L^2} \|\nabla
u\|_{L^2}^{3}\leq C \|\nabla u\|_{L^2}^2&\quad{\rm
for}\quad\sigma(T)\leq t\leq T,
\end{cases}
$$
which, combined with (\ref{3.11}), gives
\begin{eqnarray*}
\int_0^T\|\rho^{1/2}u\|_{L^4}^4dt   \leq  C \int_0^{T}\|\nabla
u\|_{L^2}^2dt
 \leq  CC_0.
\end{eqnarray*}
This, together with (\ref{3.41}) and (\ref{3.22}), finishes the
proof of (\ref{3.32}).

To prove (\ref{3.33}), operating $\sigma^m\dot
u^j[\partial_t+\divg(u\cdot)]$ with $m\geq0$ to both sides of the
$j$-th equation of (\ref{1.2}) and integrating by parts over $\R^3$,
we obtain after summing up that
\begin{eqnarray}
L&\triangleq&\frac{1}{2}\frac{d}{dt}\int\sigma^m\rho|\dot
u|^2dx-\mu\int\sigma^m\dot u^j\left[\Delta u^j_t+\divg(u\Delta
u^j)\right]dx\nonumber\\
&&-(\mu+\lambda)\int\sigma^m\dot u^j\left[\partial_j\divg
u_t+\divg(u\partial_j\divg u)\right]dx\nonumber\\
&=&\frac{m}{2}\sigma^{m-1}\sigma'\int\rho|\dot
u|^2dx-\int\sigma^m\dot u^j\left[\partial_jP_t+\divg(u
\partial_j P)\right]dx\nonumber\\
&&+\int\sigma^mM(d)_t : \nabla \dot u dx+\int\sigma^m \divg M(d)
\cdot(u\cdot\nabla \dot u) dx\triangleq\sum_{i=1}^4R_i. \label{3.42}
\end{eqnarray}

After integrating by parts and using Cauchy-Schwarz inequality, we
easily see that
\begin{eqnarray*}
L_1&\triangleq&-\mu\int\sigma^m\dot u^j\left[\Delta
u^j_t+\divg(u\Delta
u^j)\right]dx\nonumber\\
&=&\mu\int \sigma^m\left(|\nabla\dot u|^2 -\partial_k\dot
u^j\partial_k\left(u\cdot\nabla u^j\right)+\partial_k\dot
u^ju^k\Delta u^j\right) dx\nonumber\\
&=&\mu \int\sigma^m \left(|\nabla\dot u|^2-\partial_k\dot
u^j\partial_ku^l\partial_l u^j+\partial_k\dot
u^j\partial_lu^l\partial_{k} u^j-\partial_k\dot
u^j\partial_lu^k\partial_{l}u^j\right) dx\nonumber\\
&\geq&\frac{7\mu}{8}\sigma^m\|\nabla\dot u\|_{L^2}^2-C
\sigma^m\|\nabla u\|_{L^4}^4,
\end{eqnarray*}
where and in what follows, we use the Einstein convention that
repeated indices denote the summation with respect to those indices.
Similarly,
\begin{eqnarray*}
L_2&\triangleq& (\mu+\lambda)\int\sigma^m\dot
u^j\left[\partial_j\divg u_t+\divg(u\partial_j\divg
u)\right]dx\nonumber\\
&\geq&(\mu+\lambda) \sigma^m\|\divg \dot
u\|_{L^2}^2-\frac{\mu}{8}\sigma^m\|\nabla \dot u\|_{L^2}^2-C
\sigma^m\|\nabla u\|_{L^4}^4.
\end{eqnarray*}
Thus, the left-hand side of (\ref{3.42})  can be bounded from below
as follows:
\begin{equation}
L\geq \frac{1}{2}\frac{d}{dt}\left(\sigma^m\|\rho^{1/2}\dot
u\|_{L^2}^2\right)+ \frac{3\mu}{4}\sigma^m\|\nabla \dot
u\|_{L^2}^2-C \sigma^m\|\nabla u\|_{L^4}^4.\label{3.43}
\end{equation}

In view of (\ref{3.36}), we have by integration by parts and
(\ref{3.7}) that
\begin{eqnarray}
R_2 &=& \int\sigma^m\left[\partial_k\dot u^ju^k\partial_jP(\rho)-
\divg \dot u\left(\divg(P(\rho) u)+(\gamma-1)P(\rho)\divg
u\right)\right]dx\nonumber\\
&=& -\int\sigma^m\left(\partial_k\dot
u^j\partial_ju^kP(\rho)+(\gamma-1)P(\rho)(\divg \dot u)
(\divg u)\right) dx\nonumber\\
&\leq&\frac{\mu}{8}\sigma^m\|\nabla \dot
u\|_{L^2}^2+C\sigma^m\|\nabla u\|_{L^2}^2.\label{3.44}
\end{eqnarray}

Using (\ref{2.1}), (\ref{2.2}) and (\ref{3.15}), we get
\begin{eqnarray}
R_3&\leq& C\sigma^m\|\nabla d\|_{L^3}\|\nabla d_t\|_{L^6}\|\nabla
\dot u\|_{L^2}\nonumber\\
&\leq& \frac{\mu}{8}\sigma^m\|\nabla \dot
u\|_{L^2}^2+C_1C_0^{2\delta_0/3}\sigma^m\|\nabla^2d_t\|_{L^2}^2,\label{3.45}
\end{eqnarray}
and similarly,
\begin{eqnarray}
R_4&\leq& C\sigma^m\|\nabla
d\|_{L^\infty}\|\nabla^2d\|_{L^2}^{1/2}\|\nabla^2d\|_{L^6}^{1/2}\|u\|_{L^6}\|\nabla
\dot
u\|_{L^2}\nonumber\\
&\leq&  \frac{\mu}{8}\sigma^m\|\nabla \dot
u\|_{L^2}^2+C\sigma^m\|\nabla^2
d\|_{L^2}^2\|\nabla^3d\|_{L^2}^2\|\nabla
u\|_{L^2}^2\nonumber\\
&\leq&  \frac{\mu}{8}\sigma^m\|\nabla \dot
u\|_{L^2}^2+C\sigma^m\left(\|\nabla
u\|_{L^2}^4+\|\nabla^2d\|_{L^2}^4\right)\|\nabla^3d\|_{L^2}^2.\label{3.46}
\end{eqnarray}

Thus, putting (\ref{3.43})--(\ref{3.46}) into (\ref{3.42}), we know
that
\begin{eqnarray}
&&\frac{d}{dt}\left(\sigma^m\|\rho^{1/2}\dot
u\|_{L^2}^2\right)+\sigma^m\|\nabla \dot
u\|_{L^2}^2-C_1C_0^{2\delta_0/3}\sigma^m\|\nabla^2
d_t\|_{L^2}^2\nonumber\\
&&\quad\leq Cm\sigma^{m-1}\sigma'\|\rho^{1/2}\dot u\|_{L^2}^2+
C\sigma^m \| \nabla u\|_{L^2}^2+C\sigma^m
\|\nabla u\|_{L^4}^4\nonumber\\
&&\qquad +C\sigma^m \left(\|\nabla
u\|_{L^2}^4+\|\nabla^2d\|_{L^2}^4\right)\|\nabla^3d\|_{L^2}^2.\label{3.47}
\end{eqnarray}

Next, we estimate $\|\nabla^2 d_t\|_{L^2}$. To do so, noticing that
$$
d_{tt}-\Delta d_t=\left(|\nabla d|^2d-u\cdot\nabla d\right)_t,
$$
we obtain from  direct computations that
\begin{eqnarray}
&&\frac{d}{dt}\int\sigma^m|\nabla
d_t|^2dx+\int\sigma^m\left(|d_{tt}|^2+|\nabla^2d_t|^2\right)dx-m\sigma^{m-1}\sigma'\int|\nabla d_t|^2 dx\nonumber\\
&&\quad\leq C\int\sigma^m |\dot u|^2|\nabla
d|^2dx+C\int\sigma^m|u|^2|\nabla u|^2|\nabla
d|^2dx+C\int \sigma^m|u|^2|\nabla d_t|^2dx\nonumber\\
&&\qquad+C\int \sigma^m|\nabla d|^2|\nabla d_t|^2 dx+C\int
\sigma^m|d_t|^2|\nabla d|^4dx\triangleq\sum_{i=1}^5J_i,\label{3.48}
\end{eqnarray}
where we have used the facts that $|d|=1$ and $u_t=\dot
u-u\cdot\nabla u$.

For the first term on the right-hand side of (\ref{3.48}), it
follows from (\ref{2.1}) and (\ref{3.15}) that
\begin{equation}
J_1\leq C\sigma^m\|\dot u\|_{L^6}^2\|\nabla d\|_{L^3}^2\leq
C_2C_0^{2\delta_0/3}\sigma^m\|\nabla \dot u\|_{L^2}^2.\label{3.49}
\end{equation}

Using (\ref{2.1}), (\ref{2.6}), (\ref{3.11}) and (\ref{3.15}), we
see that
\begin{eqnarray}
J_2&\leq&C\sigma^m\|u\|_{L^6}^2\|\nabla d\|_{L^6}^2\|\nabla
u\|_{L^6}^2\nonumber\\
&\leq&C\sigma^m\|\nabla u\|_{L^2}^2\|\nabla
d\|_{L^6}^2\left(\|\rho\dot
u\|_{L^2}^2+\|P(\rho)-P(\tilde\rho)\|_{L^6}^2+\||\nabla
d||\nabla^2d|\|_{L^2}^2\right)\nonumber\\
&\leq&C\sigma^m\|\nabla u\|_{L^2}^2\|\nabla^2 d\|_{L^2}^2\left(
\|\rho^{1/2}\dot u\|_{L^2}^2 + C_0^{1/3}+\|\nabla
d\|_{L^3}^2\|\nabla^2d\|_{L^6}^2\right)\nonumber\\
&\leq&C \sigma^m\left(\|\nabla u\|_{L^2}^4+ \|\nabla
d\|_{L^2}^{4}\right) \left(\|\rho^{1/2}\dot
u\|_{L^2}^2+\|\nabla^3d\|_{L^2}^2\right) +C\sigma^m\|\nabla
u\|_{L^2}^2\|\nabla^2 d\|_{L^2}^{2},\label{3.50}
\end{eqnarray}
and using (\ref{2.1}) and Cauchy-Schwarz inequality, we have
\begin{eqnarray*}
\sum_{i=3}^5J_i&\leq& C\sigma^m\left(\|u\|_{L^6}^2+\|\nabla
d\|_{L^6}^2\right)\|\nabla d_t\|_{L^3}^2+C\sigma^m\|d_t\|_{L^6}^2\|\nabla d\|_{L^6}^4\nonumber\\
&\leq& C\sigma^m\left(\|\nabla u\|_{L^2}^2+\|\nabla^2
d\|_{L^2}^2\right)\|\nabla d_t\|_{L^2}\|\nabla^2 d_t\|_{L^2}+C
\sigma^m\|\nabla^2d\|_{L^2}^4\|\nabla
d_t\|_{L^2}^2\nonumber\\
&\leq&\frac{1}{4}\sigma^m\|\nabla^2d_t\|_{L^2}^2+C\sigma^m\left(\|\nabla
u\|_{L^2}^4+\|\nabla^2 d\|_{L^2}^4\right)\|\nabla
d_t\|_{L^2}^2,
\end{eqnarray*}
which, combining with (\ref{3.48})--(\ref{3.50}) gives
\begin{eqnarray}
&&\frac{d}{dt}\left(\sigma^m\|\nabla
d_t\|_{L^2}^2\right)+\sigma^m\left(\|\nabla^2d_t\|_{L^2}^2+\|d_{tt}\|_{L^2}^2\right)-C_2C_0^{2\delta_0/3}\sigma^m\|\nabla
\dot
u\|_{L^2}^2\nonumber\\
&&\quad\leq Cm\sigma^{m-1}\sigma'\|\nabla
d_t\|_{L^2}^2+C\sigma^m\|\nabla u\|_{L^2}^2\|\nabla^2 d\|_{L^2}^{2}
\nonumber\\
&&\qquad+C\sigma^m\left(\|\nabla u\|_{L^2}^4+\|\nabla^2
d\|_{L^2}^4\right)\left( \|\rho^{1/2}\dot u\|_{L^2}^2 +\|\nabla
d_t\|_{L^2}^2+\|\nabla^3d\|_{L^2}^2\right).\label{3.51}
\end{eqnarray}

Thus if $C_0$ is chosen to be such that
$$
C_0\leq\varepsilon_{4}\triangleq\min\left\{\varepsilon_3,(2C_1)^{-3/(2\delta_0)},(2C_2)^{-3/(2\delta_0)}\right\},
$$
then, choosing $m=2$ in (\ref{3.47}) and (\ref{3.51}), adding them
together, and integrating the resulting inequality over $(0,T)$, by
(\ref{3.7}) and (\ref{3.11}) we obtain
\begin{eqnarray}
&&\sup_{0\leq t\leq T}\sigma^2\left(\|\rho^{1/2}\dot
u\|_{L^2}^2+\|\nabla
d_t\|_{L^2}^2\right)\nonumber\\
&&\qquad+\int_0^T\sigma^2\left(\|\nabla \dot u\|_{L^2}^2+\|d_{tt}\|_{L^2}^2+\|\nabla^2d_t\|_{L^2}^2\right)dt\nonumber\\
&&\quad\leq C\int_0^T\sigma^2\|\nabla
u\|_{L^2}^2dt+C\int_0^T\sigma^2\|\nabla u\|_{L^2}^2\|\nabla^2
d\|_{L^2}^2dt\nonumber\\
&&\qquad+C \int_0^T\sigma^2\left(\|\nabla u\|_{L^2}^4+\|\nabla^2
d\|_{L^2}^4\right)\left( \|\rho^{1/2}\dot u\|_{L^2}^2 +\|\nabla
d_t\|_{L^2}^2+\|\nabla^3d\|_{L^2}^2\right)dt\nonumber\\
&&\qquad+ C \int_0^T\sigma\left(\|\rho^{1/2}\dot
u\|_{L^2}^2+\|\nabla d_t\|_{L^2}^2\right)dt+C
\int_0^T\sigma^2\|\nabla
u\|_{L^4}^4dt\nonumber\\
&&\quad\leq C C_0^{1/2+2\delta_0} + C A_1(T)+C
\int_0^T\sigma^2\|\nabla u\|_{L^4}^4dt,\label{3.52}
\end{eqnarray}
since it follows from (\ref{3.7}) and (\ref{3.12}) that
\begin{eqnarray*}
&&\int_0^T\sigma^2\left(\|\nabla u\|_{L^2}^4+\|\nabla^2
d\|_{L^2}^4\right)\left( \|\rho^{1/2}\dot u\|_{L^2}^2 +\|\nabla
d_t\|_{L^2}^2+\|\nabla^3d\|_{L^2}^2\right)dt\\
&&\quad\leq\sup_{0\leq t\leq T}\sigma^2\left( \|\rho^{1/2}\dot
u\|_{L^2}^2 +\|\nabla
d_t\|_{L^2}^2+\|\nabla^3d\|_{L^2}^2\right)\int_0^T\left(\|\nabla
u\|_{L^2}^4+\|\nabla^2 d\|_{L^2}^4\right)dt\\
&&\quad\leq CC_0^{2\delta_0}A_2(T)\leq CC_0^{1/2+2\delta_0}.
\end{eqnarray*}

Moreover, it follows from (\ref{3.16}) and the standard
$L^2$-estimate of elliptic system that
\begin{eqnarray*}
\|\nabla^3d\|_{L^2}&\leq& C\left( \|\nabla d_t\|_{L^2}+\|\nabla
d\|_{L^6}^3+\|\nabla d\nabla^2d\|_{L^2}\right)\nonumber\\
&&+C\left(\|\nabla u\nabla
d\|_{L^2}+\|u\nabla^2d\|_{L^2}\right) \nonumber\\
&\leq& C \left(\|\nabla d_t\|_{L^2}+\|\nabla^2 d\|_{L^2}^3+
\|\nabla^2d\|_{L^2}^{3/2}\|\nabla^3 d\|_{L^2}^{1/2}\right) \nonumber\\
&&+C \|\nabla u\|_{L^2}\|\nabla^2
d\|_{L^2}^{1/2}\|\nabla^3d\|_{L^2}^{1/2},
\end{eqnarray*}
so that,
\begin{equation}
\|\nabla^3d\|_{L^2}^2\leq C\|\nabla^2 d\|_{L^2}^2\left(\|\nabla^2
d\|_{L^2}^4+\|\nabla u\|_{L^2}^4\right)+C\|\nabla
d_t\|_{L^2}^2.\label{3.53}
\end{equation}
Therefore, keeping in mind that $\beta>1/2$, we infer from
(\ref{3.7}) and (\ref{3.52}) that
\begin{eqnarray*}
&&\sup_{0\leq t\leq T} \left(\sigma^2 \|\nabla^3d\|_{L^2}^2\right)\\
&&\quad\leq C \sup_{0\leq t\leq T} \left(\sigma^2 \|\nabla
d_t\|_{L^2}^2\right)+C \sup_{0\leq t\leq T}\left(
\sigma^2\|\nabla^2d\|_{L^2}^2\left(\|\nabla u\|_{L^2}^4+
\|\nabla^2d\|_{L^2}^4 \right)\right)\\
&&\quad\leq C C_0^{1/2+2\delta_0} + C A_1(T)+C
\int_0^T\sigma^2\|\nabla u\|_{L^4}^4dt\\
&&\qquad+C \sup_{0\leq t\leq \sigma(T)}\left(
\sigma\|\nabla^2d\|_{L^2}^2\right)\sup_{0\leq t\leq
\sigma(T)}\left(\sigma^{(3-2\beta)/4}\left(\|\nabla u\|_{L^2}^2+
\|\nabla^2d\|_{L^2}^2 \right)\right)^2\\
&&\qquad +C \sup_{\sigma(T)\leq t\leq T}\left(
\sigma\|\nabla^2d\|_{L^2}^2\right)\sup_{\sigma(T)\leq t\leq
T}\left(\sigma \left(\|\nabla u\|_{L^2}^2+ \|\nabla^2d\|_{L^2}^2
\right)\right)^2\\
&&\quad\leq C C_0^{1/2+2\delta_0} + C A_1(T)+C
\int_0^T\sigma^2\|\nabla u\|_{L^4}^4dt,
\end{eqnarray*}
which, together with (\ref{3.52}), finishes the proof of
(\ref{3.33}).\hfill$\square$

\vskip 2mm

The next lemma plays an important role in the proof of the uniform
upper bound of $\rho$. The ideas of the proof are motivated by the
ones in \cite{hof2002,HLX2010}.

\begin{lem}\label{lem3.5} Let $(\rho,u,d)$ be a smooth solution of (\ref{1.1})--(\ref{1.7}) on $\R^3\times(0,T]$
satisfying (\ref{3.7}). Then there exist  positive constants
$\varepsilon_5 $ and $C$ depending on
$\mu,\lambda,A,\gamma,\bar\rho,\tilde\rho,\beta,M_1$ and $M_2$, such
that
\begin{equation} \sup_{0\leq t\leq
\sigma(T)}\left(t^{1-\beta}\|\nabla
u\|_{L^2}^2\right)+\int_0^{\sigma(T)}t^{1-\beta}\|\rho^{1/2}\dot
u\|_{L^2}^2dt\leq C(M_1,M_2),\label{3.54}
\end{equation}
and
\begin{eqnarray}
&&\sup_{0\leq t\leq \sigma(T)}t^{2-\beta}\left(\|\rho^{1/2}\dot
u\|_{L^2}^2+\|\nabla
d_t\|_{L^2}^2\right)\nonumber\\
&&\qquad+\int_0^{\sigma(T)}t^{2-\beta}\left(\|\nabla \dot
u\|_{L^2}^2+\|d_{tt}\|_{L^2}^2+\|\nabla^2d_t\|_{L^2}^2\right)dt\leq
C (M_1,M_2),\label{3.55}
\end{eqnarray}
 provided $C_0\leq \varepsilon_5$.
\end{lem}
\pf The proof of (\ref{3.54}) is motivated by \cite{HLX2010,hof2002}.
For a fixed smooth solution $(\rho,u,d)$, we define the linear
differential operator $\mathcal{L}$ acting on the functions $w$ :
$\R^3\times[0,\infty)\to \R^3$ by
\begin{eqnarray}
(\mathcal{L}w)^j&\triangleq& \rho w^j_t+\rho u\cdot\nabla
w^j-\left(\mu\Delta w^j+(\mu+\lambda)\divg w_{x_j}\right)\nonumber\\
&=&\rho\dot w^j-\left(\mu\Delta w^j+(\mu+\lambda)\divg
w_{x_j}\right),\quad j=1,2,3,\label{3.56}
\end{eqnarray}
where $\dot w\triangleq w_t+u\cdot\nabla w.$
We  thus define $w_1$, $w_2$ and $w_3$ by
\begin{eqnarray}
&&\mathcal{L}w_1=0,\quad w_1(x,0)=w_{10}(x),\label{3.57}
\\[2mm]
&&\mathcal{L}w_2=-\nabla P(\rho),\quad w_2(x,0)=0,\label{3.58}
\end{eqnarray}
and
\begin{equation}
\mathcal{L}w_3=-\divg(\nabla d\odot\nabla v-\frac{1}{2}\nabla
d:\nabla v\mathbb{I}_3),\quad w_3(x,0)=0,\label{3.59}
\end{equation}
where $v=v(x,t)$ is the solution of (\ref{3.23}).

The estimates of $w_1,w_2$ are similar to those in
\cite{hof2002,HLX2010}. For the reader's convenience, we reproduce
the proofs here. Multiplying (\ref{3.57}) and (\ref{3.58}) by $w_1$
and $w_2$, respectively, and integrating them by parts over
$\R^3\times(0,\sigma(T))$, one easily deduces from (\ref{3.7}) and
(\ref{3.11}) that
\begin{equation}
\sup_{0\leq t\leq
\sigma(T)}\int\rho|w_1|^2dx+\int_0^{\sigma(T)}\|\nabla
w_1\|_{L^2}^2dt\leq C \|w_{10}\|_{L^2}^2,\label{3.60}
\end{equation}
and
\begin{equation}
\sup_{0\leq t\leq
\sigma(T)}\int\rho|w_2|^2dx+\int_0^{\sigma(T)}\|\nabla
w_2\|_{L^2}^2dt\leq C C_0.\label{3.61}
\end{equation}

To estimate $w_3$, multiplying (\ref{3.59}) by $w_3$ and integrating
it by parts over $\R^3\times(0,\sigma(T))$, by (\ref{2.1}) and
(\ref{3.15}) we get
\begin{eqnarray*}
&&\sup_{0\leq t\leq
\sigma(T)}\int\rho|w_3|^2dx+\int_0^{\sigma(T)}\|\nabla
w_3\|_{L^2}^2dt\nonumber\\
&&\quad\leq C\int_0^{\sigma(T)}\|\nabla d\|_{L^3}\|\nabla
v\|_{L^6}\|\nabla w_3\|_{L^2}dt\nonumber\\
&&\quad\leq \frac{1}{2}\int_0^{\sigma(T)}\|\nabla
w_3\|_{L^2}^2dt+C\int_0^{\sigma(T)}\|\nabla^2 v\|_{L^2}^2 dt,
\end{eqnarray*}
which, combined with (\ref{3.25}), yields
\begin{equation}
\sup_{0\leq t\leq
\sigma(T)}\int\rho|w_3|^2dx+\int_0^{\sigma(T)}\|\nabla
w_3\|_{L^2}^2dt\leq C \|\nabla v_0\|_{L^2}^2.\label{3.62}
\end{equation}

Next we estimate $\|\nabla w_i\|_{L^2}$ ($i=1,2,3$). Thanks to
(\ref{3.7}), an application of the standard $L^2$-estimate of
elliptic system to (\ref{3.57}) shows
$$
\|\nabla w_1\|_{L^6}\leq C\|\nabla ^2w_1\|_{L^2}\leq
C\|\rho^{1/2}\dot w_1\|_{L^2}.
$$
So, multiplying (\ref{3.57}) by $w_{1t}$ in $L^2$ and integrating by
parts, we infer from (\ref{3.7}) that
\begin{eqnarray}
&&\frac{1}{2}\frac{d}{dt}\int\left(\mu|\nabla
w_1|^2+(\mu+\lambda)(\divg w_1)^2\right)dx+\int\rho|\dot
w_1|^2dx\nonumber\\
&&\quad=\int\rho \dot w_1\cdot(u\cdot\nabla w_1)dx\nonumber\\
&&\quad \leq C \left(\int\rho|\dot
w_1|^2dx\right)^{1/2}\left(\int\rho |u|^3dx\right)^{1/3}\|\nabla
w_1\|_{L^6}\nonumber\\
&&\quad\leq C \left(\int\rho
|u|^3dx\right)^{1/3}\int\rho|\dot w_1|^2dx\nonumber\\
&&\quad\leq C_1C_0^{\delta_0/3}\int\rho|\dot w_1|^2dx\label{3.63}
\end{eqnarray}
for any $ t\in(0,\sigma(T))$. Thus if $C_0$ is chosen to be such
that
$$
C_0\leq\varepsilon_{5,1}\triangleq
\min\left\{\varepsilon_4,(2C_1)^{-3/\delta_0}\right\},
$$
then it follows from (\ref{3.63}) that
\begin{equation}
\sup_{0\leq t\leq \sigma(T)}\|\nabla
w_1\|_{L^2}^2+\int_0^{\sigma(T)}\|\rho^{1/2}\dot w_1\|_{L^2}^2dt\leq
C\|\nabla w_{10}\|_{L^2}^2.\label{3.64}
\end{equation}

On the other hand, multiplying (\ref{3.63}) by $t$ and integrating
it over $(0,\sigma(T))$, by (\ref{3.60}) we see that if $C_0\leq
\varepsilon_{5,1}$, then
\begin{eqnarray}
&&\sup_{0\leq t\leq \sigma(T)}\left(t\|\nabla
w_1\|_{L^2}^2\right)+\int_0^{\sigma(T)}t\int\rho|\dot
w_1|^2dxdt\nonumber\\
&&\qquad\leq C\int_0^{\sigma(T)}\|\nabla w_{1}\|_{L^2}^2dt\leq
C\|w_{10}\|_{L^2}^2.\label{3.65}
\end{eqnarray}

Since the solution operator $w_{10}\mapsto w_1(\cdot,t)$ is linear,
by the standard Riesz-Thorin interpolation argument (see
\cite{BL1976}), we conclude from (\ref{3.64}) and (\ref{3.65})  that
for $\beta\in (1/2,1]$,
\begin{equation}
\sup_{0\leq t\leq \sigma(T)}\left(t^{1-\beta}\|\nabla
w_1\|_{L^2}^2\right)+\int_0^{\sigma(T)}t^{1-\beta}\|\rho^{1/2}\dot
w_1\|_{L^2}^2dt\leq C \|w_{10}\|_{\dot H^\beta}^2.\label{3.66}
\end{equation}

In order to estimate $\|\nabla w_2\|_{L^2}$, analogous to the proof
of Lemma \ref{lem2.2}, we set
$$
 \tilde F\triangleq(2\mu+\lambda)\divg
w_2-(P(\rho)-P(\tilde\rho)).
$$
Then it follows from (\ref{3.58}) that
\begin{equation}
\|\nabla \tilde F\|_{L^2}\leq C\|\rho\dot w_2\|_{L^2},\quad \|\nabla
w_2\|_{L^6}\leq C\left(\|\rho\dot
w_2\|_{L^2}+\|P(\rho)-P(\tilde\rho)\|_{L^6}\right).\label{3.67}
\end{equation}
Now, multiplying (\ref{3.60}) by $w_{2t}$ in $L^2$ and integrating
by parts, by (\ref{3.36}) we obtain
\begin{eqnarray}
&&\frac{1}{2}\frac{d}{dt}\int\left(\mu|\nabla
w_2|^2+(\mu+\lambda)(\divg w_2)^2\right)dx+\int\rho|\dot
w_2|^2dx\nonumber\\
&&\quad=\frac{d}{dt}\int\left(P(\rho)-P(\tilde\rho)\right)\divg
w_{2}dx-\int\left(P(\rho)-P(\tilde\rho)\right)u\cdot\nabla\divg
w_{2}dx\nonumber\\
&&\qquad+\int\left((\gamma-1)P(\rho)+P(\tilde\rho)\right)(\divg
u)(\divg w_{2})dx +\int\rho\dot w_2\cdot\left(u\cdot\nabla
w_2\right)dx\nonumber\\
&&\quad\triangleq \frac{d}{dt}I_0+\sum_{i=1}^3I_i. \label{3.68}
\end{eqnarray}

Since $(2\mu+\lambda)\divg w_2=\tilde F+(P(\rho)-P(\tilde\rho))$,
integrating by parts and using (\ref{2.1}), (\ref{3.7}),
(\ref{3.11}) and (\ref{3.67}), we find that
\begin{eqnarray}
I_1+I_2 &\leq& C\left(\|P(\rho)-P(\tilde\rho)\|_{L^2}^2+\|\nabla
u\|_{L^2}^2+\|\nabla
w_2\|_{L^2}^2\right)\nonumber\\
&&+C\|P(\rho)-P(\tilde\rho)\|_{L^3}\|u\|_{L^6}\|\nabla\tilde F\|_{L^2} \nonumber\\
&\leq&\frac{1}{4}\|\rho^{1/2}\dot w_2\|_{L^2}^2+C\left(\|\nabla
w_2\|_{L^2}^2+ \|\nabla u\|_{L^2}^2+ C_0 \right),\label{3.69}
\end{eqnarray}
and using (\ref{3.7}), (\ref{3.11}) and (\ref{3.67}), one easily
gets
\begin{eqnarray}
I_3&\leq& C \left(\int\rho|\dot w_2|^2dx\right)^{1/2}\left(\int\rho
|u|^3dx\right)^{1/3}\|\nabla
w_2\|_{L^6}\nonumber\\
&\leq& C C_0^{\delta_0/3}\|\rho^{1/2}\dot
w_2\|_{L^2}\left(\|\rho^{1/2}\dot
w_2\|_{L^2}+\|P(\rho)-P(\tilde\rho)\|_{L^6}\right)\nonumber\\
&\leq& C_2 C_0^{\delta_0/3} \|\rho^{1/2}\dot w_2\|_{L^2}^2+C
C_0^{1/3}.\label{3.70}
\end{eqnarray}

Thus if $C_0$ is chosen to be such that
$$
C_0\leq\varepsilon_{5,2}\triangleq\min\left\{\varepsilon_{5,1},(4C_2
)^{-3/\delta_0}\right\},
$$
then we obtain after putting (\ref{3.69}), (\ref{3.70}) into
(\ref{3.68}) and integrating it over $(0,\sigma(T))$ that
\begin{eqnarray}
&&\sup_{0\leq t\leq\sigma(T)} \|\nabla w_2\|_{L^2}^2
+\int_0^{\sigma(T)} \|\rho^{1/2}\dot
w_2\|_{L^2}^2dt\nonumber\\
&&\quad \leq CC_0^{1/3}+C\int_0^{\sigma(T)}\left(\|\nabla
w_2\|_{L^2}^2+ \|\nabla u\|_{L^2}^2\right)dt\leq CC_0^{1/3},
\label{3.71}
\end{eqnarray}
where we have also used (\ref{3.11}), (\ref{3.61}) and the following
simple fact (due to (\ref{3.11})) :
$$
|I_0|\leq \frac{\mu}{4}\|\nabla w_2\|_{L^2}^2+C
\|P(\rho)-P(\tilde\rho)\|_{L^2}^2\leq\frac{\mu}{4}\|\nabla
w_2\|_{L^2}^2+C C_0.
$$

We now estimate $\|\nabla w_3\|_{L^2}$. To do so, using (\ref{2.1}),
(\ref{3.7}), (\ref{3.15}) and the standard $H^2$-regularity of
elliptic system, we first deduce from (\ref{3.59}) that
\begin{eqnarray}
&&\|\nabla w_3\|_{L^6}+\|\nabla^2w_3\|_{L^2}\nonumber\\
&&\quad\leq C\left(\|\rho\dot
w_3\|_{L^2}+\|\nabla d\nabla^2 v\|_{L^2}+\|\nabla v\nabla^2 d\|_{L^2} \right)\nonumber\\
&&\quad\leq C\left(\|\rho^{1/2}\dot w_3\|_{L^2}
+C_0^{\delta_0/3}\|\nabla^3v\|_{L^2}+\|\nabla^2v\|_{L^2}^{1/2}\|\nabla^3v\|_{L^2}^{1/2}\|\nabla^2d\|_{L^2}\right).\label{3.72}
\end{eqnarray}
Then, multiplying (\ref{3.59}) by $w_{3t}$ in $L^2$ and integrating
by parts, we find
\begin{eqnarray}
&&\frac{1}{2}\frac{d}{dt}\int\left(\mu|\nabla
w_3|^2+(\mu+\lambda)(\divg w_3)^2\right)dx+\int\rho|\dot w_3|^2dx\nonumber\\
&&\quad=\frac{d}{dt}\int \left(\nabla d\odot\nabla
v-\frac{1}{2}\nabla d:\nabla v\mathbb{I}_3\right):\nabla
w_{3}dx\nonumber\\
&&\qquad-\int \left(\nabla d\odot\nabla v_t-\frac{1}{2}\nabla
d:\nabla v_t\mathbb{I}_3\right):\nabla
w_{3}dx\nonumber\\
&&\qquad -\int \left(\nabla d_t\odot\nabla v-\frac{1}{2}\nabla
d_t:\nabla v\mathbb{I}_3\right):\nabla w_{3}dx+\int\rho\dot
w_3\cdot(u\cdot\nabla
w_3)dx\nonumber\\
&&\quad=\frac{d}{dt}J_0 +\sum_{i=1}^3J_i.\label{3.73}
\end{eqnarray}

It readily follows from (\ref{2.1}) and (\ref{3.15}) that
\begin{equation}
|J_0|\leq C\|\nabla d\|_{L^3}\|\nabla v\|_{L^6}\|\nabla w_3\|_{L^2}
\leq\frac{\mu}{4}\|\nabla
w_3\|_{L^2}^2+CC_0^{2\delta_0/3}\|\nabla^2v\|_{L^2}^2.\label{3.74}
\end{equation}
By (\ref{3.15}), (\ref{3.72}) and Cauchy-Schwarz inequality, we
have
\begin{eqnarray}
J_1&\leq&C\|\nabla
d\|_{L^3}\|\nabla v_t\|_{L^2}\|\nabla w_3\|_{L^6}\nonumber\\
&\leq& CC_0^{\delta_0/3}\|\nabla v_t\|_{L^2}\left(\|\rho^{1/2}\dot
w_3\|_{L^2}
+C_0^{\delta_0/3}\|\nabla^3v\|_{L^2}+\|\nabla^2v\|_{L^2}^{1/2}\|\nabla^3v\|_{L^2}^{1/2}\|\nabla^2d\|_{L^2}\right)\nonumber\\
&\leq& \frac{1}{4}\|\rho^{1/2} \dot w_3\|_{L^2}^2 +C  \left(\|\nabla
v_t\|_{L^2}^2+\|\nabla^3
v\|_{L^2}^2+\|\nabla^2v\|_{L^2}^2\|\nabla^2d\|_{L^2}^4\right),\label{3.75}
\end{eqnarray}
and similarly, after integrating by parts and using Cauchy-Schwarz
inequality, we have
\begin{eqnarray}
J_2&\leq& C\int\left(|d_t||\nabla^2v||\nabla w_3|+|d_t||\nabla
v||\nabla^2w_3|\right)dx\nonumber\\
&\leq&C\left(\|d_t\|_{L^2}\|\nabla^2v\|_{L^3}\|\nabla
w_3\|_{L^6}+\|d_t\|_{L^2}\|\nabla v\|_{L^\infty}\|\nabla^2
w_3\|_{L^2}\right)\nonumber\\
&\leq&C \|d_t\|_{L^2}\|\nabla^2v\|_{L^2}^{1/2}\|\nabla^3
v\|_{L^2}^{1/2} \|\nabla^2 w_3\|_{L^2}\nonumber\\
&\leq&C \|d_t\|_{L^2}\|\nabla^2v\|_{L^2}^{1/2}\|\nabla^3
v\|_{L^2}^{1/2}\nonumber\\
&& \times\left(\|\rho^{1/2}\dot w_3\|_{L^2}
+C_0^{\delta_0/3}\|\nabla^3v\|_{L^2}+\|\nabla^2v\|_{L^2}^{1/2}\|\nabla^3v\|_{L^2}^{1/2}\|\nabla^2d\|_{L^2}\right)\nonumber\\
&\leq&\frac{1}{4}\|\rho^{1/2}\dot
w_3\|_{L^2}^2+C\left(\|\nabla^3v\|_{L^2}^2+\|\nabla^2v\|_{L^2}^2\|\nabla^2d\|_{L^2}^4+\|\nabla^2v\|_{L^2}^2\|d_t\|_{L^2}^4\right).\label{3.76}
\end{eqnarray}
Finally, similar to that in (\ref{3.73}), we infer from (\ref{3.7}),
(\ref{3.15}) and (\ref{3.72}) that
\begin{eqnarray}
J_3&\leq& C \left(\int\rho|\dot w_3|^2dx\right)^{1/2}\left(\int\rho
|u|^3dx\right)^{1/3}\|\nabla
w_3\|_{L^6}\nonumber\\
&\leq& C C_0^{\delta_0/3}\|\rho^{1/2}\dot
w_3\|_{L^2}\left(\|\rho^{1/2}\dot w_3\|_{L^2}
+C_0^{\delta_0/3}\|\nabla^3v\|_{L^2}+\|\nabla^2v\|_{L^2}^{1/2}\|\nabla^3v\|_{L^2}^{1/2}\|\nabla^2d\|_{L^2}\right)\nonumber\\
&\leq& C_3 C_0^{\delta_0/3} \|\rho^{1/2}\dot w_3\|_{L^2}^2+C\left(
\|\nabla^3v\|_{L^2}^2+\|\nabla^2v\|_{L^2}^2\|\nabla^2
d\|_{L^2}^4\right).\label{3.77}
\end{eqnarray}

Thus if  $C_0$ is chosen to be such that
$$
C_0\leq \varepsilon_{5}\triangleq
\min\left\{\varepsilon_{5,2},(4C_3)^{-3/\delta_0}\right\},
$$
then  by (\ref{3.12}) and (\ref{3.29}) we deduce after putting
(\ref{3.74})--(\ref{3.77}) into (\ref{3.73}) that
\begin{eqnarray}
&&\sup_{0\leq t\leq \sigma(T)}\|\nabla w_3\|_{L^2}^2 +\int_0^{\sigma(T)}\|\rho^{1/2}\dot w_3\|_{L^2}^2dt\nonumber\\
&&\quad\leq C\sup_{0\leq t\leq \sigma(T)}\|\nabla^2
v\|_{L^2}^2+C\int_0^{\sigma(T)}\left(\|\nabla
v_t\|_{L^2}^2+\|\nabla^3
v\|_{L^2}^2\right)dt\nonumber\\
&&\qquad+ C\sup_{0\leq t\leq \sigma(T)}\|\nabla^2
v\|_{L^2}^2\int_0^{\sigma(T)}\left(\|d_t\|_{L^2}^4+\|\nabla^2d\|_{L^2}^4\right)dt\nonumber\\
&&\quad\leq C \|\nabla^2 v_0\|_{L^2}^2,\label{3.78}
\end{eqnarray}
where we have used (\ref{1.3}), (\ref{3.12}) and (\ref{3.15}) to get
that
\begin{eqnarray*}
\int_0^T\|d_t\|_{L^2}^4dt&\leq&
C\int_0^T\left(\|\nabla^2d\|_{L^2}+\|u\nabla
d\|_{L^2}+\|\nabla d\|_{L^4}^2\right)^4dt\\
&\leq&C\int_0^T\left(\|\nabla^2d\|_{L^2}+\|\nabla u\|_{L^2}\|\nabla
d\|_{L^3}+\|\nabla d\|_{L^3}\|\nabla^2d\|_{L^2}\right)^4dt\\
&\leq&C\int_0^T\left(\|\nabla^2d\|_{L^2}^4+\|\nabla u\|_{L^2}^4
\right)dt\leq CC_0^{2\delta_0}.
\end{eqnarray*}

Similarly, multiplying (\ref{3.73}) by $t$, integrating it over
$(0,\sigma(T))$, and taking (\ref{3.74})--(\ref{3.77}) into account,
we deduce from (\ref{3.30}) and (\ref{3.62}) that
\begin{eqnarray}
&&\sup_{0\leq t\leq \sigma(T)}\left(t\|\nabla
w_3\|_{L^2}^2\right)+\int_0^{\sigma(T)}t\|\rho^{1/2}\dot
w_3\|_{L^2}^2dt\nonumber\\
&&\quad\leq C\sup_{0\leq t\leq \sigma(T)}\left(t\|\nabla^2
v\|_{L^2}^2\right)+C\int_0^{\sigma(T)}t\left(\|\nabla v_t\|_{L^2}^2+\|\nabla^3v\|_{L^2}^2\right)dt\nonumber\\
&&\qquad+C\sup_{0\leq t\leq \sigma(T)}\left(t\|\nabla^2
v\|_{L^2}^2\right)\int_0^{\sigma(T)}\left(\|d_t\|_{L^2}^4+\|\nabla^2
d\|_{L^2}^4\right)dt+C\int_0^{\sigma(T)}\|\nabla
w_3\|_{L^2}^2dt\nonumber\\
&&\quad\leq C \|\nabla v_0\|_{L^2}^2.\label{3.79}
\end{eqnarray}

By the same token as that in the proof of (\ref{3.66}), we conclude
from (\ref{3.78}) and (\ref{3.79}) that
\begin{equation}
\sup_{0\leq t\leq \sigma(T)}\left(t^{1-\beta}\|\nabla
w_3\|_{L^2}^2\right)+\int_0^{\sigma(T)}t^{1-\beta}\|\rho^{1/2}\dot
w_3\|_{L^2}^2dt\leq C \|\nabla v_0\|_{\dot H^\beta}^2.\label{3.80}
\end{equation}

Now, choosing $w_{10}=u_0$ and $v_0=d_0$ so that $w_1+w_2+w_3=u$ and
$v=d$, we immediately obtain (\ref{3.54}) from (\ref{3.66}),
(\ref{3.71}) and (\ref{3.80}).

To prove (\ref{3.55}), choosing $m=2-\beta$ in (\ref{3.47}),
(\ref{3.51}) and adding them together, we infer from (\ref{3.7}),
(\ref{3.11}), (\ref{3.21}) and (\ref{3.54}) that if $C_0\leq
\varepsilon_5$, then
\begin{eqnarray}
&&\sup_{0\leq t\leq \sigma(T)}t^{2-\beta}\left(\|\rho^{1/2}\dot
u\|_{L^2}^2+\|\nabla
d_t\|_{L^2}^2\right)\nonumber\\
&&\qquad+\int_0^{\sigma(T)}t^{2-\beta}\left(\|\nabla \dot u\|_{L^2}^2+\|d_{tt}\|_{L^2}^2+\|\nabla^2d_t\|_{L^2}^2\right)dt\nonumber\\
&&\quad\leq C\int_0^{\sigma(T)}\|\nabla
u\|_{L^2}^2dt+C\int_0^{\sigma(T)}t^{2-\beta}\|\nabla
u\|_{L^2}^2\|\nabla^2 d\|_{L^2}^2dt\nonumber\\
&&\qquad +C \int_0^{\sigma(T)}t^{2-\beta}\left(\|\nabla
u\|_{L^2}^4+\|\nabla d\|_{L^2}^4\right)\left(\|\rho^{1/2}\dot
u\|_{L^2}^2+\|\nabla d_t\|_{L^2}^2+\|\nabla^3d\|_{L^2}^2\right) dt\nonumber\\
&&\qquad+ C\int_0^{\sigma(T)}t^{1-\beta}\left(\|\rho^{1/2}\dot
u\|_{L^2}^2+\|\nabla d_t\|_{L^2}^2\right)dt+C
\int_0^{\sigma(T)}t^{2-\beta}\|\nabla
u\|_{L^4}^4dt\nonumber\\
&&\quad\leq C+C\int_0^{\sigma(T)}t^{2-\beta}\|\nabla
u\|_{L^4}^4dt,\label{3.81}
\end{eqnarray}
since it follows from (\ref{3.21}) and (\ref{3.54}) that for
$\beta\in(1/2,1]$,
\begin{eqnarray*}
&&\int_0^{\sigma(T)}t^{2-\beta}\left(\|\nabla u\|_{L^2}^4+\|\nabla
d\|_{L^2}^4\right)\left(\|\rho^{1/2}\dot u\|_{L^2}^2+\|\nabla
d_t\|_{L^2}^2+\|\nabla^3d\|_{L^2}^2\right) dt\nonumber\\
&&\quad\leq \sup_{0\leq t\leq
\sigma(T)}\left(t^{1-\beta}\left(\|\nabla
u\|_{L^2}^2+\|\nabla^2d\|_{L^2}^2\right)\right)^2\nonumber\\
&&\qquad\qquad\times\int_0^{\sigma(T)}t^{1-\beta}\left(\|\rho^{1/2}\dot
u\|_{L^2}^2+\|\nabla d_t\|_{L^2}^2+\|\nabla^3d\|_{L^2}^2\right)dt\\
&&\quad\leq C.
\end{eqnarray*}

Using Lemmas \ref{lem2.1}, \ref{2.2}, (\ref{3.7}), (\ref{3.11}),
(\ref{3.12}), (\ref{3.21}) and (\ref{3.54}), we have
\begin{eqnarray}
&&\int_0^{\sigma(T)}t^{2-\beta}\|\nabla
u\|_{L^4}^4dt\leq C\int_0^{\sigma(T)}t^{2-\beta}\|\nabla u\|_{L^2}\|\nabla u\|_{L^6}^3dt\nonumber\\
&&\quad\leq C\int_0^{\sigma(T)}t^{2-\beta}\|\nabla
u\|_{L^2}\left(1+\|\rho^{1/2}\dot u\|_{L^2}^3+\|\nabla^2d\|_{L^2}^{9/2}\|\nabla^3d\|_{L^2}^{3/2}\right)dt\nonumber\\
&&\quad\leq C+
C\int_0^{\sigma(T)}t^{\beta-1/2}\left(t^{1-\beta}\|\nabla
u\|_{L^2}^2\right)^{1/2}\left(t^{2-\beta}\|\rho^{1/2}\dot
u\|_{L^2}^2\right)^{1/2}\left(t^{1-\beta}\|\rho^{1/2}\dot
u\|_{L^2}^2\right) dt\nonumber\\
&&\qquad+C\int_0^{\sigma(T)}t^{2\beta-1}\left(\|\nabla
u\|_{L^2}^4\right)^{1/4}\left(t^{1-\beta}\|\nabla^2d\|_{L^2}^2\right)^{9/4}\left(t^{1-\beta}\|\nabla^3d\|_{L^2}^2\right)^{3/4}
dt
\nonumber\\
&&\quad\leq C+ C\sup_{0\leq t\leq
\sigma(T)}\left(\left(t^{1-\beta}\|\nabla u\|_{L^2}^2
\right)^{1/2}\left(t^{2-\beta}\|\rho^{1/2}\dot u\|_{L^2}^2
\right)^{1/2}\right)\int_0^{\sigma(T)}t^{1-\beta}\|\rho^{1/2}\dot
u\|_{L^2}^2dt\nonumber\\
&&\qquad+C\sup_{0\leq t\leq
\sigma(T)}\left(t^{1-\beta}\|\nabla^2d\|_{L^2}^2 \right)^{9/4}
\left(\int_0^{\sigma(T)}\|\nabla
u\|_{L^2}^4dt\right)^{1/4}\left(\int_0^{\sigma(T)}t^{1-\beta}\|\nabla^3d\|_{L^2}^2dt\right)^{3/4}\nonumber\\
&&\quad\leq C+C\sup_{0\leq t\leq
\sigma(T)}\left(t^{2-\beta}\|\rho^{1/2}\dot u\|_{L^2}^2
\right)^{1/2},\label{3.82}
\end{eqnarray}
provided $C_0\leq\varepsilon_5$. As a result, putting (\ref{3.82})
into (\ref{3.81}) and using Young's inequality, we immediately
obtain (\ref{3.55}). The proof of Lemma \ref{lem3.6} is therefore
completed.\hfill$\square$

\vskip 2mm

We are now in a position of estimating $A_5(\sigma(T))$.
\begin{lem}\label{lem3.6}
Let $(\rho,u,d)$ be a smooth solution of (\ref{1.1})--(\ref{1.7}) on
$\R^3\times(0,T]$ satisfying (\ref{3.7}). Then there exists a
constant $\varepsilon_6>0$, depending only on
$\mu,\lambda,A,\gamma,\bar\rho,\tilde\rho,\beta,M_1$ and $M_2$, such
that
\begin{equation}
A_5(\sigma(T))\triangleq\sup_{0\leq t\leq\sigma(T)}\int\rho
|u|^3dx\leq \frac{C_0^{\delta_0}}{2},\label{3.83}
\end{equation}
provided $C_0\leq\varepsilon_6$.
\end{lem}
\pf Multiplying (\ref{1.2}) by $3|u|u$ and integrating it by parts
over $\R^3\times(0,\sigma(T))$, we have from Cauchy-Schwarz
inequality that
\begin{eqnarray}
\sup_{0\leq t\leq\sigma(T)}\int\rho |u|^3dx &\leq& \int\rho_0
|u_0|^3dx+C\int_0^{\sigma(T)}\int |P(\rho)-P(\tilde\rho)||u||\nabla
u|dxdt\nonumber\\
&& +C\int_0^{\sigma(T)}\int|u| |\nabla
u|^2dxdt+C\int_0^{\sigma(T)}\int|\nabla d|^4|u| dxdt\nonumber\\
&\triangleq& \int\rho_0 |u_0|^3dx+\sum_{i=1}^3R_i .\label{3.84}
\end{eqnarray}

First, it is easily seen from (\ref{2.1}), (\ref{3.7}) and
(\ref{3.11}) that
\begin{equation}
R_1\leq
C\int_0^{\sigma(T)}\|P(\rho)-P(\tilde\rho)\|_{L^3}\|u\|_{L^6}\|\nabla
u\|_{L^2}dt\leq C\int_0^{\sigma(T)}\|\nabla u\|_{L^2}^2dt\leq
CC_0.\label{3.85}
\end{equation}
Keeping in mind that $\beta\in(1/2,1]$ and $\delta_0\in(0,1/9]$,
using Lemmas \ref{lem2.1}, \ref{lem2.2}, (\ref{3.7}), (\ref{3.11}),
(\ref{3.12}), (\ref{3.15}) and H\"{o}lder inequality, we deduce
\begin{eqnarray}
R_2&\leq& C\int_0^{\sigma(T)}\|u\|_{L^\infty} \|\nabla u\|_{L^2}^2dt \leq C\int_0^{\sigma(T)} \|\nabla u\|_{L^2}^{5/2}\|\nabla u\|_{L^6}^{1/2}dt\nonumber\\
&\leq& C\int_0^{\sigma(T)}\|\nabla u\|_{L^2}^{5/2}\left(\|\rho\dot
u\|_{L^2}+\|\nabla d\nabla^2
d\|_{L^2}+\|P(\rho)-P(\tilde\rho)\|_{L^6}\right)^{1/2}dt \nonumber\\
&\leq& C\int_0^{\sigma(T)}\|\nabla
u\|_{L^2}^{5/2}\left(\|\rho^{1/2}\dot u\|_{L^2}+\|\nabla
d\|_{L^3}\|\nabla^3
d\|_{L^2}+C_0^{1/6}\right)^{1/2}dt\nonumber\\
&\leq&C\int_0^{\sigma(T)}\sigma^{3(2\beta-3)/8}\left(\sigma^{(3-2\beta)/4}\|\nabla
u\|_{L^2}^2\right)^{5/4}\left(\sigma^{(3-2\beta)/4}\left(\|\rho^{1/2}\dot
u\|_{L^2}^2+\|\nabla^3d\|_{L^2}^2\right)\right)^{1/4}dt\nonumber\\
&&+CC_0^{1/12}\left(\int_0^{\sigma(T)} \|\nabla
u\|_{L^2}^4dt\right)^{5/8}\nonumber\\
&\leq& CC_0^{3\delta_0/2}+C\sup_{0\leq
t\leq\sigma(T)}\left(\sigma^{(3-2\beta)/4}\|\nabla
u\|_{L^2}^2\right)^{5/4}
\nonumber\\
&&\qquad
\times\left(\int_0^{\sigma(T)}\sigma^{(2\beta-3)/2}dt\right)^{3/4}\left(\int_0^{\sigma(T)}\sigma^{(3-2\beta)/4}\left(\|\rho^{1/2}\dot
u\|_{L^2}^2+\|\nabla^3d\|_{L^2}^2\right)dt\right)^{1/4}\nonumber\\
&\leq& CC_0^{3\delta_0/2},\label{3.86}
\end{eqnarray}
and similarly,
\begin{eqnarray*}
R_3&\leq& C\int_0^{\sigma(T)}\|\nabla u\|_{L^2}^{1/2}\|\nabla u\|_{L^6}^{1/2} \|\nabla d\|_{L^3}^{2}\|\nabla^2 d\|_{L^2}^2dt\nonumber\\
&\leq& CC_0^{2\delta_0/3}\int_0^{\sigma(T)}\|\nabla
u\|_{L^2}^{1/2}\|\nabla^2
d\|_{L^2}^{2}\left(C_0^{1/6}+\|\rho^{1/2}\dot u\|_{L^2} + \|\nabla^3
d\|_{L^2}\right)^{1/2}dt\nonumber\\
&\leq& CC_0^{1/12}\left(\int_0^{\sigma(T)}\|\nabla
u\|_{L^2}^2dt\right)^{1/4}\left(\int_0^{\sigma(T)}\|\nabla^2
d\|_{L^2}^4dt\right)^{1/2}\nonumber\\
&&+C\sup_{0\leq t\leq\sigma(T)}\left(\left(\sigma^{(3-2\beta)/4}\|
\nabla u\|_{L^2}^2\right)^{1/4}\left(\sigma^{(3-2\beta)/4}\|\nabla^2
d\|_{L^2}^2
\right)\right)\nonumber\\
&&\qquad\quad\times\int_0^{\sigma(T)}\sigma^{3(2\beta-3)/8}\left(\sigma^{(3-2\beta)/4}\left(\|\rho^{1/2}\dot
u\|_{L^2}^2+\|\nabla^3d\|_{L^2}^2\right)\right)^{1/4} dt\nonumber\\
&\leq& CC_0^{3\delta_0/2}.
\end{eqnarray*}
This, combined with (\ref{3.85}), (\ref{3.86}) and (\ref{3.84})
gives
\begin{eqnarray}
\sup_{0\leq t\leq
\sigma(T)}\int\rho|u|^3dx&\leq&\int\rho_0|u_0|^3dx+CC_0^{3\delta_0/2}\nonumber\\
&\leq& CC_0^{3(2\beta-1)/(4\beta)}+CC_0^{3\delta_0/2}\leq
C_1C_0^{3\delta_0/2},\label{3.87}
\end{eqnarray}
where we have used the interpolation and Sobolev embedding
inequalities to get that
\begin{eqnarray}
\left(\int\rho_0|u_0|^3dx\right)^{1/3}&\leq& C\left(\int\rho_0
|u_0|^2dx\right)^{(2\beta-1)/(4\beta)}\left(\int
|u_0|^{6/(3-2\beta)}dx\right)^{(3-2\beta)/(12\beta)}\nonumber\\
&\leq& CC_0^{(2\beta-1)/(4\beta)}\|u_0\|_{\dot
H^\beta}^{1/(2\beta)}\leq CC_0^{(2\beta-1)/(4\beta)}\leq
CC_0^{3\delta_0/2}.\label{3.88}
\end{eqnarray}
since it follows from (\ref{3.8}) that $3\delta_0/2\leq
(2\beta-1)/(4\beta)$. So, if $C_0$ is chosen to be such that
$$
C_0\leq\varepsilon_6\triangleq
\min\left\{\varepsilon_5,(2C_1)^{-2/\delta_0}\right\},
$$
then (\ref{3.83}) immediately follows from
(\ref{3.87}).\hfill$\square$

\vskip 2mm

We can now close the estimates of $A_1(T)$ and $A_2(T)$.
\begin{lem}\label{lem3.7} Let $(\rho,u,d)$ be a smooth solution of (\ref{1.1})--(\ref{1.7}) on $\R^3\times(0,T]$ satisfying (\ref{3.7}).
Then there exists a constant $\varepsilon_7>0$, depending only on
$\mu,\lambda,A,\gamma,\bar\rho,\tilde\rho,\beta,M_1$ and $M_2$, such
that
\begin{equation}
A_1(T)+A_2(T)\leq C_0^{1/2},\label{3.89}
\end{equation}
provided $C_0\leq\varepsilon_7$.
\end{lem}
\pf To prove (\ref{3.89}), we first utilize Lemmas \ref{lem2.1},
\ref{lem2.2}, (\ref{3.11}) and (\ref{3.15}) to get that
\begin{eqnarray*}
\|F\|_{L^4}^4+\|\omega\|_{L^4}^4&\leq& C\left(\|F\|_{L^2}\|\nabla F\|_{L^2}^3+\|\omega\|_{L^2}\|\nabla\omega\|_{L^2}^3\right)\\
&\leq& C \left(\|\nabla u\|_{L^2}+C_0^{1/2}\right)\left(
\|\rho^{1/2} \dot
u\|_{L^2}^{3}+C_0^{\delta_0}\|\nabla^3d\|_{L^2}^{3}\right),
\end{eqnarray*}
from which it follows that
\begin{eqnarray}
&&\int_0^T\sigma^2\left(\|F\|_{L^4}^4+\|\omega\|_{L^4}^4\right)dt\nonumber\\
&&\quad\leq C\int_0^T\sigma^2\|\nabla u\|_{L^2} \|\rho^{1/2} \dot
u\|_{L^2}^{3}dt+ CC_0^{\delta_0}\int_0^T\sigma^2\|\nabla u\|_{L^2}
\|\nabla^3d\|_{L^2}^{3}dt\nonumber\\
&&\qquad +CC_0^{1/2}\int_0^T\sigma^2\|\rho^{1/2} \dot
u\|_{L^2}^{3}dt+ CC_0^{1/2+\delta_0}\int_0^T\sigma^2
\|\nabla^3d\|_{L^2}^{3}dt\nonumber\\
&&\quad\triangleq\sum_{i=1}^4R_i.\label{3.90}
\end{eqnarray}

In view of the definitions of $A_i(T)$ ($i=1,\ldots,5$) in
(\ref{3.1})--(\ref{3.5}), we deduce from (\ref{3.7}), (\ref{3.54})
and H\"{o}lder inequality that
\begin{eqnarray}
R_1 &\leq&C
\int_0^{\sigma(T)}\sigma^{(6\beta-7)/8}\left(\sigma^{(3-2\beta)/4}\|\nabla
u\|_{L^2}^2\right)^{1/2}\left(\sigma^2 \|\rho^{1/2} \dot
u\|_{L^2}^{2}\right) \left(\sigma^{1-\beta}\|\rho^{1/2}\dot
u\|_{L^2}^2\right)^{1/2} dt\nonumber\\
&&+ C\int_{\sigma(T)}^T\left(\sigma \|\nabla
u\|_{L^2}^2\right)^{1/2}\left(\sigma^2 \|\rho^{1/2} \dot
u\|_{L^2}^2\right)^{1/2} \left(\sigma\|\rho^{1/2}\dot
u\|_{L^2}^2\right) dt\nonumber\\
&\leq&
CA_2(\sigma(T))A_4^{1/2}(\sigma(T))\left(\int_0^{\sigma(T)}\sigma^{1-\beta}\|\rho^{1/2}\dot
u\|_{L^2}^2dt\right)^{1/2}\left(\int_0^{\sigma(T)}\sigma^{(6\beta-7)/4}dt\right)^{1/2}\nonumber\\
&&+CA_1^{1/2}(T)A_2^{1/2}(T)\int_{\sigma(T)}^T
\sigma\|\rho^{1/2}\dot
u\|_{L^2}^2 dt\nonumber\\
&\leq&CC_0^{(1+\delta_0)/2}+CC_0\leq
CC_0^{(1+\delta_0)/2}.\label{3.91}
\end{eqnarray}
Similarly, due to (\ref{3.7}) and (\ref{3.21}), we have
\begin{eqnarray}
\sum_{i=2}^4R_i&\leq&
CC_0^{\delta_0}\int_0^T\left(\sigma \|\nabla u\|_{L^2}^2\right)^{1/2}\left(\sigma^2\|\nabla^3d\|_{L^2}^2\right)^{1/2}\left(\sigma^{1-\beta}\|\nabla^3d\|_{L^2}^2\right)dt\nonumber\\
&&+CC_0^{1/2}\int_0^T\left(\sigma^2\|\rho^{1/2}\dot
u\|_{L^2}^2\right)^{1/2}\left(\sigma\|\rho^{1/2}\dot
u\|_{L^2}^2\right)dt\nonumber\\
&&+CC_0^{1/2+\delta_0}\int_0^T \left(\sigma^2\|\nabla^3d\|_{L^2}^2\right)^{1/2}\left(\sigma \|\nabla^3d\|_{L^2}^2\right) dt\nonumber\\
&\leq&C\left(C_0^{\delta_0}A_1^{1/2}(T)A_2^{1/2}(T)+C_0^{1/2}A_1(T)A_2^{1/2}(T)+C_0^{1/2+\delta_0}A_1(T)A_2^{1/2}(T)\right)\nonumber\\
&\leq& CC_0^{1/2+\delta_0}.\label{3.92}
\end{eqnarray}

Thus, substituting (\ref{3.91}), (\ref{3.92}) into (\ref{3.90})
gives
\begin{equation}
\int_0^T\sigma^2\left(\|F\|_{L^4}^4+\|\omega\|_{L^4}^4\right)dt\leq
CC_0^{(1+\delta_0)/2}.\label{3.93}
\end{equation}

In terms of the effective viscous flux $F$, we can rewrite
(\ref{3.36}) as
\begin{eqnarray}
0&=&(P(\rho)-P(\tilde\rho))_t+u\cdot\nabla
(P(\rho)-P(\tilde\rho))+\gamma P(\tilde\rho)\divg
u+\gamma(P(\rho)-P(\tilde\rho))\divg
u\nonumber\\
&=&(P(\rho)-P(\tilde\rho))_t+u\cdot\nabla
(P(\rho)-P(\tilde\rho))+\gamma P(\tilde\rho)\divg
u\nonumber\\
&&+\frac{\gamma}{2\mu+\lambda}(P(\rho)-P(\tilde\rho))F+\frac{\gamma}{2\mu+\lambda}(P(\rho)-P(\tilde\rho))^2,\label{3.94}
\end{eqnarray}
which, multiplied by $3\sigma^2(P(\rho)-P(\tilde\rho))^2$ and
integrated by parts over $\R^3\times(0,T)$, yields
\begin{eqnarray}
&&\frac{3\gamma-1}{2\mu+\lambda}\int_0^T\sigma^2\|P(\rho)-P(\tilde\rho)\|_{L^4}^4dt\nonumber\\
&&\quad\leq
C\sigma^2\|P(\rho)-P(\tilde\rho)\|_{L^3}^3+C\int_0^T\left(\sigma'\|P(\rho)-P(\tilde\rho)\|_{L^3}^3+\sigma^2\|\nabla
u\|_{L^2}^2+\sigma^2\|F\|_{L^4}^4\right)dt\nonumber\\
&&\quad \leq CC_0^{(1+\delta_0)/2}.\label{3.95}
\end{eqnarray}
where we have used (\ref{3.11}), (\ref{3.93}) and Cauchy-Schwarz
inequality.

As a result of (\ref{3.93}) and (\ref{3.95}), we deduce from the
standard $L^p$-estimate that
\begin{eqnarray}
\int_0^T\sigma^2\|\nabla u\|_{L^4}^4dt&\leq&
C\int_0^T\sigma^2\left(\|\divg
u\|_{L^4}^4+\|\omega\|_{L^4}^4\right)dt\nonumber\\
&\leq& C\int_0^T\sigma^2\left(\|F\|_{L^4}^4+\|P(\rho)-P(\tilde\rho)\|_{L^4}^4+\|\omega\|_{L^4}^4\right)dt\nonumber\\
&\leq& CC_0^{(1+\delta_0)/2}.\label{3.96}
\end{eqnarray}

Thus, it follows from (\ref{3.32}), (\ref{3.33}), (\ref{3.95}) and
(\ref{3.96}) that
$$
A_1(T)+A_2(T)\leq CC_0^{1/2+2\delta_0}+CC_0^{(1+\delta_0)/2}\leq
C_1C_0^{(1+\delta_0)/2}\leq C_0^{1/2},
$$
provided $C_0$ is chosen to be such that
$$
C_0\leq\varepsilon_7\triangleq\min\left\{\varepsilon_6,C_1^{-2/\delta_0}\right\}.
$$
The proof of Lemma \ref{lem3.7} is therefore
completed.\hfill$\square$

\vskip 2mm

Finally, we need to prove the uniform upper bound of the density.
\begin{lem}\label{lem3.8} Let $(\rho,u,d)$ be a smooth solution of (\ref{1.1})--(\ref{1.7}) on $\R^3\times(0,T]$ satisfying (\ref{3.7}).
Then there exists a constant $\varepsilon_8>0$, depending only on
$\mu,\lambda,A,\gamma,\bar\rho,\tilde\rho,\beta,M_1$ and $M_2$, such
that
\begin{equation}
\sup_{0\leq t\leq T}\|\rho\| _{L^\infty}
\leq\frac{7}{4}\bar\rho\label{3.97}
\end{equation}
provided $C_0\leq\varepsilon_8$.
\end{lem}
\pf  Let $D_t\triangleq\partial_t+u \cdot\nabla $ denote the
material derivative operator. Then, in terms of the effective viscous flux
$F$ in (\ref{2.3}), we can rewrite (\ref{1.1}) as
$$
D_t\rho=g(\rho)+b'(\rho)
$$
where
$$
g(\rho)=-\frac{A\rho}{2\mu+\lambda}(\rho^\gamma-\tilde\rho^\gamma),\quad
b(t)=-\frac{1}{2\mu+\lambda}\int_0^t\rho Fds.
$$

Thus, to apply Lemma \ref{lem2.3}, we now need to estimate $b(t)$.
To do this, we first use  (\ref{2.2}) and (\ref{a26}) to deduce that
for any $0\leq t_1<t_2\leq \sigma(T)$,
\be\ba
  |b(t_2)-b(t_1)|
 &\leq C\int_0^{\sigma(T)}\|F\|_{L^\infty}dt\\&=C\int_0^{\sigma(T)}\|\Delta^{-1}\div (\n \dot u)+\Delta^{-1}\div \div(M(d))\|_{L^\infty}dt\\ &\le  C\int_0^{\sigma(T)}\left(\| \n \dot u\|_{L^2}^{1/2}\|\dot u\|_{L^6}^{1/2}+ \|M(d)\|_{L^{3/2}}^{1/5}\|\na M(d)\|_{L^6}^{4/5}\right)dt\\ &\le   C \int_0^{\sigma(T)} \| \n \dot u\|_{L^2}^{1/2}\|\na\dot u\|_{L^2}^{1/2}dt+C\int_0^{\sigma(T)}\|\na d\|_{L^3}^{2/3}\|\na^2d\|_{L^6}^{4/3} dt\\
 &\triangleq R_1+R_2.\label{3.98}
\ea\ee
On one hand,
using (\ref{3.54}), (\ref{3.55}) and (\ref{3.89}), we find
\begin{eqnarray*}
R_1&\leq&C \int_0^{\sigma(T)} t^{(\beta-2)/4}\|\rho^{1/2}\dot
u\|_{L^2}^{1/2} \left(t^{2-\beta}\|\nabla\dot
u\|_{L^2}^2\right)^{1/4} dt\\
&\leq& C\left(\int_0^{\sigma(T)} t^{(\beta-2)/3}\|\rho^{1/2}\dot
u\|_{L^2}^{2/3} dt\right)^{3/4}\left(\int_0^{\sigma(T)}
t^{2-\beta}\|\nabla\dot
u\|_{L^2}^2dt\right)^{1/4}\\
&\leq& C\left(\int_0^{\sigma(T)}
t^{-[(2-\beta)(2/3-\delta_0)+\delta_0]}\left(t^{2-\beta}\|\rho^{1/2}\dot
u\|_{L^2}^2\right)^{1/3-\delta_0}\left(t \|\rho^{1/2}\dot
u\|_{L^2}^2\right)^{\delta_0} dt\right)^{3/4}\\
&\leq&C\left(\int_0^{\sigma(T)}
t^{-[2(2-\beta)+3(\beta-1)\delta_0]/(3(1-\delta_0))}dt\right)^{3(1-\delta_0)/4}
\left(\int_0^{\sigma(T)} t \|\rho^{1/2}\dot
u\|_{L^2}^2  dt\right)^{3\delta_0/4}\nonumber\\
&\leq& CA_1^{3\delta_0/4}(T)\leq CC_0^{3\delta_0/8},
\end{eqnarray*}
since it follows from (\ref{3.8}) that
$$
-\frac{2(2-\beta)+3(\beta-1)\delta_0}{3(1-\delta_0)}
=-\frac{9\beta-4\beta^2+1}{7\beta+1}>-1.
$$  On the other hand, it follows from (\ref{3.15}) and (\ref{3.21}) that
\begin{eqnarray*}
R_2&\leq&
CC_0^{2\delta_0/3}\int_0^{\sigma(T)}t^{-2(1-\beta)/3}\left(t^{1-\beta}
\|\nabla^3d\|_{L^2}^2\right)^{2/3} dt \\
&\leq& CC_0^{2\delta_0/3}\left(\int_0^{\sigma(T)}
t^{2(\beta-1)}dt\right)^{1/3}\left(\int_0^{\sigma(T)} t^{1-\beta}
\|\nabla^3d\|_{L^2}^2 dt\right)^{2/3}\\
&\leq& CC_0^{2\delta_0/3}.
\end{eqnarray*}

Thus, putting the estimates of $R_1,R_2$  into
(\ref{3.98}) gives
$$
|b(t_2)-b(t_1)|\leq CC_0^{\delta_0/3}.
$$
So, for $t\in[0,\sigma(T)]$ one can choose $N_0,N_1$ and $\xi^*$ in
Lemma \ref{lem2.3} as follows:
$$
N_0=CC_0^{\delta_0/3},\quad N_1=0,\quad \xi^*=\tilde\rho.
$$

Since it holds that
$$
g(\xi)=-\frac{A\xi}{2\mu+\lambda}\left(\xi^\gamma-\tilde\rho^\gamma\right)\leq
-N_1=0,\quad \forall\; \xi\geq\xi^*=\tilde\rho,
$$
we thus conclude from (\ref{2.8}) that
\begin{equation}
\sup_{0\leq t\leq\sigma(T)}\|\rho\| _{L^\infty}\leq
\max\{\bar\rho,\tilde\rho\}+N_0\leq
\bar\rho+CC_0^{\delta_0/3}\leq\frac{3}{2}\bar\rho,\label{3.99}
\end{equation}
provided $C_0$ is chosen to be such that
$$
C_0\leq\varepsilon_{8,1}\triangleq\min\left\{\varepsilon_7,
\left(\frac{\bar\rho}{2C}\right)^{3/\delta_0}\right\}.
$$

On the other hand, it follows from Lemmas \ref{lem2.1},
\ref{lem2.2}, (\ref{3.7}), (\ref{3.11}) and (\ref{3.89}) that
\begin{eqnarray*}
 \|F\|_{L^\infty}&\leq& C\|F\|_{L^2}^{1/4}\|\nabla
F\|_{L^6}^{3/4}\leq C C_0^{1/16}\|\nabla F\|_{L^6}^{3/4}\\
&\leq& C C_0^{1/16}\left(\|\nabla \dot u\|_{L^2}+\|\nabla^2
d\|_{L^2}^{1/2}\|\nabla^3d\|_{L^2}^{3/2}\right)^{3/4}
\end{eqnarray*}
for any $t\in[\sigma(T),T]$. Hence, using (\ref{3.22}) and
(\ref{3.89}), we deduce that for any $t_1,t_2\in[\sigma(T),T]$,
\begin{eqnarray*}
|b(t_2)-b(t_1)|&\leq&C\int_{t_1}^{t_2}\|F\|_{L^\infty}dt\leq\frac{A}{2\mu+\lambda}(t_2-t_1)+C\int_{\sigma(T)}^{T}\|F\|_{L^\infty}^{8/3}dt\\
&\leq&\frac{A}{2\mu+\lambda}(t_2-t_1)+CC_0^{1/6}\int_{\sigma(T)}^{T}\left(\sigma^2\|\nabla
\dot u\|_{L^2}^{2}+\sigma^{5/2}\|\nabla^2
d\|_{L^2}\|\nabla^3d\|_{L^2}^3\right)dt\\
&\leq&\frac{A}{2\mu+\lambda}(t_2-t_1)+CC_0^{1/6}.
\end{eqnarray*}
Thus, for $t\in[\sigma(T),T]$ one can choose $N_0,N_1$ and $\xi^*$
in Lemma \ref{lem2.3} as follows:
$$
N_0=CC_0^{1/6},\quad N_1=\frac{A}{2\mu+\lambda},\quad
\xi^*=\tilde\rho+1.
$$

Noting that
$$
g(\xi)=-\frac{A\xi}{2\mu+\lambda}\left(\xi^\gamma-\tilde\rho^\gamma\right)\leq
-N_1=-\frac{A}{2\mu+\lambda},\quad \forall\; \xi\geq\xi
=\tilde\rho+1,
$$
we can thus apply Lemma \ref{lem2.3} to get
\begin{equation}
\sup_{\sigma(T)\leq t\leq T}\|\rho\| _{L^\infty}\leq
\max\left\{\frac{3\bar\rho}{2},\tilde\rho+1\right\}+N_0\leq
\frac{3\bar\rho}{2}+CC_0^{1/6}\leq\frac{7}{4}\bar\rho,\label{3.100}
\end{equation}
provided $C_0$ is chosen to be such that
$$
C_0\leq \varepsilon_8\triangleq\min\left\{\varepsilon_{8,1},
\left(\frac{\bar\rho}{4C}\right)^{6}\right\}.
$$
Therefore, the combination of (\ref{3.99}) with (\ref{3.100}) ends the proof
of Lemma \ref{lem3.8}.\hfill$\square$

\vskip 2mm

Now, by virtue of Lemmas \ref{lem3.1}--\ref{lem3.8} we can complete
the proof of Proposition \ref{pro3.1}.

\vskip 2mm

{\it Proof of Proposition \ref{pro3.1}.} By Lemmas \ref{lem3.2}, \ref{lem3.6}-\ref{lem3.8}, to prove  Proposition \ref{pro3.1}, it remains to estimate the term $A_4(\si(T)).$ In fact, using (\ref{3.21}),
(\ref{3.54}) and (\ref{3.89}), we have
\begin{eqnarray*}
A_4(\si(T))
 &&\leq\sup_{0\leq t\leq\sigma(T)}
\left(\sigma^{1-\beta}\|\nabla
u\|_{L^2}^2\right)^{(2\beta+1)/(4\beta)} \sup_{0\leq t\leq\sigma(T)}
\left(\sigma \|\nabla
u\|_{L^2}^2\right)^{(2\beta-1)/(4\beta)} \nonumber\\
&&\quad+\sup_{0\leq t\leq\sigma(T)}
\left(\sigma^{1-\beta}\|\nabla^2
d\|_{L^2}^2\right)^{(2\beta+1)/(4\beta)} \sup_{0\leq t\leq\sigma(T)}
\left(\sigma \|\nabla^2
d\|_{L^2}^2\right)^{(2\beta-1)/(4\beta)} \nonumber\\
&&\quad+\left(\int_0^{\sigma(T)}\sigma^{1-\beta}\|\rho^{1/2}\dot
u\|_{L^2}^2dt\right)^{(2\beta+1)/(4\beta)}\left(\int_0^{\sigma(T)}\sigma
\|\rho^{1/2}\dot
u\|_{L^2}^2dt\right)^{(2\beta-1)/(4\beta)}\nonumber\\
&&\quad+\left(\int_0^{\sigma(T)}\sigma^{1-\beta}\|\nabla^3d\|_{L^2}^2dt\right)^{(2\beta+1)/(4\beta)}\left(\int_0^{\sigma(T)}\sigma
\|\nabla^3d\|_{L^2}^2dt\right)^{(2\beta-1)/(4\beta)}\nonumber\\
&&\leq CA_1^{(2\beta-1)/(4\beta)}(T)\leq
CC_0^{(2\beta-1)/(8\beta)}\leq \frac{C_0^{\delta_0}}{2},
\end{eqnarray*}
provided $C_0$ is chosen to be such that
$$
C_0\leq \varepsilon\triangleq\min\left\{\varepsilon_8,(2C)^{-72\beta/(2\beta-1)}\right\}.
$$
Therefore, the proof of Proposition \ref{pro3.1} is
completed.\hfill$\square$

\section{Time-dependent higher-order estimates}
In this section, we prove the global estimates on the spatial-time
derivatives of $(\rho,u,d)$ which are needed to guarantee the
existence of classical solutions. For this purpose, we assume that
the conditions of Theorem \ref{thm1.1} hold and that $(\rho,u,d)$ is
a smooth solution of (\ref{1.1})--(\ref{1.7}) on $\R^3\times(0,T]$
satisfying Proposition \ref{pro3.1}. Moreover, from now on we will
always assume that the initial energy $C_0$ satisfies (\ref{3.10}).
For simplicity, throughout this section we denote by $C$ or $C_i$
($i=1,2,3,\ldots$) the various positive constants which may depend
on
$$
\mu,\lambda,A,\gamma,\bar\rho,\tilde\rho,\beta,M_1, M_2,g\;\;{\rm
and}\;\; T,
$$
where $g\in L^2$ is the function in the compatibility condition
(\ref{1.12}) and $T>0$ is the time.

First, one easily deduces from Lemmas \ref{lem3.3}--\ref{lem3.5}
that
\begin{lem}\label{lem4.1}Assume that the conditions of Theorem \ref{thm1.1} hold. Then for any given $T>0$,
there exists a constant $C(T)>0$ such that
\begin{eqnarray}
&&\sup_{0\leq t\leq T}\left(\|\nabla^2d\|_{L^2}^2+\|\nabla
u\|_{L^2}^2\right)\nonumber\\
&&\quad+\int_0^T\left(\|\rho^{1/2}\dot u\|_{L^2}^2+\|\nabla
d_t\|_{L^2}^2+\|\nabla^3 d\|_{L^2}^2\right)dt\leq C(T)\label{4.1}
\end{eqnarray}
and
\begin{eqnarray}
&&\sup_{0\leq t\leq T} \left(\|\rho^{1/2}\dot
u\|_{L^2}^2+\|\nabla^3d\|_{L^2}^2+\|\nabla
d_t\|_{L^2}^2\right)\nonumber\\
&&\quad+\int_0^T \left(\|\nabla \dot u\|_{L^2}^2+\|d_{tt}\|_{L^2}^2
+\|\nabla^2d_{t}\|_{L^2}^2 \right)dt \leq C(T).\label{4.2}
\end{eqnarray}
\end{lem}
\pf Indeed, choosing $\beta=1$ in  (\ref{3.21}) and (\ref{3.54}),
one immediately obtains (\ref{4.1}) from (\ref{3.21}), (\ref{3.54})
and (\ref{3.89}).

In order to prove (\ref{4.2}), we first notice that (\ref{3.53}),
together with (\ref{4.1}), implies
$$
\|\nabla^3d\|_{L^2}^2\leq C\left(1+\|\nabla d_t\|_{L^2}^2\right).
$$
Taking this into account, choosing $m=0$ in (\ref{3.47}),
(\ref{3.51}) and integrating them over $(0,T)$, we have by
(\ref{4.1}) and the compatibility condition (\ref{1.12}) that
\begin{eqnarray}
&&\sup_{0\leq t\leq T} \left(\|\rho^{1/2}\dot u\|_{L^2}^2
+\|\nabla^3d\|_{L^2}^2+\|\nabla
d_t\|_{L^2}^2\right)\nonumber\\
&&\qquad+\int_0^T \left(\|\nabla \dot u\|_{L^2}^2+\|d_{tt}\|_{L^2}^2+\|\nabla^2d_t\|_{L^2}^2\right)dt\nonumber\\
&&\quad\leq C+C\int_0^T\|\nabla u\|_{L^4}^4dt\leq C+C\int_0^T\|\nabla u\|_{L^2}\|\nabla u\|_{L^6}^3dt\nonumber\\
&&\quad\leq C+C\int_0^T\left(\|\rho^{1/2}\dot u\|_{L^2}^3+\|\nabla^3
d \|_{L^2}^3\right)dt,\label{4.3}
\end{eqnarray}
where we have also used (\ref{2.1}), (\ref{2.6}) and (\ref{3.15}).
Thus, in view of (\ref{4.1}) and (\ref{4.3}), an application of
Gronwall's inequality immediately leads to
(\ref{4.2}).\hfill$\square$

\vskip 2mm

Next, we  use   the Beale-Kato-Majda-type inequality (see Lemma
\ref{lem2.4}) to derive the estimates on the spatial gradient of the
density as well as the velocity, which are very important for the
proof of the higher-order estimates of the solutions. The proofs are
mainly motivated by the ones in \cite{HLX2010,HLX2010-2}.
\begin{lem}\label{lem4.2}Assume that the conditions of Theorem \ref{thm1.1} hold. Then for any given $T>0$,
there exists a constant $C(T)>0$ such that
\begin{equation}
\sup_{0\leq t\leq T}\left(\|\nabla\rho\|_{L^2{\cap}L^6}+\|\nabla
u\|_{H^1}\right)+\int_0^T\|\nabla u\|^{3/2}_{L^\infty}dt\leq
C(T).\label{4.4}
\end{equation}
\end{lem}
\pf For $2\leq p\leq 6$, it is easily derived from (\ref{1.1}) that
\begin{eqnarray}
\frac{d}{dt}\|\nabla\rho\|_{L^p}&\leq& C\|\nabla
u\|_{L^\infty}\|\nabla\rho\|_{L^p}+C\|\nabla^2u\|_{L^p}
\nonumber \\
&\leq& C\left(1+\|\nabla
u\|_{L^\infty}\right)\|\nabla\rho\|_{L^p}+C\left(1+\|\rho\dot
u\|_{L^p}\right),\label{4.5}
\end{eqnarray}
where we have used (\ref{4.1}), (\ref{4.2}) and the $L^p$-estimate
of elliptic system to infer from (\ref{1.2}) that
\begin{eqnarray}
\|\nabla^2 u\|_{L^p}&\leq& C\left(\|\rho\dot
u\|_{L^p}+\|\nabla\rho\|_{L^p}+\|\nabla
d\nabla^2d\|_{L^p}\right)\nonumber\\
&\leq&C\left(\|\rho\dot u\|_{L^p}+\|\nabla\rho\|_{L^p}+1\right).
\label{4.6}
\end{eqnarray}

We are now in a position of estimating $\|\nabla u\|_{L^\infty}$. By
(\ref{2.1}), (\ref{4.1}) and (\ref{4.6}), we deduce from
(\ref{2.10}) with $q=6$ that
\begin{eqnarray}
\|\nabla u\|_{L^\infty} &\leq& C+C\left(\|\divg
u\|_{L^\infty}+\|\nabla\times
u\|_{L^\infty}\right)\ln\left(e+\|\rho\dot
u\|_{L^6}+\|\nabla\rho\|_{L^6}\right)\nonumber\\
&\leq& C+C\left(\|\divg u\|_{L^\infty}+\|\nabla\times
u\|_{L^\infty}\right)\ln\left(e+\|\nabla\dot
u\|_{L^2}\right)\nonumber\\
&&+ C\left(\|\divg u\|_{L^\infty}+\|\nabla\times
u\|_{L^\infty}\right)\ln\left(e+\|\nabla\rho\|_{L^6}\right).\label{4.7}
\end{eqnarray}

So if we set
$$
\Phi(t)\triangleq e+\|\nabla\rho\|_{L^6},\quad
\Psi(t)\triangleq1+\left(\|\divg  u\|_{L^\infty}+\|\nabla\times
u\|_{L^\infty}\right)\ln\left(e+\|\nabla\dot
u\|_{L^2}\right)+\|\nabla\dot u\|_{L^2},
$$
then it follows from (\ref{4.5}) with $p=6$ and (\ref{4.7}) that
$$
\Phi'(t)\leq C\Psi(t)\Phi(t)\ln\Phi(t),
$$
due to the fact that $\ln \Phi(T)\geq1$. This particularly implies
\begin{equation}
\frac{d}{dt}\ln\Phi(t)\leq  C\Psi(t)\ln\Phi(t).\label{4.8}
\end{equation}

By virtue of Lemmas \ref{lem2.1}, \ref{lem2.2}, (\ref{4.1}) and
(\ref{4.2}), we find ($\omega=\nabla\times u$)
\begin{eqnarray}
\int_0^T\Psi^{3/2}(t)dt&\leq& C+C\int_0^T\left(\|\nabla \dot
u\|_{L^2}^2+\|\divg u\|_{L^\infty}^2+\|\omega\|_{L^\infty}^2\right)dt\nonumber\\
&\leq&C+C\int_0^T\left(\|P(\rho)-P(\tilde\rho)\|_{L^\infty}^2+\|F\|_{L^\infty}^2
+\| \omega\|_{L^\infty}^2\right)dt\nonumber\\
&\leq&C+C\int_0^T\left(\|F\|_{L^2}^2+\|\omega\|_{L^2}^2+\|\nabla F\|_{L^6}^{2}+\|\nabla\omega\|_{L^6}^{2}\right)dt\nonumber\\
&\leq&C+C\int_0^T\left(\|\nabla\dot u\|_{L^2}^2+\|\nabla
d\|_{L^\infty}^2\|\nabla^3d\|_{L^2}^2\right)dt\leq C,\label{4.9}
\end{eqnarray}
and consequently, it follows from (\ref{4.8}) and Gronwall's
inequality that
\begin{equation}
\sup_{0\leq t\leq T}\|\nabla\rho(t)\|_{L^6}\leq\sup_{0\leq t\leq
T}\Phi(t)\leq C.\label{4.10}
\end{equation}

As a result of (\ref{4.7}), (\ref{4.9}) and (\ref{4.10}), we obtain
\begin{equation}
\int_0^T\|\nabla u\|^{3/2}_{L^\infty}dt\leq  C+C\int_0^T\Psi^{3/2}(t)dt\leq
C.\label{4.11}
\end{equation}
Using (\ref{4.1}) and (\ref{4.11}), we also infer from (\ref{4.5})
with $p=2$ and Gronwall's inequality that
$$
\sup_{0\leq t\leq T}\|\nabla\rho(t)\|_{L^2}\leq C,
$$
which, together with (\ref{4.2}) and (\ref{4.6}), yields $\|\nabla
u\|_{H^1}\leq C$. The proof of Lemma \ref{lem4.2} is therefore
completed.\hfill$\square$

\vskip 2mm

Basing on Lemmas \ref{lem4.1} and \ref{lem4.2}, we easily obtain
\begin{lem}\label{lem4.3}Assume that the conditions of Theorem \ref{thm1.1} hold. Then for any given
$T>0$,
\begin{eqnarray}
&&\sup_{0\leq t\leq T}\|\rho^{1/2}u_t\|_{L^2}^2+\int_0^T\|\nabla
u_t\|_{L^2}^2dt\leq C(T),\label{4.12}\\
&&\sup_{0\leq t\leq T}\left(\|\nabla \rho\|_{H^1}+\|\nabla
P\|_{H^1}\right)+\int_0^T\|\nabla^2 u\|_{H^1}^2dt\leq
C(T),\label{4.13}\\
&&\sup_{0\leq t\leq
T}\left(\|\rho_t\|_{H^1}+\|P_t\|_{H^1}\right)+\int_0^T\left(\|\rho_{tt}\|_{L^2}^2+\|P_{tt}\|_{L^2}^2\right)dt\leq
C(T).\label{4.14}
\end{eqnarray}
\end{lem}
\pf First, due to $u_t=\dot u-u\cdot\nabla u$, one easily gets
(\ref{4.12}) from Lemmas \ref{lem2.1}, \ref{lem4.1} and
\ref{lem4.2}.

Next we prove (\ref{4.13}). To do this, noting that
$P(\rho)=A\rho^\gamma$ satisfies
\begin{equation}
P_t+u\cdot\nabla P+\gamma P\divg u=0, \label{4.15}
\end{equation}
from which and (\ref{1.1}) we have by some direct computations that
\begin{eqnarray}
&&\frac{d}{dt}\left(\|\nabla^2\rho\|_{L^2}^2+\|\nabla^2
P\|_{L^2}^2\right)\nonumber\\
&&\quad\leq C\|\nabla
u\|_{L^\infty}\left(\|\nabla^2\rho\|_{L^2}^2+\|\nabla^2
P\|_{L^2}^2\right)+C\|\nabla^3
u\|_{L^2}\left(\|\nabla^2\rho\|_{L^2}+\|\nabla^2
P\|_{L^2}\right)\nonumber\\
&&\qquad+C\|\nabla^2 u\|_{L^6}\left(\|\nabla\rho\|_{L^3}+\|\nabla
P\|_{L^3}\right)\left(\|\nabla^2\rho\|_{L^2}+\|\nabla^2
P\|_{L^2}\right)\nonumber\\
&&\quad\leq C\left(1+\|\nabla
u\|_{L^\infty}\right)\left(\|\nabla^2\rho\|_{L^2}^2+\|\nabla
P\|_{L^2}^2\right)+C\|\nabla^2 u\|_{H^1}^2,\label{4.16}
\end{eqnarray}
where we have used (\ref{2.1}) and (\ref{4.4}). Using (\ref{4.1}),
(\ref{4.2}), (\ref{4.4}), (\ref{4.12}) and Lemma \ref{lem2.1}, we
deduce from (\ref{1.2}) and the standard $L^2$-estimate of elliptic
system that
\begin{eqnarray}
\|\nabla^2 u\|_{H^1}&\leq& C\left(\|\rho\dot u\|_{H^1} +\|\nabla
P\|_{H^1}+\|\nabla d\nabla^2d\|_{H^1}\right)\nonumber\\
&\leq& C\left(1+\|\nabla u_t\|_{L^2} +\|\nabla^2 P\|_{L^2}
\right).\label{4.17}
\end{eqnarray}
Thus, combining (\ref{4.17}) with (\ref{4.16}) and using Gronwall
inequality, we immediately arrive at (\ref{4.13}) since it follows
from (\ref{4.4}) and (\ref{4.12}) that $$\|\nabla
u\|_{L^\infty}^{3/2}+\|\nabla   u_t\|_{L^2}^2\in L^1(0,T).$$

Finally, it is easily seen from (\ref{1.1}) and (\ref{4.15}) that
\begin{eqnarray}
\|\rho_t\|_{H^1}+\|P_t\|_{H^1}&\leq&
C\|u\|_{L^\infty}\left(\|\nabla\rho\|_{H^1}+\|\nabla
P\|_{H^1}\right)+C\|\nabla u\|_{H^1}\nonumber\\
&&+C\|\nabla u\|_{L^6}\left(\|\nabla\rho\|_{L^3}+\|\nabla
P\|_{L^3}\right)\nonumber\\
&\leq&C\left(1+\|\nabla u\|_{H^1}^2+\|\nabla\rho\|_{H^1}^2+\|\nabla
P\|_{H^1}^2\right)\leq C,\label{4.18}
\end{eqnarray}
where we have used (\ref{4.4}), (\ref{4.13}) and Lemma \ref{lem2.1}.
Moreover, noticing that (\ref{4.15}) implies
$$
P_{tt}+u_t\cdot\nabla P+u\cdot\nabla P_t+\gamma P_t\divg u+\gamma
P\divg u_t=0,
$$
and hence, using (\ref{4.4}), (\ref{4.12}), (\ref{4.13}) and
(\ref{4.18}), one obtains
\begin{eqnarray}
\int_0^T\|P_{tt}\|_{L^2}^2dt&\leq&
C\int_0^T\left(\|u_t\|_{L^6}\|\nabla
P\|_{L^3}+\|u\|_{W^{1,\infty}}\|P_t\|_{H^1}+\|\nabla
u_t\|_{L^2}\right)^2dt\nonumber\\
&\leq& C+C\int_0^T\left(\|\nabla u\|_{H^2}^2+\|\nabla
u_t\|_{L^2}^2\right)dt\leq C.\label{4.19}
\end{eqnarray}
Analogously, one also has $\|\rho_{tt}\|_{L^2}\in L^2(0,T)$. So,
combining this with (\ref{4.18}), (\ref{4.19}) completes the proof
of (\ref{4.14}).\hfill$\square$

\vskip 2mm

Due to the weaker compatibility condition (\ref{1.12}), the methods
used in \cite{HLX2010} to derive the higher-order estimates on the
solutions $(\rho,u,d)$, which guarantee the solutions obtained are
indeed a classical one away from the initial time, cannot be applied
any more. To overcome this difficulty, we need some careful
initial-layer analysis.
\begin{lem}\label{lem4.4}Let $\sigma\triangleq\min\{1,t\}$. Then for any given
$T>0$,
\begin{eqnarray}
&&\sup_{0\leq t\leq T} \sigma\left(\|\nabla u\|_{H^2}^2+\|\nabla
u_t\|_{L^2}^2+\|\nabla^2d_t\|_{L^2}^2+\|\nabla^2d\|_{H^2}^2\right)
+\int_0^T\|\na^2d\|_{H^2}^2dt\nonumber\\
&&\quad +\int_0^T\sigma \left(\|\rho^{1/2}u_{tt}\|_{L^2}^2+\|\nabla
u_t\|_{H^1}^2+\|\nabla d_{tt}\|_{L^2}^2\right)dt\leq
C(T).\label{4.20}
\end{eqnarray}
\end{lem}
\pf Differentiating (\ref{1.2}) and (\ref{1.3}) with respect to $t$
gives
\begin{equation}
\rho u_{tt}-\mu\Delta u_t-(\mu+\lambda)\nabla{\rm div}u_t= -\rho_t
u_t-(\rho u\cdot\nabla u)_t-\nabla P_t-\divg M(d)_t\label{4.21}
\end{equation}
and
\begin{equation}
d_{tt}-\Delta d_t=(|\nabla d|^2d)_t-(u\cdot\nabla d)_t.\label{4.22}
\end{equation}
Thus, multiplying (\ref{4.21}) and (\ref{4.22}) by $u_{tt}$ and
$-\Delta d_{tt}$ respectively, and integrating the resulting
equations over $\r^3$, we obtain after adding them together that
\begin{eqnarray}
&&\frac{1}{2}\frac{d}{dt}\int \left(\mu|\nabla
u_t|^2+(\mu+\lambda)(\divg
u_t)^2+|\Delta d_{t}|^2\right)dx+\int\left( \rho|u_{tt}|^2+|\nabla d_{tt}|^2\right)dx\nonumber\\
&&\quad=-\int \left(\rho_t u_t+\rho_tu\cdot\nabla u+\rho
u_t\cdot\nabla u+\rho u\cdot\nabla u_t+\nabla P_t+\divg M(d)_t\right)\cdot u_{tt}dx\nonumber\\
&&\qquad+\int\left(u_t\cdot\nabla d+u\cdot\nabla
d_t-|\nabla d|^2d_t-2d(\nabla d)\cdot(\nabla d_t) \right)\cdot\Delta d_{tt}dx\nonumber\\
&&\quad=-\frac{d}{dt}\int \left[\frac{1}{2}\rho_t |u_t|^2+\left(
\rho_tu\cdot\nabla u+\nabla P_t+\divg M(d)_t\right)\cdot u_{t}\right]dx\nonumber\\
&&\qquad+\int \left(\rho_{tt}u\cdot\nabla u+\rho_tu_t\cdot\nabla
u+\rho_tu\cdot\nabla u_t+\divg M(d)_{tt} \right)\cdot u_{t}dx\nonumber\\
&&\qquad+\frac{1}{2}\int \rho_{tt}|u_t|^2dx-\int P_{tt}\divg u_t
dx-\int \left(\rho u_t\cdot\nabla u+\rho u\cdot\nabla
u_t\right)\cdot u_{tt}dx\nonumber\\
&&\qquad+\int\left[\left(u_t\cdot\nabla d+u\cdot\nabla d_t-|\nabla
d|^2d_t-2d(\nabla d)\cdot(\nabla d_t)\right)\cdot\Delta d_{tt} \right]dx\nonumber\\
&&\quad\triangleq\frac{d}{dt}I_0+\sum_{i=1}^5I_i.\label{4.23}
\end{eqnarray}

We now estimate each term on the right-hand side of (\ref{4.23}).
Using (\ref{1.1}) and integrating by parts, we first deduce from
Lemmas \ref{lem2.1}, \ref{lem4.1}--\ref{lem4.3} and Cauchy-Schwarz
inequality that
\begin{eqnarray*}
I_0&=&-\int\left[\rho u\cdot (u_t\cdot \nabla
u_t)+\rho_tu\cdot\nabla u \cdot u_t-P_t\divg u_t-M(d)_t:\nabla u_{t} \right]dx\nonumber\\
&\leq&C\|u\|_{L^\infty}\|\rho^{1/2} u_t\|_{L^2}\|\nabla
u_t\|_{L^2}+C\|\rho_t\|_{L^2}\|u\|_{L^\infty}\|\nabla
u\|_{L^3}\|u_t\|_{L^6}\\
&&+C\|P_t\|_{L^2}\|\nabla u_t\|_{L^2}+C\|\nabla
d\|_{L^\infty}\|\nabla
d_t\|_{L^2}\|\nabla u_t\|_{L^2}\nonumber\\
&\leq&\frac{\mu}{4}\|\nabla u_t\|_{L^2}^2+C.
\end{eqnarray*}

Since it holds that $M(d)_{tt}\leq C(|\nabla d_{tt}||\nabla
d|+|\nabla d_t|^2)$, we have by Lemmas \ref{lem2.1},
\ref{lem4.1}--\ref{lem4.3} and integration by parts that
\begin{eqnarray*}
I_1&\leq& C\|\rho_{tt}\|_{L^2}\|u\|_{L^\infty}\|\nabla
u\|_{L^3}\|u_t\|_{L^6}+C\|\rho_t\|_{L^2}\|u_t\|_{L^6}^2\|\nabla
u\|_{L^6}\nonumber\\
&&+C\|\rho_t\|_{L^6}\|u\|_{L^6}\|\nabla
u_t\|_{L^2}\|u_t\|_{L^6}+C\|\nabla d_{tt}\|_{L^2}\|\nabla
d\|_{L^\infty}\|
\nabla u_t\|_{L^2}\nonumber\\
&&+C\|\nabla d_t\|_{L^2}^{1/2}\|\nabla d_t\|_{L^6}^{3/2} \|\nabla u_t\|_{L^2}\nonumber\\
&\leq&\frac{1}{4}\|\nabla d_{tt}\|_{L^2}^2+ C
\left(\|\rho_{tt}\|_{L^2}^2+\|\nabla
u_t\|_{L^2}^2+\|\nabla^2d_t\|_{L^2}^3\right).
\end{eqnarray*}

Due to (\ref{1.1}), one has $\rho_{tt}=-\divg(\rho u)_t$, and hence,
it follows from from (\ref{4.12}), (\ref{4.14}) and  integration by
parts that
\begin{eqnarray*}
I_2+I_3&=&\int \left(\rho_t u\cdot\nabla u_t\cdot u_t+\rho
u_t\cdot\nabla u_t\cdot u_t \right)dx-\int P_{tt}\divg u_t dx\nonumber\\
&\leq&C\left(\|\rho_t\|_{L^3}\|u\|_{L^\infty}\|u_t\|_{L^6}+\|\rho
u_t\|_{L^3} \|u_t\|_{L^6}+\|P_{tt}\|_{L^2}\right)\|\nabla u_t\|_{L^2}\nonumber\\
&\leq&C \left(\|\nabla u_t\|_{L^2}^2+\|\rho^{1/2}
u_t\|_{L^2}^{1/2}\|u_t\|_{L^6}^{1/2}\|\nabla u_t\|_{L^2}^2+\|P_{tt}\|_{L^2}^2\right)\\
&\leq&C \left(1+\|\nabla u_t\|_{L^2}^4+\|P_{tt}\|_{L^2}^2\right).
\end{eqnarray*}
and similarly, by (\ref{4.4}) one gets
\begin{eqnarray*}
I_4&\leq&C\left(\|u_t\|_{L^6}\|\nabla
u\|_{L^3}+\|u\|_{L^\infty}\|\nabla u_t\|_{L^2}\right)\|\rho^{1/2}
u_{tt}\|_{L^2}\\
&\leq& \frac{1}{2}\|\rho^{1/2} u_{tt}\|_{L^2}^2+C\|\nabla
u_t\|_{L^2}^2.
\end{eqnarray*}

Finally, integrating by parts and using Lemmas \ref{lem2.1} and
\ref{lem4.1}--\ref{lem4.3}, we can estimate the last term on the
right-hand side of (\ref{4.23}) as follows:
\begin{eqnarray*}
I_5 &\leq&C\int \left(|\nabla u_t||\nabla
d|+|u_t||\nabla^2d|+|\nabla u||\nabla d_t|+|u||\nabla^2d_t|
\right)|\nabla
d_{tt}| dx\nonumber\\
&&+C\int \left( |\nabla d||\nabla^2d||d_t|+|\nabla d|^2|\nabla
d_t|+|\nabla^2d||\nabla d_t|+|\nabla d||\nabla^2d_t|\right)|\nabla
d_{tt}| dx\nonumber\\
&\leq&C \|\nabla d\|_{L^\infty}\|\nabla u_t\|_{L^2}\|\nabla
d_{tt}\|_{L^2}+C\|u_t\|_{L^6}\|\nabla^2d\|_{L^3}\|\nabla
d_{tt}\|_{L^2}\\
&&+C\|\nabla u\|_{L^3}\|\nabla d_t\|_{L^6}\|\nabla
d_{tt}\|_{L^2}+C\| u\|_{L^\infty}\|\nabla^2
d_t\|_{L^2} \|\nabla d_{tt}\|_{L^2}\\
&&+C\|\nabla d\|_{L^\infty}\|\nabla^2d\|_{L^3}\|d_t\|_{L^6}\|\nabla
d_{tt}\|_{L^2}+C \|\nabla d\|_{L^\infty}^2 \|\nabla
d_t\|_{L^2}\|\nabla
d_{tt}\|_{L^2}\\
&&+C\|\nabla^2d\|_{L^3}\|\nabla d_t\|_{L^6}\|\nabla
d_{tt}\|_{L^2}+C\|\nabla d\|_{L^\infty}\|\nabla^2d_t\|_{L^2}
\|\nabla
d_{tt}\|_{L^2}\nonumber\\
&\leq&\frac{1}{4}\|\nabla d_{tt}\|_{L^2}^2+C \left(1+\|\nabla
u_t\|_{L^2}^2+\|\nabla^2d_t\|_{L^2}^2\right).
\end{eqnarray*}

By virtue of the estimates of $I_i$ ($i=0,1,\ldots,5$), we deduce
after multiplying (\ref{4.23}) by $\sigma(t)$, integrating it over
$[0,T]$ and using Gronwall's inequality that
\begin{equation}
\sup_{0\leq t\leq T}\sigma\left(\|\nabla u_t\|_{L^2}^2+\|\nabla^2
d_{t}\|_{L^2}^2\right)+\int_{0}^T\sigma\left(\|
\rho^{1/2}u_{tt}\|_{L^2}^2+\|\nabla d_{tt}\|_{L^2}^2\right) dt\leq
C,\label{4.24}
\end{equation}
where we have also used (\ref{4.2}), (\ref{4.12}) and (\ref{4.14}).
This, together with (\ref{4.13}) and (\ref{4.17}), gives
\begin{equation}
\sup_{0\leq t\leq T}\left(\sigma\|\nabla u\|_{H^2}^2\right)\leq
\sup_{0\leq t\leq T}\sigma\left(1+\|\nabla u_t\|_{L^2}^2+\|\nabla
P\|_{H^1}^2\right)\leq C.\label{4.25}
\end{equation}

On the other hand, by Lemmas \ref{lem4.1}--\ref{lem4.3} we obtain by
applying the standard $L^2$-estimate to (\ref{4.21}) that
\begin{eqnarray}
\|\nabla^2 u_t\|_{L^2}&\leq& C\|\mu\Delta
u_t+(\mu+\lambda)\nabla\divg
u_t\|_{L^2}\nonumber\\
&\leq&C\|\left(\n u_{t} +\rho u\cdot\nabla u+\nabla P+\divg M(d)\right)_t \|_{L^2}\nonumber\\
&\leq&C \|\rho^{1/2}
u_{tt}\|_{L^2}+C\|\rho_t\|_{L^3}\|u_t\|_{L^6} +C\|\rho_t\|_{L^3}\|u\|_{L^\infty}\|\nabla u\|_{L^6}\nonumber\\
&&  +C\|u_t\|_{L^6}\|\nabla u\|_{L^3}+C\|u\|_{L^\infty}\|\nabla
u_t\|_{L^2}+C\|\nabla P_t\|_{L^2}\nonumber\\
&&+C\|\nabla d\|_{L^\infty}\|\nabla^2 d_t\|_{L^2}+C\|\nabla d_t\|_{L^6}\|\nabla^2d\|_{L^3} \nonumber\\
&\leq& C\left(1+\|\nabla u_t\|_{L^2}+\|\rho^{1/2}
u_{tt}\|_{L^2}+\|\nabla^2d_t\|_{L^2}\right),\label{4.26}
\end{eqnarray}
which, together with (\ref{4.24}), yields
\begin{equation}
\int_0^T\sigma\|\nabla^2 u_t\|_{L^2}^2dt\leq C.\label{4.27}
\end{equation}

In a similar manner as the derivation of (\ref{4.25}), we also infer
from (\ref{1.3}) that
\begin{eqnarray*}
\|\nabla^4d\|_{L^2}&\leq&
C\left(\|\nabla^2d_t\|_{L^2}+\|\nabla^2(u\cdot\nabla
d)\|_{L^2}+\|\nabla^2(|\nabla d|^2d)\|_{L^2}\right)\nonumber\\
&\leq&C \left(1+\|\nabla^2d_t\|_{L^2}+\|\nabla u\|_{H^1}\|\nabla
d\|_{H^2}+\|\nabla d\|_{H^2}^4\right)\nonumber\\
&\leq&C+C\|\nabla^2d_t\|_{L^2},
\end{eqnarray*}
from which,   (\ref{4.24}), and (\ref{4.2}), it follows that
\begin{equation}
\sup_{0\leq t\leq T}\left(\sigma\|\nabla^2d\|_{H^2}^2\right)+\int_0^T\|\na^2d\|_{H^2}^2dt\leq
C.\label{4.28}
\end{equation}
Therefore, combining (\ref{4.24}), (\ref{4.25}), (\ref{4.27}) and
(\ref{4.28}) finishes the proof of (\ref{4.20}).\hfill$\square$

\vskip 2mm

To prove the H\"{o}lder continuity of the first-order derivatives of
density and pressure, we need the following lemma which is concerned
with the $W^{1,q}$-estimate ($q\in(3,6)$) on the gradients of
density and pressure.
\begin{lem}\label{lem4.5} For fixed $q\in(3,6)$, it holds for any $T>0$ that
\begin{equation}
\sup_{0\leq t\leq T}\left(\|\nabla \rho\|_{W^{1,q}}+\|\nabla
P\|_{W^{1,q}}\right)+\int_0^T\left(\|\nabla
u_t\|_{L^{q}}^{p_0}+\|\nabla^2 u\|_{W^{1,q}}^{p_0}\right)dt\leq
C(T),\label{4.29}
\end{equation}
where
\begin{equation}
\la{4.30}1\leq p_0<\frac{4q}{5q-6}\in(1,2).
\end{equation}
\end{lem}
\pf Operating $\na^2$ to both sides of (\ref{1.1}), (\ref{4.15}),
and multiplying them by  $q|\na^2 \rho|^{q-2}\na^2\rho$, $q|\na^2
P(\rho)|^{q-2}\na^2 P(\rho)$, respectively, we obtain after
integrating by parts over $\R^3$ and using (\ref{4.4}), (\ref{4.13})
and Lemma \ref{lem2.1} that for any $q\in(3,6)$,
\begin{eqnarray}
&&\frac{d}{dt}\left(\|\nabla^2\rho\|_{L^q}^q+\|\nabla^2
P\|_{L^q}^q\right)\nonumber\\
&&\quad\leq C\|\nabla u\|_{L^\infty}\left(\|\nabla^2
\rho\|_{L^q}^q+\|\nabla^2
P\|_{L^q}^q\right)+C\|\nabla^2u\|_{W^{1,q}}\left(\|\nabla^2
\rho\|_{L^q}^{q-1}+\|\nabla^2
P\|_{L^q}^{q-1}\right)\nonumber\\
&&\qquad+C\|\nabla^2 u\|_{L^q}\left(\|\nabla
\rho\|_{L^\infty}\|\nabla^2 \rho\|_{L^q}^{q-1}+\|\nabla
P\|_{L^\infty}\|\nabla^2
P\|_{L^q}^{q-1}\right)\nonumber\\
&&\quad\leq C\left(1+\|\nabla
u\|_{H^2}\right)\left(1+\|\nabla^2\rho\|_{L^q}^q+\|\nabla^2P\|_{L^q}^q\right)\nonumber\\
&&\quad+C\|\nabla^2u\|_{W^{1,q}}
\left(\|\nabla^2\rho\|_{L^q}^{q-1}+\|\nabla^2P\|_{L^q}^{q-1}\right).\label{4.31}
\end{eqnarray}

Applying the standard $L^{p}$-estimate to the elliptic system
(\ref{1.2}) yields that
\begin{eqnarray}
\|\nabla^2 u\|_{W^{1,q}}&\leq& C \|\nabla^2u\|_{L^q}+C\|\nabla(\rho
u_t+\rho u\cdot \nabla u+\nabla P+\divg M(d))\|_{L^{q}} \nonumber\\
&\leq&C \|\nabla
u\|_{H^2}+C\left(\|\nabla\rho\|_{L^q}\|u_t\|_{L^\infty}+\|\nabla
u_t\|_{L^{q}}\right)\nonumber\\
&& +C\left(\|u\|_{L^\infty}\|\nabla^2u\|_{L^q}+\|\nabla
u\|_{L^\infty}\|\nabla
u\|_{L^q}+\|\nabla\rho\|_{L^q}\|u\|_{L^\infty}\|\nabla
u\|_{L^\infty}\right)\nonumber\\
&&+C\|\nabla^2
P\|_{L^{q}}+C\left(\|\nabla d\|_{L^{\infty}}\|\nabla^3 d\|_{L^{q}}+\|\nabla^2d\|_{L^\infty}\|\nabla^2 d\|_{L^q}\right)\nonumber\\
&\leq& C\left(1+\|\nabla u_t\|_{L^2}+\|\nabla u_t\|_{L^q}+\|\nabla
u\|_{H^2}+\|\nabla^2d\|_{H^2}+\|\nabla^2
P\|_{L^{q}}\right),\label{4.32}
\end{eqnarray}
where we have also used Lemmas \ref{lem2.1} and
\ref{lem4.1}--\ref{lem4.3}. Putting (\ref{4.32}) into (\ref{4.31})
gives
\begin{eqnarray}
\frac{d}{dt}\left(\|\nabla^2 \rho\|_{L^q}^q+\|\nabla^2
P\|_{L^q}^q\right)&\leq& C\left(1+\|\nabla u_t\|_{L^2}+\|\nabla
u_t\|_{L^q}+\|\nabla
u\|_{H^2}+\|\nabla^2d\|_{H^2}\right)\nonumber\\
&&\times\left(1+\|\nabla^2\rho\|_{L^q}^q+\|\nabla^2P\|_{L^q}^q\right).\label{4.33}
\end{eqnarray}

For  $p_0$ being the one in (\ref{4.30}), it follows from
(\ref{4.12}) and (\ref{4.20}) that
\begin{eqnarray}
&&\int_0^T\left(\|\nabla u_t\|_{L^2}+\|\nabla
u\|_{H^2}+\|\nabla^2d\|_{H^2}\right)^{p_0}dt\nonumber\\
&&\quad\leq C+C\sup_{0\leq t\leq T}\left(\sigma\|\nabla
u\|_{H^2}^2+\sigma\|\nabla^2d\|_{H^2}^2\right)^{p_0/2}\int_0^T\sigma^{-p_0/2}dt\leq
C,\label{4.34}
\end{eqnarray}
and moreover, using Lemma \ref{lem2.1}, H\"{o}lder inequality and
(\ref{4.20}), we find that
\begin{eqnarray}
&&\int_0^T\|\nabla u_t\|_{L^q}^{p_0}dt\leq C\int_0^T\|\nabla
u_t\|_{L^2}^{p_0(6-q)/(2q)}\|\nabla u_t\|_{L^6}^{p_0(3q-6)/(2q)}dt\nonumber\\
&&\quad\leq C\int_0^T\sigma^{-p_0/2}\left(\sigma\|\nabla
u_t\|_{L^2}^2\right)^{p_0(6-q)/(4q)}\left(\sigma\|\nabla
u_t\|_{H^1}^2\right)^{p_0(3q-6)/(4q)}dt\nonumber\\
&&\quad\leq C\left(\sup_{0\leq t\leq T}\sigma\|\nabla
u_t\|_{L^2}^2\right)^{p_0(6-q)/(4q)}\int_0^T\sigma^{-p_0/2}\left(\sigma\|\nabla
u_t\|_{H^1}^2\right)^{p_0(3q-6)/(4q)}dt\nonumber\\
&&\quad\leq C\left(\int_0^T
\sigma^{-2p_0q/(4q-p_0(3q-6))}dt\right)^{(4q-p_0(3q-6))/(4q)}\left(\int_0^T\sigma\|\nabla
u_t\|_{H^1}^2dt\right)^{p_0(3q-6)/(4q)}\nonumber\\
&&\quad\leq C,\label{4.35}
\end{eqnarray}
since $0<2p_0q/(4q-p_0(3q-6))<1$ and $0<p_0(3q-6)/(4q)<1$.

By virtue of (\ref{4.34}) and (\ref{4.35}), we deduce from
(\ref{4.33}) and Gronwall inequality that
$$
\sup_{0\leq t\leq T}\left(\|\nabla\rho\|_{W^{1,q}}+\|\nabla
P\|_{W^{1,q}}\right)\leq C,\quad\forall\;q\in(3,6),
$$
which, together with (\ref{4.32}), (\ref{4.34}) and (\ref{4.35}),
yields the desired estimate of (\ref{4.29}).\hfill$\square$

\vskip 2mm

Finally, we need the following initial-layer analysis which
particularly implies $(u_t,\nabla^2 u)$ are H\"{o}lder continuous
away from the initial time $t=0$.
\begin{lem}\label{lem4.6}For any given $T>0$, it holds that
\begin{equation}
\sup_{0\leq t\leq T}\sigma\left(\|\rho^{1/2}
u_{tt}\|_{L^2}+\|\nabla^2 u_t\|_{L^2}+\|\nabla^2
u\|_{W^{1,q}}\right)+\int_0^T\sigma^2\|\nabla u_{tt}\|_{L^2}^2dt\leq
C(T).\label{4.36}
\end{equation}
\end{lem}
\pf Differentiating (\ref{4.21})   with respect to $t$ gives
\begin{eqnarray*}
&&\rho u_{ttt}+\rho u\cdot\nabla u_{tt}-\mu\Delta
u_{tt}-(\mu+\lambda)\nabla\divg u_{tt}\\
&&\quad=2\divg (\rho u)u_{tt}+\divg(\rho u)_tu_t-2(\rho
u)_t\cdot\nabla u_t-(\rho_{tt}u+2\rho_tu_t)\cdot\nabla u\\
&&\qquad-\rho u_{tt}\cdot\nabla u-\nabla P_{tt}-\divg M(d)_{tt},
\end{eqnarray*}
which, multiplied by $u_{tt}$ in $L^2$ and integrated by parts over
$\R^3$, yields
\begin{eqnarray}
&&\frac{1}{2}\frac{d}{dt}\int \rho|u_{tt}|^2dx+\int \left(\mu|\nabla
u_{tt}|^2+(\mu+\lambda)(\divg u_{tt})^2\right)dx\nonumber\\
&&\quad=-4\int \rho u\cdot\nabla u_{tt}\cdot u_{tt}dx-\int (\rho
u)_t\cdot\left(\nabla(u_t\cdot
u_{tt})+2\nabla u_t\cdot u_{tt}\right)dx\nonumber\\
&&\qquad-\int \left(\rho_{tt} u+2\rho_tu_t\right)\cdot \nabla u\cdot
u_{tt}dx-\int \rho u_{tt}\cdot\nabla u\cdot
u_{tt}dx\nonumber\\
&&\qquad+\int P_{tt}\divg u_{tt}dx+\int M(d)_{tt}:\nabla
u_{tt}dx\triangleq\sum_{i=1}^6J_i.\label{4.37}
\end{eqnarray}

The right-hand side of (\ref{4.37}) will be estimated term by term
as follows, using Lemma \ref{lem2.1}, \ref{lem4.1}--\ref{lem4.5} and
Cauchy-Schwarz inequality as well.
\begin{eqnarray*}
J_1&\leq&\|u\|_{L^\infty}\|\rho^{1/2} u_{tt}\|_{L^2}\|\nabla
u_{tt}\|_{L^2}\leq \delta \|\nabla
u_{tt}\|_{L^2}^2+C(\delta)\|\rho^{1/2} u_{tt}\|_{L^2}^2,\\[2mm]
J_2&\leq&C\left(\|\rho
u_{t}\|_{L^3}+\|u\|_{L^\infty}\|\rho_t\|_{L^3}\right)\left(\|\nabla
u_t\|_{L^2}\|u_{tt}\|_{L^6}+\|u_t\|_{L^6}\|\nabla
u_{tt}\|_{L^2}\right)\\
&\leq&C\left(1+\|\rho^{1/2} u_t\|_{L^2}^{1/2}\|\nabla
u_t\|_{L^2}^{1/2}\right)\|\nabla u_{t}\|_{L^2}\|\nabla
u_{tt}\|_{L^2}\\
&\leq& \delta\|\nabla u_{tt}\|_{L^2}^2+C(\delta)\left(1+\|\nabla
u_t\|_{L^2}^{3}\right),\\[2mm]
J_3&\leq&
C\left(\|\rho_{tt}\|_{L^2}\|u\|_{L^\infty}+\|\rho_t\|_{L^3}\|u_t\|_{L^6}
\right)\|\nabla
u\|_{L^3}\|u_{tt}\|_{L^6}\\
&\leq& \delta\|\nabla u_{tt}\|_{L^2}^2
+C(\delta)\left(\|\rho_{tt}\|_{L^2}^2+\|\nabla
u_t\|_{L^2}^2\right),\\[2mm]
J_4&\leq&C\|\rho^{1/2} u_{tt}\|_{L^2}\|\nabla
u\|_{L^3}\|u_{tt}\|_{L^6}\leq \delta\|\nabla
u_{tt}\|_{L^2}^2+C(\delta)\|\rho^{1/2} u_{tt}\|_{L^2}^2,\\[2mm]
J_5&\leq& C \|P_{tt}\|_{L^2} \|\nabla u_{tt}\|_{L^2}
\leq\delta\|\nabla u_{tt}\|_{L^2}^2+C(\delta) \|P_{tt}\|_{L^2}^2
\end{eqnarray*}
and
\begin{eqnarray*}
J_6&\leq&C\left(\|\nabla d\|_{L^\infty}\|\nabla
d_{tt}\|_{L^2}+\|\nabla d_t\|_{L^2}^{1/2}\|\nabla
d_t\|_{L^6}^{3/2}\right)\|\nabla u_{tt}\|_{L^2}\\
&\leq&\delta\|\nabla u_{tt}\|_{L^2}^2+C(\delta)\left(\|\nabla
d_{tt}\|_{L^2}^2+\|\nabla^2d_t\|_{L^2}^3\right),
\end{eqnarray*}
since direct computations give
$$
|M(d)_{tt}|\leq C\left(|\nabla d||\nabla d_{tt}|+|\nabla
d_t|^2\right).
$$

Thus, putting the estimates of $J_1,\ldots,J_6$ into (\ref{4.37})
and multiplying the resulting inequality by $\sigma^2$, we have by
choosing $\delta>0$ small enough and using Lemmas
\ref{lem4.1}--\ref{lem4.4} that
\begin{eqnarray}
&&\sup_{0\leq t\leq
T}\left(\sigma^2\|\rho^{1/2}u_{tt}\|_{L^2}^2\right)+\int_0^T\sigma^2\|\nabla
u_{tt}\|_{L^2}^2dt\nonumber\\
&&\quad\leq C+C\int_0^T \sigma\|\rho^{1/2}
u_{tt}\|_{L^2}^2dt+C\int_0^T\left(\|P_{tt}\|_{L^2}^2+\|\rho_{tt}\|_{L^2}^2+\sigma\|\nabla
d_{tt}\|_{L^2}^2\right)dt\nonumber\\
&&\qquad+C\sup_{0\leq t\leq T}\sigma^{1/2}\left(\|\nabla
u_t\|_{L^2}+\|\nabla^2d_t\|_{L^2}\right) \int_0^T\left(\|\nabla
u_t\|_{L^2}^2+\|\nabla^2d_t\|_{L^2}^2\right)dt\nonumber\\
&&\quad\leq C.\label{4.38}
\end{eqnarray}

As a result of (\ref{4.38}), (\ref{4.26}) and Lemma \ref{lem4.4}, we
also see that
\begin{equation}
\sigma\|\nabla^2 u_t\|_{L^2} \leq C\sigma\left(1+\|\nabla
u_t\|_{L^2}+\|\rho^{1/2}
u_{tt}\|_{L^2}+\|\nabla^2d_t\|_{L^2}\right)\leq C,\label{4.39}
\end{equation}
and thus, it follows from (\ref{4.32}), (\ref{4.39}), Lemmas
\ref{lem4.4} and \ref{lem4.5} that for any $q\in(3,6)$,
\begin{eqnarray}
\sigma\|\nabla^2 u\|_{W^{1,q}}&\leq& C\sigma\left(1+\|\nabla
u_t\|_{L^2}+\|\nabla u_t\|_{L^q}+\|\nabla
u\|_{H^2}+\|\nabla^2d\|_{H^2}+\|\nabla^2
P\|_{L^{q}}\right)\nonumber\\
&\leq&C\sigma\left(1+\|\nabla u_t\|_{H^1}+\|\nabla
u\|_{H^2}+\|\nabla^2d\|_{H^2}+\|\nabla^2
P\|_{L^{q}}\right)\nonumber\\
&\leq& C.\label{4.40}
\end{eqnarray}
Therefore, combining (\ref{4.38})--(\ref{4.40}) immediately leads to
(\ref{4.36}).\hfill$\square$

\section{Proof of Theorem \ref{thm1.1}}
With all the a priori estimates at hand, we are now ready to prove
our main results. To this end, we first need the following local
existence theorem of classical solutions of (\ref{1.1})--(\ref{1.7})
with large initial data.

\begin{pro}\label{pro5.1}Assume that the initial data
$(\rho_0,u_0,d_0)$ satisfy the conditions (\ref{1.11}), (\ref{1.12})
of Theorem \ref{thm1.1}. Then there exist a positive time $T_0>0$
and a unique classical solution $(\rho,u,d)$ of
(\ref{1.1})--(\ref{1.7}) on $\R^3\times(0,T_0]$, satisfying
$\rho\geq0$, $|d|=1$, and for any $\tau\in(0,T_0)$,
\begin{equation}
\left\{
\begin{array}{lll}
\left(\rho-\tilde\rho,P(\rho)-P(\tilde\rho)\right)\in
C([0,T_0];H^1\cap W^{1,q}),
\\[2mm]
\left(\rho-\tilde\rho,P(\rho)-P(\tilde\rho)\right)\in L^\infty(0,T_0;H^2\cap W^{2,q})\\[2mm]
u\in C([0,T];D^1\cap D^2)\cap L^\infty(\tau,T_0; D^3\cap
D^{3,q}),\\[2mm]
u_t\in L^\infty(\tau,T_0;D^1\cap D^2)\cap H^1(\tau,T_0;D^1),
\\[2mm]
\nabla d\in C([0,T_0]; H^2)\cap L^\infty(\tau,T_0; H^3),\\[2mm]
d_t\in C([0,T_0]; H^1)\cap L^\infty(\tau,T_0; H^2).
\end{array}\right.\label{5.1}
\end{equation}
\end{pro}
\pf As that in \cite{HWW2011}, we can use the Galerkin's
approximation method to construct the approximate solutions $u^m$ to
the momentum equation, then use this approximate $u^m$ and the
equations of conservation mass and angular momentum to get
$\rho^m,d^m$. The existence of a smooth approximate solution
$(\rho^m,u^m,d^m)$ follows from the fixed point theorem, similar to
that on the compressible Navier-Stokes equations (see, for example,
\cite{Pa1986,HWW2011}). Now, in order to prove the convergence for
the approximate solutions and to obtain a smooth solution of
(\ref{1.1})--(\ref{1.7}), it is essential to derive some uniform a
priori estimate for $(\rho^m,u^m,d^m)$.

It has been shown in \cite[Theorem 2.1 or (3.21),
(\ref{3.27})]{HWW2011} that there exists a small $T_0>0$, independent
of $m$ and the lower bound of density, such that
\begin{eqnarray}
&&\sup_{0\leq t\leq T_0}\left(\|\sqrt{\rho^m}
u^m_t\|_{L^2}^2+\|\rho^m-\tilde\rho\|_{H^1\cap W^{1,q}}^2+\|\nabla
u^m\|_{H^1}+\|d^m_t\|_{H^1}^2+\|\nabla
d^m\|_{H^2}^2\right)\nonumber\\
&&\quad +\int_0^{T_0}\left(\|u^m\|_{D^{2,q}}^2+\|\nabla
u^m_t\|_{L^2}^2+\|\nabla^4d^m\|_{L^2}^2+\|\nabla^2d^m_t\|_{L^2}^2\right)\leq
\tilde C\exp\left(\tilde C\|g\|_{L^2}^2\right),\label{5.2}
\end{eqnarray}
where $\tilde C>0$ may depend on $T_0$, but is independent of the
size of domain.

In view of the bounds in (\ref{5.2}), we can proceed to derive more
higher-order estimates on the solutions $(\rho^m,u^m,d^m)$ in the
same way as those carried out in Lemmas \ref{lem4.3}--\ref{lem4.6}.
To summarize up, for any $\tau\in(0,T_0]$ the approximate solutions
$(\rho^m,u^m,d^m)$ satisfy
\begin{eqnarray}\la{5.3}
&&\sup_{0\leq t\leq T_0}\left(\|\nabla \rho^m\|_{W^{1,q}}+\|\nabla
P^m\|_{W^{1,q}}\right)+\int_0^{T_0}\left(\|\nabla
u^m_t\|_{L^{q}}^{p_0}+\|\nabla^2 u^m\|_{W^{1,q}}^{p_0}\right)dt\nonumber\\
&&\quad+\sup_{\tau\leq t\leq
T_0}\left(\|\sqrt{\rho^m}u^m_{tt}\|_{L^2}^2+\|\nabla^2
u^m\|_{H^1\cap W^{1,q}}^2+\|\nabla
u^m_t\|_{H^1}^2\right)+\int_{\tau}^{T_0}\|\nabla
u^m_{tt}\|_{L^2}^2dt\nonumber\\
&&\quad+\sup_{\tau\leq t\leq T_0}
\left(\|\nabla^2d^m\|_{H^2}^2+\|\nabla^2d^m_t\|_{L^2}^2\right)+\int_{\tau}^{T_0}
\|\nabla d^m_{tt}\|_{L^2}^2 dt\leq \tilde C
\end{eqnarray}
for some positive constant $\tilde C$ which may depend on $T_0$, but
is independent of $m$ and the lower bound of density.

With the help of the local estimates (\ref{5.2}) and (\ref{5.3}), we
easily deduce after taking a subsequence
$(\rho^{m_j},u^{m_j},d^{m_j})$ and passing to limit as $j\to\infty$
that $(\rho^{m_j},u^{m_j},d^{m_j})$ would converges to a solution
$(\rho, u,d)$ of (\ref{1.1})--(\ref{1.7}) on $\R^3\times(0,T_0)$
satisfying (\ref{5.1}) due to the lower semi-continuity. The
uniqueness of strong/classical solutions can be proved in the same
manner as that in \cite{HW2011}. This finishes the proof of
Proposition \ref{pro5.1}.\hfill$\square$
 \vskip 2mm

{\it Proof of Theorem \ref{thm1.1}.} By Proposition \ref{pro5.1},
there exists a small time $T_0>0$ such that the Cauchy problem
(\ref{1.1})--(\ref{1.7}) has a unique classical solution $(\rho,
u,d)$ on $\R^3\times (0,T_0]$. We shall make use of the a priori
estimates, Proposition \ref{pro3.1} and Lemmas
\ref{lem4.1}--\ref{lem4.6}, to extend the local classical solution
$(\rho,u,d)$ to all time.

First, in view of the definitions in (\ref{3.1})--(\ref{3.5}), it is
easily seen from (\ref{3.20}) and (\ref{3.88}) that
$$
A_1(0)+A_2(0)=0,\quad A_3(0)\leq C_0^{\delta_0},\quad
A_4(0)+A_5(0)\leq C_0^{\delta_0}\quad{\rm and}\quad \rho_0\leq
\bar\rho,
$$
due to $C_0\leq\varepsilon$. Thus, there exists a $T_1\in(0,T_0]$
such that (\ref{3.7}) hold for $T=T_1$.

To be continued, we set
\begin{equation}
T^*=\sup\left\{T\;\left|\;(\ref{3.7})\;\; {\rm
holds}\right.\right\}.\label{5.4}
\end{equation}
Then $T^*\geq T_1>0$. Hence, for any $0<\tau<T\leq T^*$, it follows
from (\ref{1.1}), (\ref{4.4}), (\ref{4.12}), (\ref{4.13}),
(\ref{4.20}), (\ref{4.36}) and Lemma \ref{lem2.1} that
\begin{eqnarray*}
&&\int_\tau^T\int\left|\partial_t\left(\rho |u_t|^2\right)
\right|dxdt+\int_\tau^T\int\left|\partial_t\left(\rho |u\nabla
u|^2\right) \right| dx
dt\\
&&\quad\leq C\int_\tau^T\int\left(|\rho_t||u_t|^2+\rho |u_t||
u_{tt}|\right)dxdt\\
&&\qquad+C\int_\tau^T\int\left(|\rho_t||u|^2|\nabla u|^2+\rho|u||u_t||\nabla u|^2+\rho|u|^2|\nabla u||\nabla u_t|\right)dxdt\\
&&\quad\leq C\int_\tau^T\int\left(\rho|\divg u||u_t|^2+|u||\nabla
\rho||u_t|^2+\rho |u_t||u_{tt}|\right)dxdt\\
&&\qquad+C\int_\tau^T\int\left( |\divg u||u|^2|\nabla u|^2+|\nabla \rho||u|^3|\nabla u|^2+ |u|^2|\nabla u||\nabla u_t|+\rho|u||u_t||\nabla u|^2\right)dxdt\\
&&\quad \leq C\int_\tau^T\left(\|\nabla u\|_{L^\infty}\|\rho^{1/2}
u_t\|_{L^2}^2+\|u\|_{L^6}\|\nabla
\rho\|_{L^2}\|u_t\|_{L^6}^2+\|\rho^{1/2}u_t\|_{L^2}\|\rho^{1/2}u_{tt}\|_{L^2}\right)dt\\
&&\qquad+C\int_\tau^T\left(\|\nabla u\|_{H^2}^5 +\|\nabla
u\|_{H^2}^5\|\nabla \rho\|_{L^2} +\|\nabla u\|_{H^2}^3\|\nabla
u_t\|_{L^2}+\|\rho^{1/2} u_t\|_{L^2}\|\nabla
u\|_{H^2}^3\right)dt\\
&&\quad\leq C(\tau,T),
\end{eqnarray*}
which, together with (\ref{4.4}) and (\ref{4.12}), yields
$$
\rho^{1/2} u_t,\quad \rho^{1/2} u\cdot\nabla u\in C([\tau, T];L^2).
$$
and consequently,
\begin{equation}
\label{5.5} \rho^{1/2} \dot u\in  C ([\tau, T];L^2).
\end{equation}

Next, we claim that
\begin{equation}
T^*=\infty.\label{5.6}
\end{equation}
Otherwise, $T^*<\infty$. Then by Proposition \ref{pro3.1},
(\ref{3.9}) holds for $T= T^*$. So, it follows from Lemmas
\ref{lem4.1}--\ref{lem4.6} and (\ref{5.5}) that $(\rho,
u,d)(x,T^*))$ satisfies (1.11) and (\ref{1.12}), where
$g(x)=\rho^{1/2}\dot u(x, T^*)$. Thus, Proposition \ref{pro5.1}
implies that there exists some $T^{**}>T^*$ such that (\ref{3.7})
holds for $T=T^{**}$, which contradicts (\ref{5.4}). Hence,
(\ref{5.6}) holds. Proposition \ref{pro5.1} and Lemmas
\ref{lem4.1}--\ref{lem4.6} thus show that $(\rho, u,d)$ is in fact a
unique classical solution on $\R^3\times(0, T]$ for any
$0<T<\infty$.

In order to complete the proof of the existence part, it remains to
prove that $(\rho,u,d)$ is continuous in $t$, especially, to prove
that
\begin{equation}
\left(\rho-\tilde\rho,P(\rho)-P(\tilde\rho)\right)\in C([0,T];
D^2\cap D^{2,q}),\quad q\in(3,6),\label{5.7}
\end{equation}
since by virtue of Lemmas \ref{lem4.1}--\ref{lem4.5} one easily deduces from (\ref{1.1})--(\ref{1.3}) that
\begin{equation}
\begin{cases}
(\rho-\tilde\rho,P(\rho)-P(\tilde\rho))\in C([0,T];H^1\cap
W^{1,q}),\\[2mm]
(\rho-\tilde\rho,P(\rho)-P(\tilde\rho))\in C([0,T];H^2\cap
W^{2,q}-{\rm weakly}),\\[2mm]
(u,\nabla d)\in C([0,T];D^1\cap D^2).
\end{cases}\label{5.8}
\end{equation}

To prove (\ref{5.7}), we denote by $D_{ij}\triangleq\partial_{ij}^2$
with $i,j=1,2,3$. Then it follows from (\ref{1.1}) that
$$
\pa_t D_{ij}\rho+\divg(uD_{ij}\rho)=-\divg(\rho
D_{ij}u)-\divg(\pa_i\rho\cdot\pa_ju+\pa_j\rho\cdot\pa_iu)
$$
holds in $\mathcal{D}'(\R^3\times(0,T))$. Now, let $j_\nu(x)$ is the
standard mollifying kernel with width $\nu$ and set
$\rho^\nu\triangleq\rho\ast j_\nu$. Then, we infer from the above
equation that
\begin{equation}
\pa_t D_{ij}\rho^\nu+\divg(uD_{ij}\rho^\nu)=-\divg(\rho D_{ij}u)\ast
j_\nu-\divg(\pa_i\rho\cdot\pa_ju+\pa_j\rho\cdot\pa_iu)\ast j_\nu
+R_\nu,\label{5.9}
\end{equation}
where $R_\nu\triangleq\divg(uD_{ij}\rho^\nu)-\divg(uD_{ij}\rho)\ast
j_\nu$ satisfies (cf. \cite[Lemma 2.3]{Li1996})
\begin{equation}
\int_0^T\|R_\nu\|_{L^2\cap L^q}^{p_0}dt\leq
C\int_0^T\|u\|_{W^{1,\infty}}^{p_0}\|D_{ij}\rho\|_{L^2\cap
L^q}^{p_0}dt\leq C\label{5.10}
\end{equation}
due to (\ref{4.29}), where $p_0>1$ being the same one as in
(\ref{4.30}).

Multiplying (\ref{5.9}) by $q|D_{ij}\rho^\nu|^{q-2}D_{ij}\rho^\nu$
and integrating by parts over $\R^3$, we obtain
\begin{eqnarray*}
\frac{d}{dt}
\|D_{ij}\rho^\nu\|_{L^q}^q&=&-(q-1)\int|D_{ij}\rho^\nu|^q\divg u
dx-q\int\divg(\rho D_{ij}u)\ast j_\nu
|D_{ij}\rho^\nu|^{q-2}D_{ij}\rho^\nu dx\nonumber\\
&&-q\int\left(\divg(\pa_i\rho\cdot\pa_ju+\pa_j\rho\cdot\pa_iu)\ast
j_\nu\right) |D_{ij}\rho^\nu|^{q-2}D_{ij}\rho^\nu
dx\nonumber\\
&&+q\int R_\nu |D_{ij}\rho^\nu|^{q-2}D_{ij}\rho^\nu dx,
\end{eqnarray*}
which, combining with  (\ref{4.29}) and (\ref{5.10}), yields
\begin{eqnarray*}
&&\sup_{0\leq t\leq T} \|\nabla^2\rho^\nu\|_{L^q}+
\int_0^T\left|\frac{d}{dt}
\|\nabla^2\rho^\nu\|_{L^q}^q\right|^{p_0}dt\\
&&\quad\leq C+C\int_0^T\left(\|\nabla
u\|_{W^{2,q}}^{p_0}+\|R_\nu\|_{L^2\cap L^q}^{p_0}\right)dt\\
&&\quad\leq C.
\end{eqnarray*}
This, together with Ascoli-Arzela theorem, gives
$$
\|\nabla^2\rho^\nu(\cdot,t)\|_{L^q}\to
\|\nabla^2\rho(\cdot,t)\|_{L^q}\quad {\rm in}\quad C([0,T])\quad{\rm
as}\quad \nu\to0.
$$
In particular, we have
\begin{equation}
\|\nabla^2\rho(\cdot,t)\|_{L^q}\in C([0,T]).\label{5.11}
\end{equation}

Similarly, one also has
$$
\|\nabla^2\rho(\cdot,t)\|_{L^2}\in C([0,T]),
$$
which, combined with (\ref{5.11}) and (\ref{5.8})$_2$, implies
\begin{equation}
\nabla^2\rho \in C([0,T];L^2\cap L^q).\label{5.12}
\end{equation}
In the same way, we can also prove that $ \nabla^2P(\rho) \in
C([0,T];L^2\cap L^q)$. This, together with (\ref{5.12}), finishes
the proof of (\ref{5.7}), and thus, the existence of a classical
solution $(\rho,u,d)$ of (\ref{1.1})--(\ref{1.7}) has been proved.

Next we prove the large-time behavior of $(\rho,u,d)$. This can be
done as the ones in \cite{HLX2010,HuLi}, however, for completeness
we sketch the proof here. Multiplying (\ref{3.36}) by
$4(P(\rho)-P(\tilde\rho))^3$ and integrating it over $\R^3$, we get
that
$$
\frac{d}{dt}\|P(\rho)-P(\tilde\rho)\|_{L^4}^4=\int\left(|P(\rho)-P(\tilde\rho)|^4\divg
u-3\gamma P(\rho)(P(\rho)-P(\tilde\rho))^3\divg u\right)dx,
$$
which, together with (\ref{3.95}) and (\ref{3.96}), shows
$$
\int_1^\infty\left|\frac{d}{dt}\|P(\rho)-P(\tilde\rho)\|_{L^4}^4\right|dt
\leq C\int_1^\infty\left(\|\nabla
u\|_{L^4}^4+\|P(\rho)-P(\tilde\rho)\|_{L^4}^4 \right)dt\leq C.
$$
As a result, we have
$$
\|P(\rho)-P(\tilde\rho)\|_{L^4}\to0\quad{\rm as}\quad t\to\infty.
$$
This, together with (\ref{3.11}) and the uniform upper bound of
density, shows that
\begin{equation}
\lim\limits_{t\rightarrow \infty}\|\rho-\tilde\rho \|_{L^p}
=0\label{5.13}
\end{equation}
holds for any $p$ as in \eqref{a116}.

To study the large-time behavior of the velocity, we set
$$
M(t)\triangleq\frac{\mu}{2}\|\nabla
u\|_{L^2}^2+\frac{\mu+\lambda}{2}\|\divg u\|_{L^2}^2.
$$
Then, multiplying (\ref{1.2}) by $\dot u$ in $L^2$ and integrating
by parts over $\R^3$, we obtain
\begin{equation}
-\int\left(\mu\Delta u +(\mu+\lambda)\nabla\divg u\right)\cdot \dot
udx= \int\left((P(\rho)-P(\tilde\rho))\divg\dot u+ M(d):\nabla \dot
u-\rho|\dot u|^2\right)dx.\label{5.14}
\end{equation}
Recalling the definition of `` $\dot{}$ '', we deduce after
integrating by parts that
\begin{eqnarray*}
-\int \mu\Delta u\cdot\dot u dx&=&\frac{\mu}{2}\frac{d}{dt}\|\nabla
u\|_{L^2}^2+\mu\int\partial_k
u^j\partial_k\left(u_i\partial_iu^j\right)dx\\
&=&\frac{\mu}{2}\frac{d}{dt}\|\nabla
u\|_{L^2}^2+\mu\int\left(\partial_k u^j\partial_k
u_i\partial_iu^j-\frac{1}{2}|\nabla u|^2(\divg u)\right)dx,
\end{eqnarray*}
and similarly,
\begin{eqnarray*}
-\int (\mu+\lambda)\nabla\divg u\cdot\dot u
dx&=&\frac{\mu+\lambda}{2}\frac{d}{dt}\|\divg u\|_{L^2}^2+O(\|\nabla
u\|_{L^3}^3),
\end{eqnarray*}
which, inserted into (\ref{5.14}), yields
\begin{eqnarray*}
|M'(t)|&\leq& C \left(\|\nabla u\|_{L^3}^3+\| \rho^{1/2}\dot
u\|_{L^2}^2+\left(\|P(\rho)-P(\tilde\rho)\|_{L^2}+\|\nabla
d\|_{L^3}\|\nabla d\|_{L^6}\right)\|\nabla\dot
u\|_{L^2}\right)\nonumber\\
&\leq&C\left(\|\nabla u\|_{L^2}\|\nabla u\|_{L^4}^2+\|
\rho^{1/2}\dot u\|_{L^2}^2+  \|\nabla\dot u\|_{L^2}+\|\nabla^2
d\|_{L^2}^2 \|\nabla\dot u\|_{L^2}\right).
\end{eqnarray*}
Thus, by (\ref{3.89}) and (\ref{3.96}) we see that
\begin{eqnarray*}
\int_1^\infty|M'(t)|^2dt&\leq& C\int_1^\infty\left(\|\rho^{1/2}\dot
u\|_{L^2}^4+ \|\nabla u\|_{L^4}^{4}+\|\nabla\dot
u\|_{L^2}^2\right)dt\nonumber\\
&\leq& C\left(\sup_{t\geq1}\|\rho^{1/2}\dot
u\|_{L^2}^2\right)\int_1^\infty \|\rho^{1/2}\dot u\|_{L^2}^2dt\leq
C,
\end{eqnarray*}
from which and (\ref{3.89}) it follows that
\begin{equation}
\|\nabla u(t)\|_{L^2}\to0\quad{\rm as}\quad t\to\infty.\label{5.15}
\end{equation}
As a result, we also have
\begin{equation}
\int\rho^{1/2}|u|^4dx\leq \left(\int\rho
|u|^2dx\right)^{1/2}\|u\|_{L^6}\leq C\|\nabla u\|_{L^2}\to0\quad{\rm
as}\quad t\to\infty.\label{5.16}
\end{equation}
Moreover, using (\ref{3.11}), (\ref{3.15}) and (\ref{3.89}), we
infer from (\ref{2.6}) that
$$
\sup_{t\geq1}\|\nabla u(t)\|_{L^6}\leq C,
$$
which, together with (\ref{5.15}) and the interpolation inequality,
leads to
\begin{equation}
\|\nabla u(t)\|_{L^r}\to 0\quad{\rm as}\quad
t\to\infty,\quad\forall\; r\in[2,6).\label{5.17}
\end{equation}

Finally, applying $\nabla$ to both sides of (\ref{1.3}) and taking
the $L^2$-inner product, by (\ref{3.11}) and (\ref{3.89}) we deduce
after integrating by parts over $\R^3\times(1,\infty)$ that
\begin{eqnarray*}
\int_1^\infty\left|\frac{d}{dt}\|\nabla^2d\|_{L^2}^2\right|dt&\leq&
C\int_1^\infty\left(\|\nabla
d_t\|_{L^2}^2+\|\nabla^3d\|_{L^2}^2+\|\nabla d\|_{L^6}^6+\|\nabla
d\|_{L^3}^2\|\nabla^2d\|_{L^6}^2\right)dt\\
&&+C\int_1^\infty\left(\|\nabla u\|_{L^2}^2\|\nabla
d\|_{L^\infty}^2+\|u\|_{L^6}^2\|\nabla^2d\|_{L^2}\|\nabla^2d\|_{L^6}\right)\\
&\leq& C\int_1^\infty\left(\|\nabla
d_t\|_{L^2}^2+\|\nabla^3d\|_{L^2}^2+\|\nabla^2 d\|_{L^2}^2+\|\nabla
u\|_{L^2}^2\right)dt\\
&\leq&C,
\end{eqnarray*}
which, together with (\ref{3.11}), gives
$$
\|\nabla^2 d(t)\|_{L^2}\to 0\quad{\rm as}\quad t\to\infty.
$$
This, combining with  \eqref{3.11} and  \eqref{3.89}, yields that for any $k\in(2,6)$,
\begin{equation}
\|\nabla d(t)\|_{W^{1,k}}\to 0\quad{\rm as}\quad t\to\infty.\label{5.18}
\end{equation}
Combining (\ref{5.13}) with (\ref{5.15})--(\ref{5.18}) immediately
proves (\ref{1.16}). The proof of Theorem \ref{thm1.1} is therefore
completed.\hfill$\square$

\end{document}